\documentclass[10pt, a4paper]{amsart}
\usepackage{amssymb,amsmath,amscd,epsfig}

\newtheorem{theorem}{Theorem}
\theoremstyle{plain}

\newtheorem{corollary}{Corollary}

\newtheorem{definition}{Definition}
\newtheorem{example}{Example}

\newtheorem{lemma}{Lemma}

\newtheorem{proposition}{Proposition}
\newtheorem{remark}{Remark}

\numberwithin{equation}{section}

\begin{document}
\title[On representations in quaternion
algebras]{On representations of 2-bridge knot groups in quaternion
algebras}
\author[H.Hilden]{Hugh M. Hilden}
\address[H.Hilden]{ Departament of Mathematics, University of Hawaii,
Honolulu, HI 96822, USA}
\author[M.T.Lozano]{ Mar\'{\i}a Teresa Lozano*}
\address[M.T.Lozano]{ IUMA, Departamento de Matem\'{a}ticas, Universidad de
Zaragoza, Zaragoza 50009, Spain }
\thanks{*This research was supported by grant MTM2007-67908-C02-01}
\author[J.M.Montesinos]{ Jos\'{e} Mar\'{\i}a Montesinos-Amilibia**}
\address[J.M.Montesinos]{Departamento de Geometr\'{\i}a y Topolog\'{\i}a,
Universidad Complutense, Madrid 28040, Spain}
\thanks{**This research was supported by grant MTM2006-00825}
\date{May, 2008}
\dedicatory{}

\begin{abstract}
Representations of two bridge knot groups in the isometry group of some complete Riemannian 3-manifolds
 as $E^{3}$ (Euclidean 3-space), $H^{3}$ (hyperbolic 3-space) and $
E^{2,1}$ (Minkowski 3-space), using quaternion algebra theory, are studied.
We study the different representations of a 2-generator group in which the generators are send to conjugate elements, by
analyzing the points of an algebraic variety, that we call the \emph{variety of affine
c-representations of }$G$. Each point in this variety correspond to a
representation in the unit group of a quaternion algebra and their affine
deformations.
\end{abstract}

\maketitle

\section{Introduction}

The representations of knot groups in special groups have been used to
define invariants of knots, as the Alexander polynomial (\cite{Rham1967}),
A-polynomial (\cite{CL1996}), peripheral polynomials (\cite{HLM2005}%
), etc. Moreover the representations of a knot group $G(K)$ in the group of
isometries of a geometric manifold contain the holonomies of the orbifold and cone-manifold
structures in $S^{3}$ (or in a manifold obtained by Dehn-surgery in $K$) with
the knot $K$ as singular locus. We are interested in representations of
knot-groups in the isometry group of some complete Riemannian 3-manifolds.
In this paper we focus our attention on representations of 2-generator
groups mapping the generators to conjugate elements, because this case has
immediate application to two-bridge knot groups and to a special case of representations of two bridge link groups.

It turns out that the quaternion algebra theory is useful to study in a
comprehensive way the groups of isometries of some 3-dimensional Riemannian
3-manifolds as $E^{3}$ (Euclidean 3-space), $H^{3}$ (hyperbolic 3-space) and $%
E^{2,1}$ (Minkowski 3-space). The group of orientation preserving isometries
of hyperbolic 3-space $H^{3}$ is $PSL(2,\mathbb{C})\cong SL(2,\mathbb{C}%
)/\{\pm I\},$ where $SL(2,\mathbb{C})$ is the group of unit quaternions in
the quaternion algebra $M(2,\mathbb{C})=\left( \frac{-1,1}{\mathbb{C}}%
\right)$. We explain in this paper how the group of orientation preserving isometries of the Euclidean 3-space $E^{3}$ is associated with the Hamilton
quaternions $H=\left( \frac{-1,-1}{\mathbb{R}}\right) $, and the group of isometries
preserving orientation of the Minkowski 3-space $E^{1,2}$ is related with
the quaternion algebra $M(2,R)=\left( \frac{-1,1}{\mathbb{R}}\right) ,$ in such a way
that we can study the different representations of a 2-generator group by
analyzing the points of an algebraic variety, that we call the \emph{variety of affine
c-representations of }$G$. The points in this variety correspond to
representations in the unit group of a quaternion algebra and their affine
deformations. These representations are important in order to relate knot groups with affine crystallographic groups and Lorentz structures. (For the definition and relevance of these concepts see, for instance, \cite{FG1983}, \cite{CDGM2003}, \cite{G1985}.)

In Section \ref{Sec1} we review some concepts and results on
Quaternion algebra theory. (See \cite{Lam} as general reference). In
Section \ref{Sec2} we characterize a pair of unit quaternions up to
conjugation in different quaternion algebras.
In Section
\ref{srepresentationinU}, Theorem \ref{Teorema uno} describes a new
algorithm to obtain the ideal defining the algebraic \emph{variety
}$V(\mathcal{I}_{G}^{c})$ \emph{of c-representations of}$G$ \emph{in }$SL(2,\mathbb{C})$ giving explicitly the c-representation
associated to each point in the variety $V(\mathcal{I} _{G}^{c})$
and it corresponding minimal quaternion algebra. Theorem
\ref{Teorema dos} gives the complete classification of
c-representations of $G$ in $S^{3}$ and $SL(2,\mathbb{R})$. We apply
Theorems \ref{Teorema uno} and \ref{Teorema dos} to the group of the
Trefoil knot (a no hyperbolic knot) and to the group of the Figure
Eight knot (a hyperbolic knot). Finally in Section \ref{Sec4} we
characterize a pair of conjugate affine isometries and we obtain the
ideal defining the \emph{variety }$V_{a}(\mathcal{I}_{aG}^{c})$
\emph{of affine c-representations of }$G$. As an example we obtain
the ideal $\mathcal{I}_{aG}^{c}$ for the Trefoil knot and the Figure
Eight knot.

\section{Quaternion algebras\label{Sec1}}
\subsection{Definitions}

\begin{definition}
Let $k$ be a number field and let be $\mu ,\nu \in k^{\ast }$, where $%
k^{\ast }$ is the group of invertible elements in $k$. The quaternion
algebra $H=\left( \frac{\mu ,\nu }{k}\right)$ is the $k$-algebra on two
generators $i,j$ with the defining relations:
\begin{equation*}
i^{2}=\mu ,\qquad j^{2}=\nu \qquad \text{and}\qquad ij=-ji.
\end{equation*}
\end{definition}

An easy consequence of the definition is that \ $H$ is a four dimensional
vector space over $k$, with basis $\left\{ 1,i,j,ij\right\} $. The square of
$ij$ is also an invertible element in $k$:%
\begin{equation*}
(ij)(ij)=(ij)(-ji)=-ijji=-\mu \nu \in k^{\ast }.
\end{equation*}

\begin{proposition}
\label{pM2k} For any $k$, the quaternion algebra $\left( \frac{-1,1}{k}%
\right)$ is isomorphic to the algebra $M(2,k)$ of $2\times 2$ matrices
over $k:$
\begin{equation*}
\begin{array}{ccc}
\left( \frac{-1,1}{k}\right) =\langle 1,i,j,ij\rangle  & \longrightarrow  &
M(2,k) \\
i & \rightarrow  & I=\left(
\begin{array}{cc}
0 & 1 \\
-1 & 0%
\end{array}%
\right)  \\
j & \rightarrow  & J=\left(
\begin{array}{cc}
0 & 1 \\
1 & 0%
\end{array}%
\right)
\end{array}%
\end{equation*}\qed
\end{proposition}

\begin{proposition}
\label{pM2C} For the complex field $\mathbb{C}$, there is only one
quaternion algebra, up to isomorphism, the algebra $M(2,\mathbb{C})$. In
fact, given two nonzero complex numbers, $\mu ,\nu $, there exists $x,y\in
\mathbb{C}$, such that $\mu =x^{2}$ and $\nu =y^{2}$. Then the map%
\begin{equation*}
\begin{array}{ccc}
\left( \frac{\mu ,\nu }{\mathbb{C}}\right) =\langle 1,i,j,ij\rangle
& \longrightarrow & M(2,\mathbb{C)=}\left(
\frac{-1,1}{\mathbb{C}}\right)
=\langle 1,I,J,IJ\rangle \\
i & \rightarrow & \sqrt{-1} xI \\
j & \rightarrow & yJ
\end{array}
\end{equation*}
is an algebra isomorphism.\qed
\end{proposition}

Easy consequences of the above result are the following:

\begin{proposition}
\label{pabk} Any quaternion algebra $H=\left( \frac{\mu ,\nu }{k}\right) $, $%
k\subset \mathbb{C}$, is a subalgebra of the algebra $M(2,\mathbb{C)=}\left( \frac{-1,1}{%
\mathbb{C}}\right) $.\qed
\end{proposition}

\begin{proposition}
For the real field $\mathbb{R}$, there are only two quaternion algebras, up
to isomorphism, the algebra $M(2,\mathbb{R})=\left( \frac{-1,1}{\mathbb{R}}
\right)$ and the Hamilton quaternions $\mathbb{H}=\left( \frac{-1,-1}{
\mathbb{R}}\right)$.\qed
\end{proposition}

Given a quaternion $A=\alpha +\beta i+\gamma j+\delta ij$, $A\in H=\left(
\frac{\mu ,\nu }{k}\right) =\langle 1, i, j, ij\rangle $, we
use the notation $A^{+}=\alpha$, and $A^{-}=\beta i+\gamma j+\delta ij$.
Then, $A=A^{+}+A^{-}$. Compare \cite{BH1995}.

\begin{definition}
The \emph{conjugate} of $A$, is by definition, the quaternion
\begin{equation*}
\overline{A}:=A^{+}-A^{-}=\alpha -\beta i-\gamma j-\delta ij.
\end{equation*}

 The \emph{norm} of $A$ is $N(A):=A\overline{A}=\overline{A}A$.
Observe that
\begin{equation*}
N(A)=(\alpha +\beta i+\gamma j+\delta ij)(\alpha -\beta i-\gamma j-\delta
ij)=\alpha ^{2}-\beta ^{2}\mu -\gamma ^{2}\nu +\delta ^{2}\mu \nu \in k.
\end{equation*}

The \emph{trace} of $A$ is $T(A)=A+\overline{A}$. Observe that
\begin{equation*}
T(A)=2A^{+}=2\alpha \in k.
\end{equation*}
\end{definition}

The map

\begin{equation*}
\begin{array}{cccc}
\langle \quad ,\quad \rangle : & H\times H & \longrightarrow & k \\
& (A,B) & \rightarrow & \frac{1}{2}(A\overline{B}+B\overline{A})=\frac{1}{2}%
T(A\overline{B})=(A\overline{B})^{+}%
\end{array}%
\end{equation*}%
is a \emph{symmetric bilinear form} on $H$. The associated quadratic form%
\begin{equation*}
\begin{array}{ccc}
H & \longrightarrow & k \\
A & \rightarrow & \frac{1}{2}T(A\overline{A})=\frac{1}{2}2N(A)=N(A)%
\end{array}%
\end{equation*}%
is the \emph{norm form} $N$ on $H$. We denote by $(H,N)$ the quadratic
structure in $H$. Observe that the associated matrix to the norm form on $H$
in the basis $\left\{ 1,\quad i,\quad j,\quad ij\right\}$, is the matrix $%
\left(
\begin{array}{cccc}
1 & 0 & 0 & 0 \\
0 & -\mu & 0 & 0 \\
0 & 0 & -\nu & 0 \\
0 & 0 & 0 & \mu \nu%
\end{array}%
\right) .$

\subsection{Examples}

We are mostly interested in the following three examples: $M(2,\mathbb{C})$, $\mathbb{H}$ and  $M(2,\mathbb{R})$.

\subsubsection{ $M(2,\mathbb{C})$}
The quaternion algebra $M(2,\mathbb{C)=}\left( \frac{-1,1}{\mathbb{C}}%
\right) $ (See Proposition \ref{pM2k}). The trace is the usual trace of the
matrix, and the norm is the determinant of the matrix. Proposition \ref{pM2C}
shows that it is possible to change the basis to obtain any possible pair $(\mu ,\nu )$ of complex numbers defining the structure of quaternion algebra
$\left( \frac{\mu ,\nu }{\mathbb{C}}\right) $. In particular
\begin{equation*}
M(2,\mathbb{C)=}\left( \frac{-1,1}{\mathbb{C}}\right) =\left( \frac{1,1}{\mathbb{C}}\right) =\left( \frac{-1,-1}{\mathbb{C}}\right)
\end{equation*}

Next we analyze, using Proposition \ref{pM2C}, the basis and the corresponding matrix
associated to the norm form for this three presentations.
\begin{enumerate}
\item $M(2,\mathbb{C)=}\left( \frac{-1,1}{\mathbb{C}}\right) $.
As in Proposition \ref{pM2k}, we consider  the basis $\left\{ Ident,I,J,IJ\right\} $, where
\begin{equation*}
I=\left(
\begin{array}{cc}
0 & 1 \\
-1 & 0%
\end{array}
\right) ,\quad J=\left(
\begin{array}{cc}
0 & 1 \\
1 & 0%
\end{array}
\right) ,\quad IJ=\left(
\begin{array}{cc}
1 & 0 \\
0 & -1
\end{array}
\right)
\end{equation*}
For $A=\alpha +\beta I+\gamma J+\delta IJ=
\begin{pmatrix}
\alpha +\delta & \beta +\gamma \\
-\beta +\gamma & \alpha -\delta
\end{pmatrix}
\in M(2,\mathbb{C})$, $\alpha ,\beta ,\gamma ,\delta \in \mathbb{C}$,
$N(A)=\alpha ^{2}+\beta ^{2}-\gamma ^{2}-\delta ^{2}$. Therefore the matrix
associated to the quadratic form in the basis $\left\{ Ident,I,
J,IJ\right\}$, is the matrix $\left(
\begin{array}{cccc}
1 & 0 & 0 & 0 \\
0 & 1 & 0 & 0 \\
0 & 0 & -1 & 0 \\
0 & 0 & 0 & -1
\end{array}
\right) $.

\item $M(2,\mathbb{C)=}\left( \frac{1,1}{\mathbb{C}}\right) $. Using Proposition \ref{pM2C}, $x=1$, $y=1$, we consider the basis $
\left\{ Ident,I_{0},J_{0},I_{0}J_{0}\right\} $ , where
\begin{gather*}
I_{0}=\sqrt{-1}I=\left(
\begin{array}{cc}
0 & \sqrt{-1} \\
-\sqrt{-1} & 0
\end{array}
\right) ,\quad J_{0}=J=\left(
\begin{array}{cc}
0 & 1 \\
1 & 0
\end{array}
\right) ,\\
 I_{0}J_{0}=\left(
\begin{array}{cc}
\sqrt{-1} & 0 \\
0 & -\sqrt{-1}
\end{array}
\right)\qquad
\end{gather*}
For $A=\alpha +\beta I_{0}+\gamma J_{0}+\delta I_{0}J_{0}=
\begin{pmatrix}
\alpha +\delta \sqrt{-1} & \gamma +\beta \sqrt{-1} \\
\gamma -\beta \sqrt{-1} & \alpha -\delta \sqrt{-1}
\end{pmatrix}%
\in M(2,\mathbb{C})$, where $\alpha ,\beta ,\gamma ,\delta \in \mathbb{C}$, $%
N(A)=\alpha ^{2}-\beta ^{2}-\gamma ^{2}+\delta ^{2}$. Therefore the matrix
associated to the quadratic form in the basis $\left\{ Ident,
I_{0}, J_{0}, I_{0}J_{0}\right\}$, is the matrix $\left(
\begin{array}{cccc}
1 & 0 & 0 & 0 \\
0 & -1 & 0 & 0 \\
0 & 0 & -1 & 0 \\
0 & 0 & 0 & 1
\end{array}
\right) $.

\item $M(2,\mathbb{C)=}\left( \frac{-1,-1}{\mathbb{C}}\right) $. Using Proposition \ref{pM2C}, $x=\sqrt{-1}$, $y=\sqrt{-1}$, we consider  the basis
$\left\{ Ident,I_{1},J_{1},I_{1}J_{1}\right\} $ , where
\begin{gather*}
I_{1}=\sqrt{-1} \sqrt{-1} I=-I=\left(
\begin{array}{cc}
0 & -1 \\
1 & 0
\end{array}
\right) ,\quad J_{1}=\sqrt{-1} J=\left(
\begin{array}{cc}
0 & \sqrt{-1} \\
\sqrt{-1} & 0
\end{array}
\right) ,\\
 I_{1}J_{1}=\left(
\begin{array}{cc}
-\sqrt{-1} & 0 \\
0 & \sqrt{-1}
\end{array}
\right)
\end{gather*}
For $A=\alpha +\beta I_{1}+\gamma J_{1}+\delta I_{1}J_{1}=
\begin{pmatrix}
\alpha -\delta \sqrt{-1} & -\beta +\gamma \sqrt{-1} \\
\beta +\gamma \sqrt{-1} & \alpha +\delta \sqrt{-1}
\end{pmatrix}
\in M(2,\mathbb{C})$, where $\alpha ,\beta ,\gamma ,\delta \in \mathbb{C}$, $
N(A)=\alpha ^{2}+\beta ^{2}+\gamma ^{2}+\delta ^{2}$. Therefore the matrix
associated to the quadratic form in the basis $\left\{ Ident,
I_{1}, J_{1}, I_{1}J_{1}\right\}$, is the matrix $\left(
\begin{array}{cccc}
1 & 0 & 0 & 0 \\
0 & 1 & 0 & 0 \\
0 & 0 & 1 & 0 \\
0 & 0 & 0 & 1
\end{array}
\right) $.
\end{enumerate}

\subsubsection{$\mathbb{H}$}
The quaternion algebra $\mathbb{H}=\left(
\frac{-1,-1}{\mathbb{R}}\right) $is the algebra of Hamilton quaternions. It is isomorphic to an $\mathbb{R}$ -subalgebra of the
$\mathbb{C}$-algebra $M(2,\mathbb{C})\mathbb{=}\left(
\frac{-1,-1}{\mathbb{C}}\right) $.
\begin{equation}
\begin{array}{ccc}
\Phi :\mathbb{H}=\left( \frac{-1,-1}{\mathbb{R}}\right) =\langle
1,i,j,ij\rangle & \longrightarrow & M(2,\mathbb{C}) \\
i & \rightarrow & I_{1}=\left(
\begin{array}{cc}
0 & -1 \\
1 & 0
\end{array}
\right) \\
j & \rightarrow & J_{1}=\left(
\begin{array}{cc}
0 & \sqrt{-1} \\
\sqrt{-1} & 0
\end{array}
\right)
\end{array}
\label{inclusionhamilton}
\end{equation}

For $A=\alpha +\beta i+\gamma j+\delta ij\in \mathbb{H}$, $N(A)=\alpha
^{2}+\beta ^{2}+\gamma ^{2}+\delta ^{2}$, where $\alpha ,\beta ,\gamma
,\delta \in \mathbb{R}$. Therefore the matrix associated to the quadratic
form in the basis $\left\{ 1,i,j,ij\right\}$, is the Identity matrix $
\left(
\begin{array}{cccc}
1 & 0 & 0 & 0 \\
0 & 1 & 0 & 0 \\
0 & 0 & 1 & 0 \\
0 & 0 & 0 & 1
\end{array}
\right) $. Note that $\Phi (A)=\Phi (\alpha +\beta i+\gamma j+\delta ij)=
\begin{pmatrix}
\alpha -\delta \sqrt{-1} & -\beta +\gamma \sqrt{-1} \\
\beta +\gamma \sqrt{-1} & \alpha +\delta \sqrt{-1}
\end{pmatrix}
\in M(2,\mathbb{C})$.

\subsubsection{$M(2,\mathbb{R})$}
The quaternion algebra $M(2,\mathbb{R})=\left(
\frac{-1,1}{\mathbb{R}} \right) $. See Proposition \ref{pM2k}. It is
isomorphic to an $\mathbb{R}$ -subalgebra of the $\mathbb{C}$-algebra
$M(2,\mathbb{C})\mathbb{=}\left( \frac{-1,1}{\mathbb{C}}\right) $.
\begin{equation}
\begin{array}{ccc}
\Psi :M(2,\mathbb{R})=\left( \frac{-1,1}{\mathbb{R}}\right) =\langle
1,I,J,IJ\rangle & \longrightarrow & M(2,\mathbb{C}) \\
I & \rightarrow & I=\left(
\begin{array}{cc}
0 & 1 \\
-1 & 0%
\end{array}
\right) \\
J & \rightarrow & J=\left(
\begin{array}{cc}
0 & 1 \\
1 & 0%
\end{array}
\right)
\end{array}
\label{inclusionm2r}
\end{equation}

For $A=\alpha +\beta I+\gamma J+\delta IJ\in M(2,\mathbb{R})$, $N(A)=\alpha
^{2}+\beta ^{2}-\gamma ^{2}-\delta ^{2}$, where $\alpha ,\beta ,\gamma
,\delta \in \mathbb{R}$. Therefore the matrix associated to the quadratic
form in the basis $\left\{ 1,I,J,IJ\right\}$, is the matrix $\left(
\begin{array}{cccc}
1 & 0 & 0 & 0 \\
0 & 1 & 0 & 0 \\
0 & 0 & -1 & 0 \\
0 & 0 & 0 & -1%
\end{array}%
\right) $. Note that $\Psi (A)=\Psi (\alpha +\beta I+\gamma J+\delta IJ)=%
\begin{pmatrix}
\alpha +\delta & \beta +\gamma \\
-\beta +\gamma & \alpha -\delta%
\end{pmatrix}%
\in M(2,\mathbb{C})$.

\subsection{Pure and unit quaternions}

There are two important subsets in the quaternion algebra $H=\left( \frac{
\mu ,\nu }{k}\right) $. The pure quaternions $H_{0}$ (a 3-dimensional vector
space), and the unit quaternions $U_{1}$ (a multiplicative group).

The \emph{pure quaternions} $H_{0}=\left\{ A\in H:A^{+}=0\right\} $ form a
three dimensional vector space over $k$ generated by $\left\{ i,j,ij\right\} $.
 The symmetric bilinear form $\langle \quad ,\quad \rangle $ restricts to $H_{0}$ defining the quadratic space $(H_{0},N)$. Because $\overline{A^{-}}=-A^{-}$, we have that
\begin{equation*}
\left\langle A^{-},B^{-}\right\rangle =\frac{1}{2}\left( A^{-}\overline{B^{-}}+B^{-}\overline{A^{-}}\right) =\frac{1}{2}(-A^{-}B^{-}-B^{-}A^{-})=-(A^{-}B^{-})^{+}
\end{equation*}

Therefore two elements $A^{-},B^{-}\in H_{0}$ are orthogonal for the
symmetric bilinear form $\langle \quad ,\quad \rangle $ , ( $\left\langle
A^{-},B^{-}\right\rangle =0$) if and only if $A^{-},B^{-}$ anticommute. The
elements $i,j,ij$ anticommute. Then $\left\{ i,j,ij\right\} $ is a
orthogonal basis in $H_{0}$. In this basis the associated matrix to the
restricted symmetric bilinear form $\langle \quad ,\quad \rangle _{|H_{0}}$
is $\left(
\begin{array}{ccc}
-\mu & 0 & 0 \\
0 & -\nu & 0 \\
0 & 0 & \mu \nu
\end{array}
\right) $.

Denote by $O(H_{0},N)$ the \emph{orthogonal group} for the quadratic space $(H_{0},N)$:
\begin{equation*}
O(H_{0},N)=\left\{ \varphi :H_{0}\overset{\cong }{\longrightarrow }
H_{0}:N(\varphi (A^{-}))=N(A^{-}),A^{-}\in H_{0}\right\}
\end{equation*}

The isomorphism $(H_{0},N)\cong \left( k^{3},\left(
\begin{array}{ccc}
-\mu & 0 & 0 \\
0 & -\nu & 0 \\
0 & 0 & \mu \nu
\end{array}
\right) \right) $, defined by the basis $\left\{ i,j,ij\right\} $, induced
the isomorphism
\begin{equation*}
O(H_{0},N)\cong \left\{ M\in GL(3,k):M^{t}\left(
\begin{array}{ccc}
-\mu & 0 & 0 \\
0 & -\nu & 0 \\
0 & 0 & \mu \nu
\end{array}
\right) M=\left(
\begin{array}{ccc}
-\mu & 0 & 0 \\
0 & -\nu & 0 \\
0 & 0 & \mu \nu
\end{array}
\right) \right\}
\end{equation*}

The \emph{special orthogonal group} is defined by
\begin{equation*}
SO(H_{0},N)=\left\{ M\in O(H_{0},N):\det (M)=1\right\} .
\end{equation*}

Consider the group $U$ of invertible elements in $H$. There exits a short exact sequence of groups
\begin{equation*}
1\longrightarrow k^{\ast }\longrightarrow U\overset{c}{\longrightarrow }
SO(H_{0},N)\longrightarrow 1
\end{equation*}
such that $c(A)$ acts in $H_{0}$ by conjugation : $c(A)(B^{-})=AB^{-}A^{-1}$. See \cite[p. 63]{Lam}.

The \emph{unit quaternions} $U_{1}$ are the elements in the group $U$ with
norm 1.

We will also consider the \emph{group }$U_{\pm 1}$\emph{\ of quaternions
with norm }$\pm 1.$ It is a subgroup of $U$ and either it coincides with $U_{1}$ or has $U_{1}$ as a subgroup of index 2.
\begin{equation*}
U_{1}\vartriangleleft U_{\pm 1}<U
\end{equation*}

\subsubsection{The matrices and action of $c(A)$ as linear map in $H_{0}$}
Let $H_{0}$ be the vector space of pure quaternions in a quaternion algebra $H=\left( \frac{\mu ,\nu }{k}\right) $. Consider a coordinate system $\left\{
X,Y,Z\right\} $ in $H_{0}$ in the basis $\left\{ -ij,j,i\right\} $, such
that a pure quaternion is represented by $Zi+Yj-Xij$. For a unit quaternion $A=\alpha +\beta i+\gamma j+\delta ij,$ $c(A)$ acts in $H_{0}$ as a linear map ,
therefore there exists an associated $3\times 3$  matrix $m(\mu ,\nu; \alpha ,\beta ,\gamma ,\delta )$ such that
\begin{equation*}
\begin{pmatrix}
X^{\prime } \\
Y^{\prime } \\
Z^{\prime }
\end{pmatrix}
=m(\mu ,\nu; \alpha ,\beta ,\gamma ,\delta )
\begin{pmatrix}
X \\
Y \\
Z
\end{pmatrix}
\end{equation*}
where $(Z^{\prime }i+Y^{\prime }j-X^{\prime }ij)=c(A)\left(
(Zi+Yj-Xij)\right)
=(\alpha +\beta i+\gamma j+\delta ij)(Zi+Yj-Xij)(\alpha -\beta i-\gamma j-\delta ij).$
Then a straightforward computation yields the following explicit matrix
\begin{multline}
m(\mu ,\nu; \alpha,\beta ,\gamma ,\delta )=\\
=
\begin{pmatrix}
\alpha^{2}+\mu \beta^{2}+\nu \gamma^{2}+\mu \nu \delta^{2} & -2\alpha \beta +2\nu
\gamma \delta  & 2\alpha \gamma +2\mu \beta \delta  \\
-2\mu \alpha \beta -2\mu \nu \gamma \delta  & \alpha^{2}+\mu \beta^{2}-\nu
\gamma^{2}-\mu \nu \delta^{2} & -2\mu \beta \gamma -2\mu \alpha \delta  \\
2\nu \alpha \gamma -2\mu \nu \beta \delta  & -2\nu \beta \gamma +2\nu \alpha \delta  &
\alpha^{2}-\mu \beta ^{2}+\nu \gamma^{2}-\mu \nu \delta^{2}
\end{pmatrix}
\label{emtxyz}
\end{multline}

\begin{remark}
We use the coordinate system $\left\{
X,Y,Z\right\} $ in $H_{0}$ in the basis $\left\{ -ij,j,i\right\} $, in order to have the usual pictures in the particular case of the Minkowski space, as we will see in Example \ref{eminkowski}.
\end{remark}

\begin{example}{$M(2,\mathbb{C})$.}
In the case $H=M(2,\mathbb{C})=\left( \frac{-1,-1}{\mathbb{C}}\right) $
there exits a short exact sequence of groups
\begin{equation}
1\longrightarrow \left\{ \pm 1\right\} \longrightarrow U_{1}\overset{c}{%
\longrightarrow }SO(H_{0},N)\longrightarrow 1
\end{equation}
where $U_{1}=\{A\in U:N(A)=1\}=SL(2,\mathbb{C})$. $U_{1}\neq U_{\pm 1}$.
Here $(H_{0},N)$ is the 3-dimensional complex space $\mathbb{C}^{3}$, with
the symmetric bilinear form defined by the identity matrix in the basis $\left\{ Ident,I_{1},J_{1},I_{1}J_{1}\right\} $, because $\mu =\nu =-1$, and $SO(H_{0},N)\cong SO(3,\mathbb{C)}$. The above short exact sequence defines
the 2-fold covering
\begin{equation*}
U_{1}\cong SL(2,\mathbb{C})\overset{c}{\longrightarrow }SO(3,\mathbb{C)}
\end{equation*}
and the isomorphism $PSL(2,C)\cong SO(3,\mathbb{C)}$.
\end{example}

\begin{example}{$\mathbb{H}$.}
In the case $\mathbb{H}=\left( \frac{-1,-1}{\mathbb{R}}\right) $ there exits
a short exact sequence of groups
\begin{equation}
1\longrightarrow \left\{ \pm 1\right\} \longrightarrow U_{1}\overset{c}{
\longrightarrow }SO(H_{0},N)\longrightarrow 1  \label{eSO3}
\end{equation}
where $U_{1}=\left\{ A\in U:N(A)=1\right\} =U_{\pm 1}$ is the group of unit
quaternions. Here $U_{1}\cong S^{3}$, $(H_{0},N)$ is the Euclidean
3-dimensional real space $E^{3}$, because $\mu =\nu =-1$, and $SO(H_{0},N)\cong SO(3,\mathbb{R)}$. The above short exact sequence defines
the 2-fold covering
\begin{equation*}
U_{1}\cong S^{3}\overset{c}{\longrightarrow }SO(3,\mathbb{R)}
\end{equation*}
\end{example}

 Observe that the inclusion $\Phi :\mathbb{H}=\left( \frac{-1,-1}{
\mathbb{R}}\right) \longrightarrow M(2,\mathbb{C})$, maps the unit
quaternions $S^{3}$ isomorphically on $SU(2)\subset SL(2,\mathbb{C})$. Let $
A=\alpha +\beta i+\gamma j+\delta ij\in U_{1}\cong S^{3}$, then $
1=N(A)=\alpha ^{2}+\beta ^{2}+\gamma ^{2}+\delta ^{2}$, and
\begin{equation*}
\text{ }\Phi (A)=\Phi (\alpha +\beta i+\gamma j+\delta ij)=
\begin{pmatrix}
\alpha -\delta \sqrt{-1} & -\beta +\gamma \sqrt{-1} \\
\beta +\gamma \sqrt{-1} & \alpha +\delta \sqrt{-1}
\end{pmatrix}
=
\begin{pmatrix}
z_{1} & -z_{2} \\
\overline{z_{2}} & \overline{z_{1}}
\end{pmatrix}
\end{equation*}
which is an element of $SU(2)=\left\{ M\in M(2,\mathbb{C});M^{-1}=\overline{
M^{t}}\right\} \subset SL(2,\mathbb{C})$. Reciprocally, every element of $
SU(2)$ is of the form $
\begin{pmatrix}
z_{1} & -z_{2} \\
\overline{z_{2}} & \overline{z_{1}}
\end{pmatrix}
$, and therefore is the image by $\Phi $ of an element of $S^{3}$.

Next, we analyze the action of $c(A)$ as a linear map, $A\in U_{1}\cong
S^{3} $, in $(\mathbb{H}_{0},N)\cong E^{3}$.
\begin{equation*}
c(A):E^{3}\longrightarrow E^{3}
\end{equation*}

Consider $A=\alpha +\beta i+\gamma j+\delta ij\in U_{1}$, then $
1=N(A)=\alpha ^{2}+\beta ^{2}+\gamma ^{2}+\delta ^{2}$. Therefore, because $\alpha^{2}\leq 1$, we can
assume that $A^{+}=\alpha =\cos (\frac{\theta }{2})$, and $
N(A^{-})=-A^{-}A^{-}=\beta ^{2}+\gamma ^{2}+\delta ^{2}=1-\alpha ^{2}=1-\cos
^{2}(\frac{\theta }{2})=\sin ^{2}(\frac{\theta }{2})$.

\begin{proposition}\label{proposicion5}
The action \ of $c(A)$ on $H_{0}\cong E^{3}$ is a right rotation with angle $
\theta $ around the oriented axis $A^{-}$. (We assume $0\leq \theta \leq \pi
$)
\end{proposition}

\begin{proof}
Any element $A=A^{+}+A^{-}=\alpha +A^{-}$ in the group $U_{1}$ is conjugate
in $U_{1}$ to one with the same trace $2\alpha $, say $\alpha =\cos (\frac{
\theta }{2})$ and $\gamma =\delta =0$.

Then, it is enough to consider
\begin{equation*}
A=\cos (\frac{\theta }{2})+\sin (\frac{\theta }{2})i
\end{equation*}
The matrix of the action of $c(A)$ as a linear map (see (\ref{emtxyz})) on
the basis $\left\{ -ij, j, i\right\} $ of $E^{3}$ is:
\begin{multline*}
m(-1,-1;\cos (\frac{\theta }{2}),\sin (\frac{\theta }{2}),0,0) =\\
=\begin{pmatrix} \cos ^{2}({\frac{\theta }{2}})-\sin
^{2}({\frac{\theta }{2}})  & -2\cos
({\frac{\theta }{2}})\sin({\frac{\theta }{2}}) & 0 \\
2\cos ({\frac{\theta }{2}})\sin ({\frac{\theta }{2}})& \cos
^{2}({\frac{\theta }{2}})-\sin ^{2}({\frac{\theta }{2}}) & 0 \\
0 & 0 & 1
\end{pmatrix}
=
\begin{pmatrix}
\cos \theta & -\sin \theta & 0 \\
\sin \theta & \cos \theta & 0 \\
0 & 0 & 1
\end{pmatrix}
\end{multline*}
\end{proof}

\begin{example}\label{eminkowski}{$M(2,\mathbb{R})$.}
In $M(2,\mathbb{R})=\left( \frac{-1,1}{\mathbb{R}}\right) $, The
3-dimensional vector space of pure quaternions with the norm form is
isomorphic to the Minkowski space $E^{1,2}$, because the matrix of the
quadratic form $\langle \quad ,\quad \rangle $ in the basis $\left\{
I,J,IJ\right\} $ is
\begin{equation*}
\left(
\begin{array}{ccc}
-\mu  & 0 & 0 \\
0 & -\nu  & 0 \\
0 & 0 & \mu \nu
\end{array}%
\right) =\left(
\begin{array}{ccc}
1 & 0 & 0 \\
0 & -1 & 0 \\
0 & 0 & -1
\end{array}
\right)
\end{equation*}
Choose a component $\mathcal{N}_{+}$ of the complement of  $0$ in the
nullcone or light cone
\begin{equation*}
\mathcal{N}=\left\{ q\in E^{1,2}:\left\langle q,q\right\rangle =0\right\}
\end{equation*}
Denote by $SO^{+}(1,2)$ the index 2 subgroup of $SO(1,2)$ preserving the
component $\mathcal{N}_{+}$, it is also the connected component of the group
$SO(1,2)$ containing the identity matrix.
\begin{equation*}
SO^{+}(1,2)=\left\{
\begin{array}{c}
M\in SL(3,\mathbb{R}):M^{t}
\begin{pmatrix}
1 & 0 & 0 \\
0 & -1 & 0 \\
0 & 0 & -1
\end{pmatrix}
M=%
\begin{pmatrix}
1 & 0 & 0 \\
0 & -1 & 0 \\
0 & 0 & -1%
\end{pmatrix}
, \\
\quad M%
\begin{pmatrix}
1 \\
0 \\
0%
\end{pmatrix}
=%
\begin{pmatrix}
x>0 \\
y \\
z%
\end{pmatrix}
\end{array}
\right\}
\end{equation*}
In this case there also exist a short exact sequence of groups
\begin{equation}
1\longrightarrow \left\{ \pm 1\right\} \longrightarrow U_{1}\overset{c}{
\longrightarrow }SO^{+}(1,2)\longrightarrow 1  \label{eSO12}
\end{equation}
\end{example}

\begin{figure}[ht]
\epsfig{file=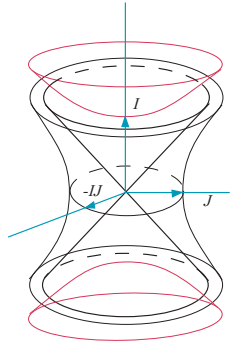,height=4.5cm}
\caption{The Minkowski space
$E^{1,2}$.}
\label{fminkowski}
\end{figure}

Observe that the inclusion $\Psi :M(2,\mathbb{R})=\left( \frac{-1,1}{\mathbb{
R}}\right) \longrightarrow M(2,\mathbb{C})$, maps the unit quaternions $U_{1}$
isomorphically onto the subgroup $SL(2,\mathbb{R})\subset SL(2,\mathbb{C})$. Observe that $U_{1}\neq
U_{\pm 1}$. Let $A=\alpha +\beta I+\gamma J+\delta IJ\in U_{1}$, then $
1=N(A)=\alpha ^{2}+\beta ^{2}-\gamma ^{2}-\delta ^{2}$, and
\begin{equation*}
\text{ }\Psi (A)=\Psi (\alpha +\beta I+\gamma J+\delta IJ)=
\begin{pmatrix}
\alpha +\delta  & \beta +\gamma  \\
-\beta +\gamma  & \alpha -\delta
\end{pmatrix}
\end{equation*}
which is an element of $SL(2,\mathbb{R})$ because it is a real matrix with
determinant equal to 1:
\begin{eqnarray*}
Det
\begin{pmatrix}
\alpha +\delta  & \beta +\gamma  \\
-\beta +\gamma  & \alpha -\delta
\end{pmatrix}
&=&(\alpha +\delta )(\alpha -\delta )-(\beta +\gamma )(-\beta +\gamma ) \\
&=&\alpha ^{2}+\beta ^{2}-\gamma ^{2}-\delta ^{2}=N(A)=1
\end{eqnarray*}
Reciprocally, every element $\begin{pmatrix}
a  & b  \\
c  & d
\end{pmatrix} \in SL(2,\mathbb{R})$ is the image by $\Psi $ of the element $A=\alpha +\beta I+\gamma J+\delta IJ$
 where $\alpha =\frac{a+d}{2}$, $\beta =\frac{b-c}{2}$, $\gamma =\frac{b+c}{2}$, $\delta
 =\frac{a-d}{2}$, and $A\in U_{1}$
 because $\alpha ^{2}+\beta ^{2}-\gamma ^{2}-\delta ^{2}=ad-bc=1$. Note that
considering a coordinate system $\left\{ X,Y,Z\right\} $ in $H_{0}$
in the basis $\left\{ -IJ,J,I\right\} $, the unit pure quaternions
$U_{1}\cap H_{0}$ constitute the two sheeted hyperboloid:
$Z^{2}-Y^{2}-X^{2}=1$. The pure quaternions with norm $-1$ are
the points in the deSitter sphere: $Z^{2}-Y^{2}-X^{2}=-1.$ See
Figure \ref{fminkowski}.

To study the element $c(A)$, $A\in U_{1}$ as a linear map in $E^{1,2}$, we
consider three cases according to the value of $N(A^{-})(>,<,=)0$. Note that
$N(A)=N(A^{+})+N(A^{-})=(A^{+})^{2}+N(A^{-})=1$. Then $
N(A^{-})=1-(A^{+})^{2} $.

\begin{enumerate}
\item[Case 1] $N(A)=1$, $N(A^{-})>0$. The vector $A^{-}$ is inside the light
cone. We say that $A^{-}$ is a \emph{time-like vector}. Then $A=\cos \left( \frac{\theta }{2}\right) +\sin \left( \frac{\theta
}{2}\right) P^{-}$ where $N(P^{-})=1$. Because the kernel of $c$ is $\left\{
\pm 1\right\} $,we assume that the $Z-$coordinate of $P^{-}$ is positive.
Then $P^{-}$ is a point in the upper sheet of the two sheeted hyperboloid
defined by the unit pure quaternions. Up to conjugation (mapping $P^{-}$ to $
I$)
\begin{equation*}
A=\cos (\frac{\theta }{2})+\sin (\frac{\theta }{2})I
\end{equation*}%
The matrix of the action of $c(A)$ as a linear map (see (\ref{emtxyz})) on
the basis $\left\{ -IJ,J,I\right\} $ of $E^{1,2}$ is:
\begin{multline*}
m(-1,1;\cos (\frac{\theta }{2}),\sin (\frac{\theta }{2}),0,0) =\\
=
\begin{pmatrix}
\cos ^{2}(\frac{\theta }{2})-\sin ^{2}(\frac{\theta }{2})& -2\cos
(\frac{\theta }{2})\sin (\frac{\theta }{2})& 0 \\
2\cos (\frac{\theta }{2})\sin (\frac{\theta }{2}) & \cos
^{2}(\frac{\theta}{2})-\sin ^{2}(\frac{\theta}{2}) & 0 \\
0 & 0 & 1
\end{pmatrix}
=
\begin{pmatrix}
\cos \theta  & -\sin \theta  & 0 \\
\sin \theta  & \cos \theta  & 0 \\
0 & 0 & 1%
\end{pmatrix}
\end{multline*}
Then the action of $c(\cos (\frac{\theta }{2})+\sin (\frac{\theta }{2})I)$
is a (Euclidean) positive rotation around the oriented axis $I$ with angle $
\theta $. In the general case, $A=\cos \left( \frac{\theta }{2}\right) +\sin
\left( \frac{\theta }{2}\right) P^{-}$, $c(A)$ acts as a (hyperbolic)
rotation around the oriented axis $P^{-}$ with angle $\theta $, preserving
its orthogonal plane $(P^{-})^{\bot }=\left\{ B^{-}:\left\langle
P^{-},B^{-}\right\rangle =0\right\} $. To understand this rotation consider
in $RP^{2}$ the conic defined by the nullcone as the boundary of the
hyperbolic plane $H^{2}$. See Figure \ref{rothyper}.

\begin{figure}[ht]
\epsfig{file=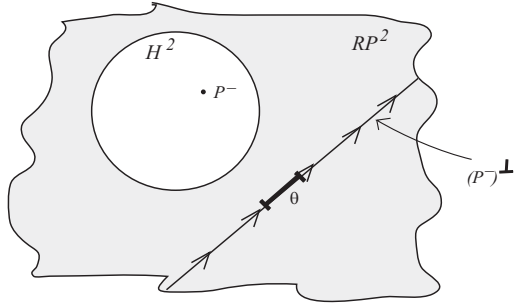,height=4cm}\caption{Action of $A=\cos \left( \frac{\protect\theta }{2}
\right) +\sin \left( \frac{\protect\theta }{2}\right) P^{-}$ in $RP^{2}$.}\label{rothyper}
\end{figure}

\item[Case 2] $N(A)=1$, $N(A^{-})<0$. Then $A=\alpha +A^{-}$, $\alpha >1$.
Thus $A=\cosh \left( \frac{d}{2}\right) +\sinh \left( \frac{d}{2}\right)
P^{-}$ where $N(P^{-})=-1$. The vector $P^{-}$ belongs to the deSitter
sphere, and we assume that its $Y-$coordinate is nonnegative. The vectors $A^{-}$ and $P^{-}$ are \emph{space-like vectors}. Up to
conjugation we have
\begin{equation*}
A=\cosh \left( \frac{d}{2}\right) +\sinh \left( \frac{d}{2}\right) J
\end{equation*}
The matrix of the action of $c(A)$ as a linear map (see (\ref{emtxyz})) on
the basis $\left\{ -IJ,\quad J,\quad I\right\} $ of $E^{1,2}$ is:
\begin{multline*}
m(-1,1;\cosh (\frac{d}{2}),0,\sinh (\frac{d}{2}),0) =\\
=
\begin{pmatrix}
\cosh ^{2}(\frac{d}{2})+\sinh ^{2}(\frac{d}{2})& 0 & 2\cosh
(\frac{d}{2})\sinh (\frac{d}{2}) \\
0 & 1 & 0 \\
2\cosh (\frac{d}{2})\sinh (\frac{d}{2}) & 0 & \cosh ^{2}(\frac{d}{2})+\sinh ^{2}(\frac{d}{2})
\end{pmatrix}
\\=
\begin{pmatrix}
\cosh (d) & 0 & -\sinh (d) \\
0 & 1 & 0 \\
-\sinh (d) & 0 & \cosh (d)
\end{pmatrix}
\end{multline*}
The action is a hyperbolic rotation around the oriented axis
$P^{-}$. The action on the plane $\left\{ -IJ,I\right\} \subset
E^{1,2}$, a plane orthogonal to $P^{-}=J$, is depicted in Figure
\ref{girohyper}.

\begin{figure}[ht]
\epsfig{file=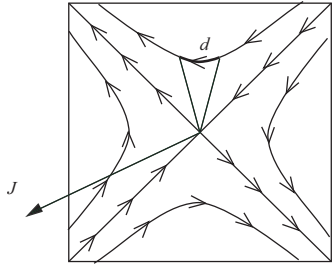,height=3.5cm}\caption{Action of $A=\cosh \left( \frac{d}{2}\right)
+\sinh \left( \frac{d}{2}\right) J$ in $\left\{ -IJ,I\right\} .$}\label{girohyper}
\end{figure}

The action of $A=\cosh
\left( \frac{d}{2}\right) +\sinh \left( \frac{d}{2}\right) J$ in $RP^{2}$ is
depicted in Figure \ref{fcase2}.

\begin{figure}[ht]
\epsfig{file=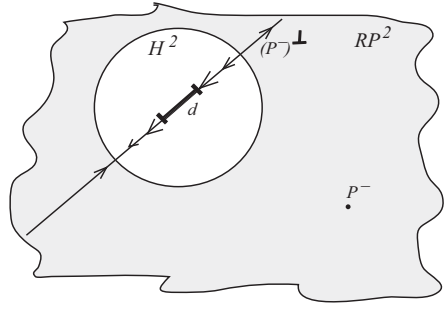,height=4cm}\caption{The action of $A=\cosh
\left( \frac{d}{2}\right) +\sinh \left( \frac{d}{2} \right) J$ in
$RP^{2}$.}\label{fcase2}
\end{figure}

\item[Case 3] $N(A)=1$, $N(A^{-})=0$. Then $A^{-}$ belongs to the nullcone, it is a \emph{nullvector}.
Up to conjugation $A=1+I+J$. Then the matrix associated to $c(A)$ is
\begin{equation*}
m(-1,1;1,1,1,0)=
\begin{pmatrix}
1 & -2 & 2 \\
2 & -1 & 2 \\
2 & -2 & 3
\end{pmatrix}
\end{equation*}
The action is a parabolic transformation fixing $A^{-}$. See Figure
\ref {fcase3}.

\begin{figure}[ht]
\epsfig{file=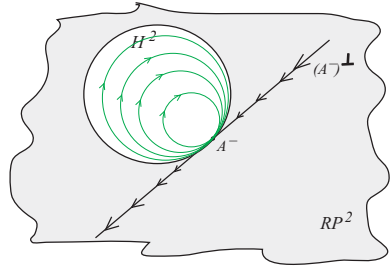,height=3.5cm}\caption{The action of $A=1+I+J$
in $RP^{2}$.}\label{fcase3}
\end{figure}

\end{enumerate}

\subsection{Scalar and vector products on pure quaternions}
\begin{definition}
 The symmetric bilinear form $\langle \quad ,\quad \rangle _{|H_{0}}$
defines the\emph{\ scalar product \ }$s_{p}$ of two pure
quaternions.
\begin{equation*}
\begin{array}{cccc}
s_{p}: & H_{0}\times H_{0} & \longrightarrow & k \\
& (A^{-},B^{-}) & \rightarrow & \left\langle
A^{-},B^{-}\right\rangle
\end{array}
\end{equation*}
where
\begin{equation*}
\left\langle A^{-},B^{-}\right\rangle =\frac{1}{2}\left(
A^{-}\overline{B^{-} }+B^{-}\overline{A^{-}}\right) =\frac{1}{2}
(-A^{-}B^{-}-B^{-}A^{-})=-(A^{-}B^{-})^{+}
\end{equation*}
\end{definition}

\begin{definition}
 The \emph{vector product} $(\times )$ \ of two pure quaternions is
given by
\begin{equation*}
\begin{array}{cccc}
(\times ) & H_{0}\times H_{0} & \longrightarrow & H_{0} \\
& (A^{-},B^{-}) & \rightarrow & A^{-}\times B^{-}:=(A^{-}B^{-})^{-}
\end{array}
\end{equation*}
\end{definition}

Therefore the product $A^{-}B^{-}$ of two pure quaternions decomposes as the
difference between its vector product and its scalar product.
\begin{equation*}
A^{-}B^{-}=-\left\langle A^{-},B^{-}\right\rangle +A^{-}\times B^{-}
\end{equation*}

 Lets compute the vector product of $A^{-}=xi+yj+zij$, and $%
B^{-}=x^{\prime }i+y^{\prime }j+z^{\prime }ij$.
\begin{multline*}
A^{-}\times B^{-} =A^{-}B^{-}+\left\langle A^{-},B^{-}\right\rangle = \\
=(xi+yj+zij)(x^{\prime }i+y^{\prime }j+z^{\prime }ij)+(x,y,z)
\begin{pmatrix}
-\mu & 0 & 0 \\
0 & -\nu & 0 \\
0 & 0 & \mu \nu
\end{pmatrix}
\begin{pmatrix}
x^{\prime } \\
y^{\prime } \\
z^{\prime }
\end{pmatrix}
\\
=-\nu (yz^{\prime }-y^{\prime }z)i -\mu (zx^{\prime }-z^{\prime
}x)j+(xy^{\prime }-x^{\prime }y)ij =
\begin{vmatrix}
-\nu i & -\mu j & ij \\
x & y & z \\
x^{\prime } & y^{\prime } & z^{\prime }
\end{vmatrix}
\end{multline*}
The above formula implies that the vector product is anticommutative:%
\begin{equation*}
A^{-}\times B^{-}=-B^{-}\times A^{-}.
\end{equation*}

\begin{proposition}
The vector product $A^{-}\times B^{-}$ is orthogonal to the plane generated
by $A^{-}$ and $B^{-}.$
\end{proposition}

\begin{proof}
Suppose $A^{-}=xi+yj+zij$, and $B^{-}=x^{\prime }i+y^{\prime }j+z^{\prime }ij
$. Then
\begin{multline*}
\left\langle A^{-}\times B^{-},A^{-}\right\rangle  =\\
=(-\nu (yz^{\prime }-y^{\prime }z),-\mu (zx^{\prime }-z^{\prime
}x),xy^{\prime }-x^{\prime }y)
\begin{pmatrix}
-\mu  & 0 & 0 \\
0 & -\nu  & 0 \\
0 & 0 & \mu \nu
\end{pmatrix}
\begin{pmatrix}
x \\
y \\
z
\end{pmatrix}
\\
=(\mu \nu x(yz^{\prime }-y^{\prime }z)+\mu \nu y(zx^{\prime
}-z^{\prime }x)+\mu \nu z(xy^{\prime }-x^{\prime }y))=0;
\end{multline*}
\begin{multline*}
\left\langle A^{-}\times B^{-},B^{-}\right\rangle  =\\
=(-\nu (yz^{\prime }-y^{\prime }z),-\mu (zx^{\prime }-z^{\prime
}x),xy^{\prime }-x^{\prime }y)
\begin{pmatrix}
-\mu  & 0 & 0 \\
0 & -\nu  & 0 \\
0 & 0 & \mu \nu
\end{pmatrix}
\begin{pmatrix}
x^{\prime } \\
y^{\prime } \\
z^{\prime }
\end{pmatrix}
\\
=(\mu \nu x^{\prime }(yz^{\prime }-y^{\prime }z)+\mu \nu y^{\prime
}(zx^{\prime }-z^{\prime }x)+\mu \nu z^{\prime }(xy^{\prime
}-x^{\prime }y))=0
\end{multline*}
\end{proof}

\begin{corollary}
\label{cambmmortoambm}The pure quaternion $\left( A^{-}B^{-}\right)
^{-}$ is orthogonal to $A^{-}$ and $B^{-}.$\qed
\end{corollary}

\section{A pair of conjugate unit quaternions\label{Sec2}}

We are interested in subgroups of the group of unit quaternions
generated by two conjugate elements. Therefore we shall analyze the
properties of a pair of conjugate unit quaternions.

\begin{lemma}\label{lema1}
Let $A,B\in U_{1}$, $A,B$ conjugate elements in $U_{1}\subset H=\left( \frac{%
\mu ,\nu }{k}\right) $. Then $A^{+}=B^{+}$.
\end{lemma}

\begin{proof}
If $A,B$ are conjugate elements in $U_{1}$, there exist an element $C\in
U_{1}$ such that%
\begin{equation*}
CAC^{-1}=B
\end{equation*}%
Recall that for any quaternion $A\in H$, $T(A)=2A^{+}$. Then it is enough to
prove that $T(A)=T(B)$.%
\begin{eqnarray*}
T(B)
&=&T(CAC^{-1})=T(CA\overline{C})=CA\overline{C}+\overline{CA\overline{C}}=
CA\overline{C}+C\overline{A}\,\overline{C
}\\
&=&C(A+\overline{A})\overline{C}
=CT(A)\overline{C}=T(A)C\overline{C}=T(A)
\end{eqnarray*}
\end{proof}

Next we prove that the conjugation by a unit quaternion preserves the scalar
product of pure quaternions.

\begin{lemma}\label{lema2}
Let $(A,B),(A_{0},B_{0})$ be two pairs of unit quaternions, that are
conjugate (there exist a unit quaternion $C$ such that $CA\overline{C}=A_{0}$%
, and $CB\overline{C}=B_{0}$). Then $
(A^{-}B^{-})^{+}=(A_{0}^{-}B_{0}^{-})^{+}$.
\end{lemma}

\begin{proof}
First note that $(CA\overline{C})^{-}=CA^{-}\overline{C}$. This is
because
\begin{equation*}
CA\overline{C}=C(A^{+}+A^{-})\overline{C}=A^{+}+CA^{-}\overline{C}
\end{equation*}
and
\begin{equation*}
CA\overline{C}=(CA\overline{C})^{+}+(CA\overline{C})^{-}=A^{+}+(CA\overline{C
})^{-}
\end{equation*}
yields to  $(CA\overline{C})^{-}=CA^{-}\overline{C}$.

Using this property, one obtain
\begin{multline*}
(A_{0}^{-}B_{0}^{-})^{+}=((CA\overline{C})^{-}(CB\overline{C}
)^{-})^{+}=(CA^{-}\overline{C}CB^{-}\overline{C})^{+}\\
=(CA^{-}B^{-}\overline{C })^{+}=(A^{-}B^{-})^{+}
\end{multline*}
where the last equality is true by applying Lemma \ref{lema1} to the
conjugate quaternions $A^{-}B^{-}$ and $CA^{-}B^{-}\overline{C }$.
\end{proof}

The following two theorems show that a pair of conjugate (in
$U_{1}$) unit quaternions in a quaternion algebra over $\mathbb{R}$
are determined, up to conjugation, by two real numbers: the trace of
them and the scalar product of their pure parts. First we analyze
the Hamilton quaternions.

\begin{theorem}
\label{torema1A} A pair $A,B$ of unit quaternions in
$\mathbb{H}=\left( \frac{-1,-1}{\mathbb{R}}\right) $ with the same
trace ($2x$) are determined, up to conjugation in $U_{1}$, by the
real number $y=-(A^{-}B^{-})^{+}$.
\end{theorem}

\begin{proof}
Write $A=x+A^{-}$, $B=x+B^{-}$. Note that
\begin{eqnarray*}
N(A^{-}) &=&N(A-A^{+})=N(A)-N(A^{+})=1-x^{2} \\
N(B^{-}) &=&N(B-B^{+})=N(B)-N(B^{+})=1-x^{2}
\end{eqnarray*}%
implies that the pure quaternions $A^{-}$ and $B^{-}$ have the same
norm, thus both are in the sphere of radius $\sqrt{1-x^{2}}$. Using
conjugation by unit quaternions in $E^{3}(\cong H_{0})$, Proposition \ref{proposicion5} and Lemma
\ref{lema2}, we can choose one pair $A,B$ with simpler coordinates
as follows. Consider  the coordinates system $\left\{ X,Y,Z\right\} $ in $E^{3}$
in the basis $\left\{ -ij,j,i\right\} $. First, conjugating by a unit quaternion $C$, such that
$CB^{-}\overline{C}$ is in the positive $Z$-axis: $CB^{-}\overline{C}=\gamma i$,
$\gamma
>0$.
Let $D$ be a unit
quaternion such that the conjugation in $E^{3}(\cong H_{0})$ by $D$
is a rotation around the $Z$-axis mapping $CA^{-}\overline{C}$ to a
vector in the plane $YZ$  with positive $Y$- coordinate. See Figure
\ref{AmBm}. Then $DCB^{-}\overline{C}\overline{D}=CB^{-}\overline{C}$ and $DCA\overline{CD}=\alpha
i+\beta j$, $\beta \geqslant 0$.

\begin{figure}[ht]
\epsfig{file=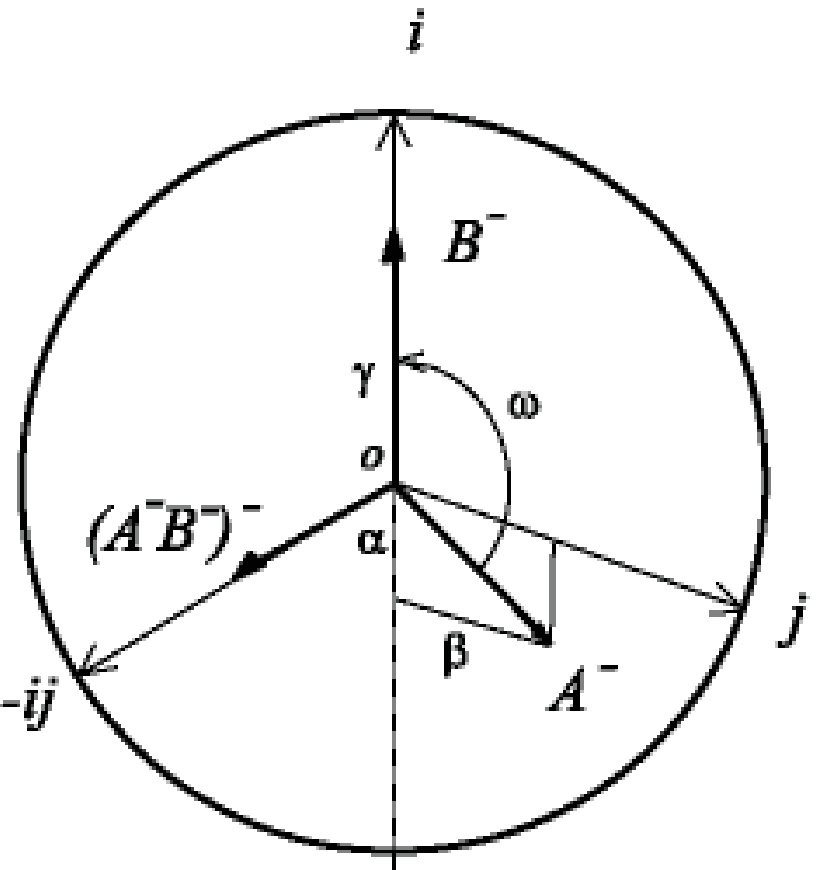,height=4cm}\caption{The vectors $A^{-}$ and $B^{-}$ in $E^{3}$}\label{AmBm}
\end{figure}

Therefore we assume
\begin{eqnarray*}
A^{-}&=&\alpha i+\beta j,\quad \beta \geqslant 0\\
B^{-}&=&\gamma i,\quad \gamma >0
\end{eqnarray*}
Note that
\begin{eqnarray*}
N(A^{-})&=&N((\alpha i+\beta j)=(\alpha i+\beta j)(-\alpha i-\beta j)=(\alpha
^{2}+\beta ^{2})\\
N(B^{-})&=&N(\gamma i)=(\gamma i)(-\gamma i)=\gamma ^{2}
\end{eqnarray*}
Then
\begin{eqnarray}
\alpha ^{2}+\beta ^{2}&=&1-x^{2}  \label{ealfabeta1}\\
\gamma ^{2}&=&1-x^{2}\label{egama}
\end{eqnarray}

From the computation of $A^{-}B^{-}$
\begin{equation*}
A^{-}B^{-}=(\alpha i+\beta j)(\gamma i)=-\gamma \alpha -\gamma \beta ij
\end{equation*}
one obtains
\begin{equation}
y=-(A^{-}B^{-})^{+}=\gamma \alpha  \label{eygamaalfa}
\end{equation}%
\begin{equation*}
(A^{-}B^{-})^{-}=-\gamma \beta ij
\end{equation*}%
Solving the equations \ref{ealfabeta1}, \ref{egama} and \ref{eygamaalfa} we
obtain
\begin{eqnarray*}
\gamma &=& \sqrt{1-x^{2}}\\
\alpha &=&\frac{y}{\gamma }=\frac{y}{\sqrt{1-x^{2}}}\\
\beta &=&+\sqrt{1-x^{2}-\frac{y^{2}}{\gamma
^{2}}}=\frac{1}{\sqrt{1-x^{2}}}\sqrt{(1-x^{2})^{2}-y^{2}}
\end{eqnarray*}
Then, up to conjugation in $U_{1}$
\begin{equation}
\fbox{$
\begin{array}{l}
A=x+\frac{1}{+\sqrt{1-x^{2}}}\left( yi+\sqrt{(1-x^{2})^{2}-y^{2}}j\right) \\
B=x+\sqrt{1-x^{2}}i,\quad \sqrt{1-x^{2}}>0%
\end{array}
$}  \label{eABteo1}
\end{equation}%
Observe also that
\begin{equation*}
(A^{-}B^{-})^{-}=-\sqrt{(1-x^{2})^{2}-y^{2}}ij
\end{equation*}%
The angle $\omega $ subtended by $A^{-}$ and $B^{-}$ (see Figure \ref{AmBm})
is
\begin{equation*}
\cos \omega =\frac{\left\langle A^{-},B^{-}\right\rangle }{\sqrt{
\left\langle A^{-},A^{-}\right\rangle }\sqrt{\left\langle
B^{-},B^{-}\right\rangle }}=\frac{-(A^{-}B^{-})^{+}}{\gamma \gamma }=\frac{y%
}{1-x^{2}}
\end{equation*}
\end{proof}

\begin{remark}
The multiplicative group generated by the elements $A$ and $B$ (\ref{eABteo1}
) of the above Theorem \ref{torema1A} belongs to the quaternion algebra $
\left( \frac{-1,-1}{k}\right) $, where
\begin{equation*}
k=Q\left( x,y,\sqrt{1-x^{2}},\sqrt{\left( 1-x^{2}\right) ^{2}-y^{2}}\right)
\end{equation*}
\end{remark}

The case \ $M(2,\mathbb{R})=\left( \frac{-1,1}{\mathbb{R}}\right) $ is
studied in the following theorem.

\begin{theorem}
\label{teorema2A} A pair $A,B$ of conjugate (in $U_{1}$) independent unit
quaternions in $M(2,\mathbb{R})=\left( \frac{-1,1}{\mathbb{R}}\right) $ with
the same trace ($2x)$ are determined, up to conjugation in $U_{\pm 1}$, by
the real number $y=-(A^{-}B^{-})^{+}$ .
\end{theorem}

\begin{proof}
Write $A=x+A^{-}$, $B=x+B^{-}$. There exists an element $C\in U_{1}$ such
that
\begin{equation*}
CA\overline{C}=C(x+A^{-})\overline{C}=x+CA^{-}\overline{C}=B=x+B^{-}
\end{equation*}%
Then, the vectors $A^{-},B^{-}\in E^{1,2}$ are also conjugate
\begin{equation*}
CA^{-}\overline{C}=B^{-}
\end{equation*}%

The exact sequence (\ref{eSO12}) shows that conjugation by $C$ acts in $%
E^{1,2}$ as an element $c(C)$ of $SO^{+}(1,2)$.

Using conjugation and Lemma \ref{lema2} we can choose one pair $A,B$
with convenient coordinates as follows.

Note that
\begin{equation*}
N(A^{-})=N(B^{-})=1-x^{2}
\end{equation*}%
Then we distinguish three cases according to the value of $N(A^{-})(>,<,=)0$,
which geometrically fixes the position of $A^{-}$ and $B^{-}$ with respect
to the nullcone. If we write $\left\{ X,Y,Z\right\} $ to denote the
coordinates in the basis $\left\{ -IJ,J,I\right\} $, $A^{-}\ $and $B^{-}$
belong to the hyperboloid
\begin{equation*}
-X^{2}-Y^{2}+Z^{2}=1-x^{2}=N(A^{-})=N(B^{-})
\end{equation*}%
which is a two-sheeted hyperboloid if $x^{2}<1,$ the nullcone if $x^{2}=1$
and a one-sheeted hyperboloid if $x^{2}>1$.

\begin{enumerate}
\item[Case 1] $N(A^{-})=N(B^{-})>0$. $(x^{2}=1-N(A^{-})<1$. Since $A^{-}$
and $B^{-}$ are conjugate by an element of $U_{1}$, they are placed in the same
sheet of a two-sheeted hyperboloid in the interior of the nullcone. Up to
conjugation by $j\in U_{\pm 1}$, if necessary, we assume that the time-like vectors $A^{-}$ and $B^{-}$ are in the upper component of the hyperboloid. Because $SO^{+}(1,2)$
acts transitively on the rays inside the nullcone, we may assume that
\begin{equation*}
A^{-}=\gamma I,\quad \gamma >0
\end{equation*}
A suitable conjugation by an elliptic rotation around the $Z$ axis$,$ maps $
B^{-}$ to a vector in the $YZ$ plane with positive $Y$ coordinate. Say%
\begin{equation*}
B^{-}=\alpha I+\beta J,\quad \alpha >\beta >0
\end{equation*}%
The following computations give the expressions of $\alpha ,\beta $ and $
\gamma $ as functions of $x$ and $y$:
\begin{eqnarray*}
1 &=&A\overline{A}=(A^{+}+A^{-})(A^{+}-A^{-})=x^{2}-A^{-}A^{-}\quad  \\
&\Longrightarrow &\quad 0<-A^{-}A^{-}=1-x^{2}=-(\gamma I)(\gamma I)=\gamma
^{2}
\end{eqnarray*}
\begin{equation*}
\Longrightarrow \gamma =+\sqrt{1-x^{2}}
\end{equation*}
\begin{equation*}
0<-B^{-}B^{-}=1-x^{2}=-(\alpha I+\beta J)(\alpha I+\beta J)=-(-\alpha
^{2}+\beta ^{2})
\end{equation*}
\begin{equation*}
\Longrightarrow 1-x^{2}=\alpha ^{2}-\beta ^{2}
\end{equation*}
\begin{eqnarray*}
A^{-}B^{-} &=&(\gamma I)(\alpha I+\beta J)=-\gamma \alpha +\gamma \beta IJ \\
&\Longrightarrow &(A^{-}B^{-})^{+}=-\gamma \alpha ,\qquad
(A^{-}B^{-})^{-}=\gamma \beta IJ
\end{eqnarray*}
\begin{equation*}
\Longrightarrow y=\gamma \alpha
\end{equation*}
Then, up to conjugation
\begin{equation}
\fbox{$
\begin{array}{c}
A=x+\sqrt{1-x^{2}}I,\quad \sqrt{1-x^{2}}>0 \\
B=x+\frac{1}{+\sqrt{1-x^{2}}}\left( yI+\sqrt{y^{2}-(1-x^{2})^{2}}J\right)
\end{array}
$}  \label{eABcaso1}
\end{equation}
\begin{equation*}
(A^{-}B^{-})^{-}=+\sqrt{y^{2}-(1-x^{2})^{2}}IJ
\end{equation*}
The hyperbolic distance $d$ between the projection of $A^{-}$ and $B^{-}$ on
the hyperbolic plane (pure unit quaternions in the upper component) (see
Figure \ref{AmBm2}) is
\begin{equation*}
\cosh d=\frac{\left\langle A^{-},B^{-}\right\rangle }{\sqrt{\left\langle
A^{-},A^{-}\right\rangle }\sqrt{\left\langle B^{-},B^{-}\right\rangle }}=
\frac{-(A^{-}B^{-})^{+}}{\gamma \gamma }=\frac{y}{1-x^{2}}
\end{equation*}

\begin{figure}[ht]
\epsfig{file=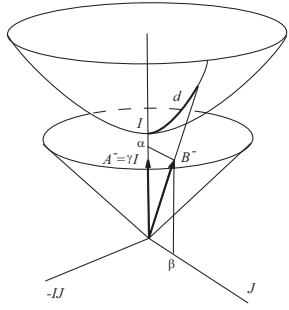,height=4cm}\caption{The vectors $A^{-},$ $B^{-}$ in
$E^{1,2}$, case 1.}\label{AmBm2}
\end{figure}

The above quantity is larger than 1 and therefore is a hyperbolic cosine:
\begin{equation*}
\frac{y}{1-x^{2}}=\frac{\gamma \alpha }{\gamma ^{2}}=\frac{\alpha }{\gamma },%
\text{but }\gamma ^{2}=1-x^{2}=\alpha ^{2}-\beta ^{2}\Longrightarrow \gamma
^{2}<\alpha ^{2}\Longrightarrow 0<\gamma <\alpha
\end{equation*}

\item[Case 2] $N(A^{-})=N(B^{-})<0$. $(x^{2}=1-N(A^{-})>1).$ The vectors $
A^{-},$ $B^{-}$ are space-like vectors (outside the nullcone). Recall that $A^{-},$ $B^{-}$
belong to a one-sheeted hyperboloid. There are three subcases according to
the relative position of the plane $\Pi $ defined by $A^{-},$ $B^{-}$ and
the nullcone. We shall prove that these cases correspond to $y^{2}(>,<,=)\gamma ^{4}$.

\item[Subcase 2A] The plane $\Pi $ intersects the nullcone. See Figure \ref{fsubcase2a1}. Considering the projective plane $RP^{2}$ defined by the
lines through the origin in $E^{1,2}$, the nullcone defines a conic $\mathcal{C}$. The polars of $A^{-}$ and $B^{-}$ with respect to the conic $\mathcal{C}$ intersect in a point $Q$ (pole of $\Pi $) outside $\mathcal{C}$. Up to conjugation we can suppose that $Q=\infty $. Then there are two
possibilities : $A^{-},$ $B^{-}$ are either on the same side or on different side of $\mathcal{C}$. See Figure \ref{fsubcase2a2} for the situation in $E^{1,2}$.

\begin{figure}[ht]
\epsfig{file=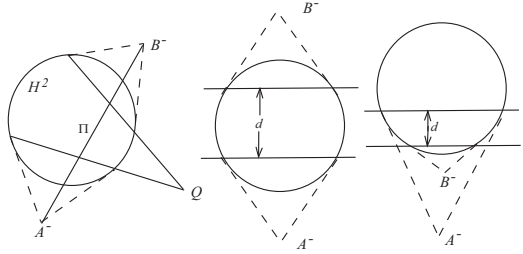,height=3.5cm}\caption{Subcase 2A.}\label{fsubcase2a1}
\end{figure}

Up to conjugation
\begin{eqnarray*}
A^{-} &=&\gamma J,\quad \gamma >0 \\
B^{-} &=&\alpha I+\beta J,\quad \left\{
\begin{array}{c}
\beta >\alpha >0 \\
\beta <\alpha <0
\end{array}
\right\}
\end{eqnarray*}

\begin{figure}[ht]
\epsfig{file=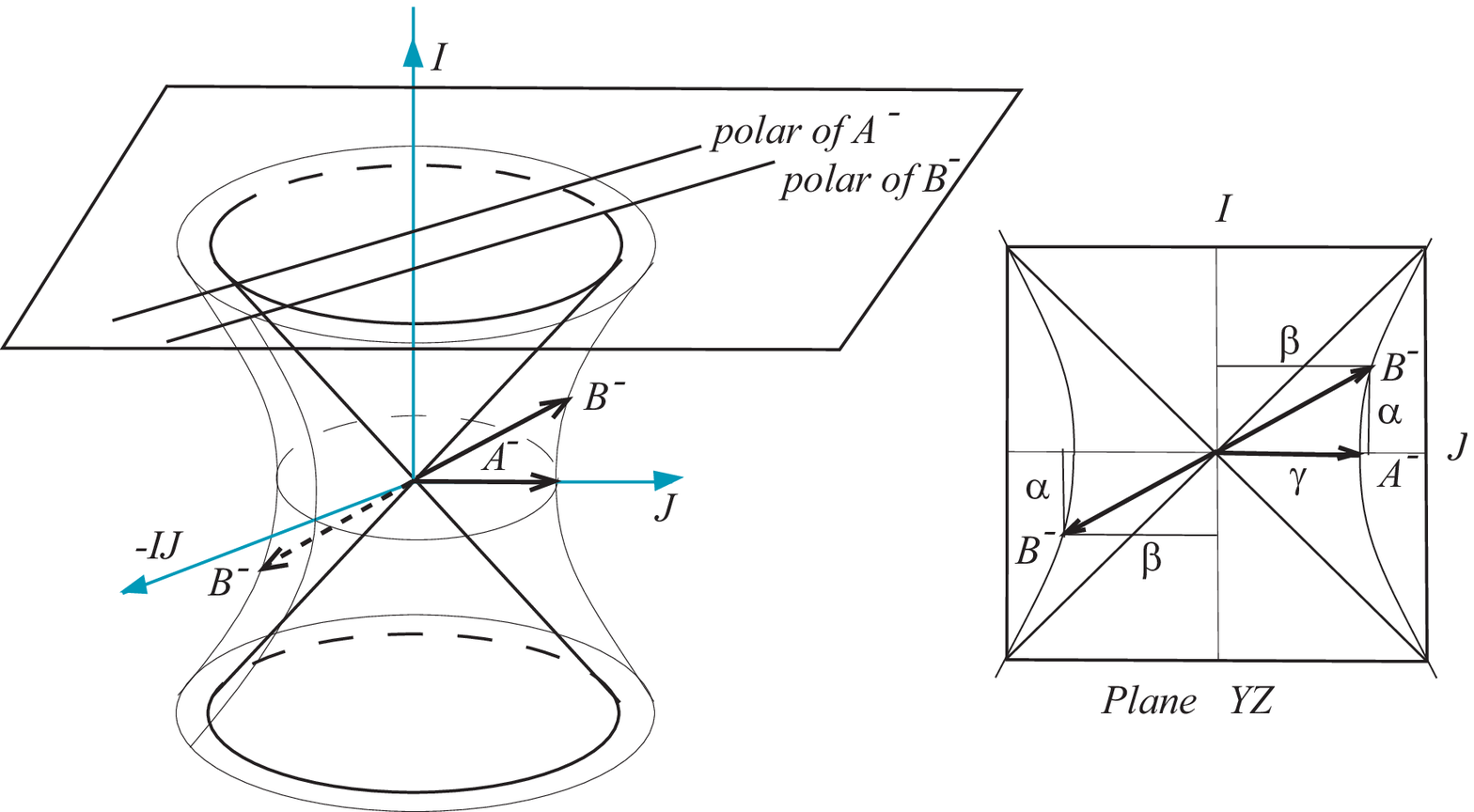,height=4.5cm}\caption{The vectors $A^{-},$ $B^{-}$ in
$E^{1,2}$, subcase 2A.}\label{fsubcase2a2}
\end{figure}

The following computations
give the expressions of $\alpha ,\beta $ and $\gamma $ as functions of $x$
and $y:$
\begin{eqnarray*}
\quad 0 &>&N(A^{-})=-A^{-}A^{-}=1-x^{2}=-(\gamma J)(\gamma J)=-\gamma ^{2} \\
0 &>&N(B^{-})=-B^{-}B^{-}=1-x^{2}=-(\alpha I+\beta J)(\alpha I+\beta
J)=-(-\alpha ^{2}+\beta ^{2})
\end{eqnarray*}
\begin{equation*}
\Longrightarrow
\begin{array}{c}
\gamma =+\sqrt{x^{2}-1} \\
-\alpha ^{2}+\beta ^{2}=x^{2}-1
\end{array}%
\end{equation*}
\begin{eqnarray*}
A^{-}B^{-} &=&(\gamma J)(\alpha I+\beta J)=\gamma \beta -\gamma \alpha IJ \\
&\Longrightarrow &(A^{-}B^{-})^{+}=\gamma \beta ,\qquad
(A^{-}B^{-})^{-}=-\gamma \alpha IJ
\end{eqnarray*}
\begin{equation*}
\Longrightarrow y=-\gamma \beta
\end{equation*}
\begin{equation*}
\Longrightarrow \left\{
\begin{array}{c}
y<0,\ \beta =\frac{y}{-\gamma }>0,\ \Longrightarrow \ \alpha >0\
\Longrightarrow \ \alpha =+\frac{1}{\gamma }\sqrt{y^{2}-\gamma ^{4}} \\
y>0,\ \beta =\frac{y}{-\gamma }<0,\ \Longrightarrow \ \alpha <0\
\Longrightarrow \ \alpha =-\frac{1}{\gamma }\sqrt{y^{2}-\gamma ^{4}}
\end{array}
\right\}
\end{equation*}%
Therefore, the two possibilities are given by the sign of $y$. For $y<0$.%
\begin{equation}
\frame{$
\begin{array}{c}
A=x+\sqrt{x^{2}-1}J,\quad \sqrt{x^{2}-1}>0 \\
B=x+\frac{1}{+\sqrt{x^{2}-1}}\left( \sqrt{y^{2}-(x^{2}-1)^{2}}I-yJ\right)
\\
(A^{-}B^{-})^{-}=-\sqrt{y^{2}-(x^{2}-1)^{2}}IJ \\
\cosh d=\frac{-y}{x^{2}-1}
\end{array}
$}  \label{eABcaso2A1}
\end{equation}%
\qquad \qquad For $y>0$.%
\begin{equation}
\frame{$%
\begin{array}{c}
A=x+\sqrt{x^{2}-1}J,\quad \sqrt{x^{2}-1}>0 \\
B=x-\frac{1}{+\sqrt{x^{2}-1}}\left( \sqrt{y^{2}-(x^{2}-1)^{2}}I+yJ\right)
\\
(A^{-}B^{-})^{-}=+\sqrt{y^{2}-(x^{2}-1)^{2}}IJ \\
\cosh d=\frac{y}{x^{2}-1}%
\end{array}%
$}  \label{eABcaso2A2}
\end{equation}
Where $d$ is the distance between the polars of $A^{-}$ and $B^{-}$. Note
that
\begin{equation*}
\cosh ^{2}d=\frac{\left\langle A^{-},B^{-}\right\rangle ^{2}}{\left\langle
A^{-},A^{-}\right\rangle \left\langle B^{-},B^{-}\right\rangle }=\frac{y^{2}
}{(x^{2}-1)^{2}}>1\Longleftrightarrow y^{2}>\gamma ^{4}
\end{equation*}%
because the polars of $A^{-}$ and $B^{-}$ do not intersect inside the
conic.

\item[Subcase 2B] The plane $\Pi $ does not intersect the nullcone. See
Figure \ref{subcase2b1} for the situation in projective plane $RP^{2}$
defined by the lines through the origin in $E^{1,2}$. The polars of $A^{-}$
and $B^{-}$ with respect to the conic $\mathcal{C}$ intersect in a point $Q$
(pole of $\Pi $) inside $\mathcal{C}$.

\begin{figure}[ht]
\epsfig{file=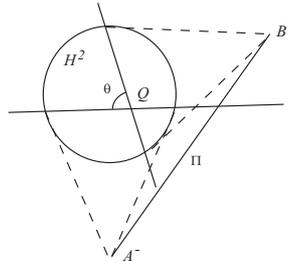,height=3.5cm}\caption{Subcase 2B.}\label{subcase2b1}
\end{figure}

 Up to conjugation we
can suppose that $Q$ is the center of the conic $\mathcal{C}$. The angle $
\theta $ is given by
\begin{equation*}
\cos \theta =\frac{\left\langle A^{-},B^{-}\right\rangle }{\sqrt{
\left\langle A^{-},A^{-}\right\rangle }\sqrt{\left\langle
B^{-},B^{-}\right\rangle }}
\end{equation*}%
See Figure \ref{subcase2b2} for the situation in $E^{1,2}$.

\begin{figure}[ht]
\epsfig{file=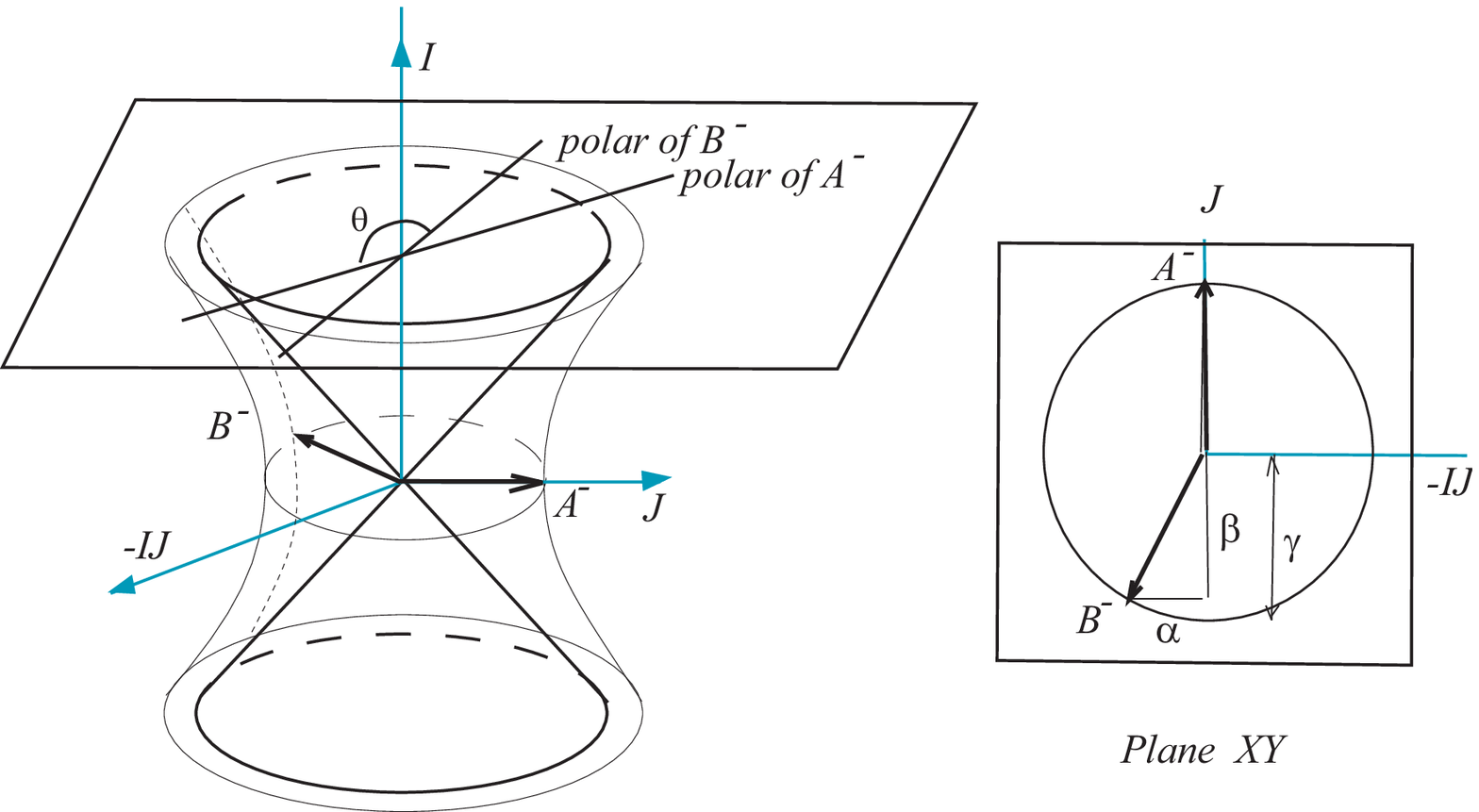,height=5cm}\caption{The vectors \ $A^{-}$, $B^{-}$ in $E^{1,2}$,
subcase 2B.}\label{subcase2b2}
\end{figure}

The vectors $A^{-}$ and $%
B^{-}$ are in the $XY$ plane. Up to conjugation,
\begin{eqnarray*}
A^{-} &=&\gamma J,\quad \gamma >0 \\
B^{-} &=&\alpha IJ+\beta J
\end{eqnarray*}%
Then%
\begin{eqnarray*}
\quad 0 &>&N(A^{-})=-A^{-}A^{-}=1-x^{2}=-(\gamma J)(\gamma J)=-\gamma ^{2} \\
0 &>&N(B^{-})=-B^{-}B^{-}=1-x^{2}=-(\alpha IJ+\beta J)(\alpha IJ+\beta
J)=-\alpha ^{2}-\beta ^{2}
\end{eqnarray*}%
\begin{equation*}
\Longrightarrow
\begin{array}{c}
\gamma =+\sqrt{x^{2}-1} \\
\alpha ^{2}+\beta ^{2}=x^{2}-1%
\end{array}%
\end{equation*}%
\begin{eqnarray*}
A^{-}B^{-} &=&(\gamma J)(\alpha IJ+\beta J)=-\gamma \alpha I+\gamma \beta  \\
&\Longrightarrow &(A^{-}B^{-})^{+}=\gamma \beta ,\qquad
(A^{-}B^{-})^{-}=-\gamma \alpha I
\end{eqnarray*}%
\begin{equation*}
\Longrightarrow y=-\gamma \beta
\end{equation*}%
\begin{equation*}
\cos^{2} \theta =\frac{\left\langle A^{-},B^{-}\right\rangle^{2} }{
\left\langle A^{-},A^{-}\right\rangle \left\langle
B^{-},B^{-}\right\rangle }=\frac{y^{2}}{(x^{2}-1)^{2}}<1\Longleftrightarrow
y^{2}<(x^{2}-1)^{2}=\gamma ^{4}
\end{equation*}%
\begin{eqnarray*}
\beta  &=&\frac{-y}{\sqrt{x^{2}-1}} \\
\alpha ^{2} &=&x^{2}-1-\beta ^{2}=\frac{(x^{2}-1)^{2}-y^{2}}{(x^{2}-1)}>0 \\
\alpha  &=&\pm \sqrt{\frac{(x^{2}-1)^{2}-y^{2}}{(x^{2}-1)}}
\end{eqnarray*}%
There are two possibilities for $\alpha $. See Figure \ref{subcase2b3} .

\begin{figure}[ht]
\epsfig{file=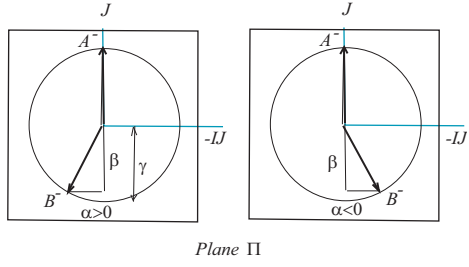,height=3.5cm}\caption{The two possibilities for
$\protect\alpha $ in the plane $XY$.}\label{subcase2b3}
\end{figure}

Conjugation by $J\in U_{\pm
1}$, acts in $A^{-},B^{-}$ as follows. Note that $J^{2}=1\Longrightarrow J^{-1}=J$.
\begin{equation*}
J(A^{-},B^{-})J=(J\gamma JJ,\ J(\alpha IJ+\beta J)J)=(\gamma J,-\alpha
IJ+\beta J)
\end{equation*}
Therefore, up to conjugation we can choose $\alpha >0$. Then, in this subcase
\begin{equation}
\fbox{$
\begin{array}{c}
A=x+\sqrt{x^{2}-1}J,\quad \sqrt{x^{2}-1}>0 \\
B=x+\frac{1}{+\sqrt{x^{2}-1}}\left( -yJ+\sqrt{(x^{2}-1)^{2}-y^{2}}IJ\right)
\\
(A^{-}B^{-})^{-}=-\sqrt{(x^{2}-1)^{2}-y^{2}}IJ \\
\cos \theta =\frac{-y}{x^{2}-1}\Longrightarrow \left\{
\begin{array}{c}
y<0\Longleftrightarrow 0<\theta <\frac{\pi }{2} \\
y>0\Longleftrightarrow \frac{\pi }{2}<\theta <\pi
\end{array}
\right\}
\end{array}%
$}  \label{eABcaso2B}
\end{equation}

\item[Subcase 2C] The plane $\Pi $ is tangent to the nullcone. See Figure
\ref{subcase2c1}. Assume $\Pi $ is tangent to $\mathcal{N}$ and contains the
$Y$ axis. Then, we can assume that $A^{-},B^{-}$ are in the plane $X=Z$,
where $\left\{ X,Y,Z\right\} $ correspond to coordinates in the basis $%
\left\{ -IJ,J,I\right\} $.

\begin{figure}[ht]
\epsfig{file=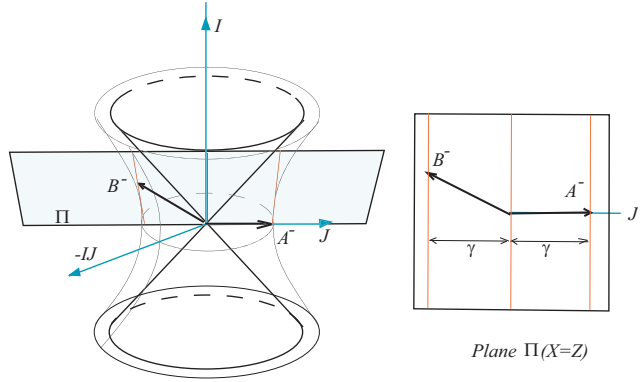,height=5cm}\caption{Subcase 2C.}\label{subcase2c1}
\end{figure}

Up to conjugation by a
parabolic transformation fixing $\Pi $, we may assume that
\begin{eqnarray*}
A^{-} &=&\gamma J,\quad \gamma >0 \\
B^{-} &=&\lambda (IJ+I)+\gamma ^{\prime }J
\end{eqnarray*}%
Then%
\begin{eqnarray*}
\quad 0 >N(A^{-})&=&-A^{-}A^{-}=1-x^{2}=-(\gamma J)(\gamma J)=-\gamma ^{2} \\
0 >N(B^{-})&=&-B^{-}B^{-}=1-x^{2}=-(\lambda (IJ+I)+\gamma ^{\prime
}J)(\lambda (IJ+I)+\gamma ^{\prime }J)\\
&=&-(\gamma ^{\prime })^{2}
\end{eqnarray*}
implies that $\gamma ^{\prime }=\pm \gamma $ , then we write $\gamma
^{\prime }=\varepsilon \gamma $ , where $\varepsilon =\pm 1$ . This also can
be obtained geometrically, because $A^{-},B^{-}$ are both in the one-sheeted
hyperboloid $-X^{2}-Y^{2}+Z^{2}=1-x^{2}=N(A^{-})=N(B^{-})=-\gamma ^{2}.$
Then, $A^{-},B^{-}$ are in
\begin{equation*}
\begin{Bmatrix}
-X^{2}-Y^{2}+Z^{2}=-\gamma ^{2} \\
X=Z
\end{Bmatrix}
=
\begin{Bmatrix}
Y^{2}=\gamma ^{2} \\
X=Z
\end{Bmatrix}
=
\begin{Bmatrix}
Y=\pm \gamma \\
X=Z
\end{Bmatrix}
\end{equation*}
Assume $\lambda >0$, then the conjugation by the unit quaternion
\begin{equation*}
C=\frac{1+\lambda }{2\sqrt{\lambda }}+\frac{\lambda -1}{2\sqrt{\lambda }}J
\end{equation*}
acts in $H_{0}$ as a hyperbolic rotation about the Y-axis $J$. See Figure \ref
{girohyper}. The matrix of this transformation is
\begin{equation*}
c(C)=
\begin{pmatrix}
\frac{1+\lambda ^{2}}{2\lambda } & 0 & \frac{1-\lambda ^{2}}{2\lambda } \\
0 & 1 & 0 \\
\frac{1-\lambda ^{2}}{2\lambda } & 0 & \frac{1+\lambda ^{2}}{2\lambda }%
\end{pmatrix}%
\end{equation*}%
Therefore%
\begin{eqnarray*}
c(C)(\gamma J) &=&\gamma J \\
c(C)(\lambda (IJ+I)+\varepsilon \gamma J) &=&(IJ+I)+\varepsilon \gamma J
\end{eqnarray*}%
If $\lambda <0$ conjugation by
\begin{equation*}
D=\frac{1-\lambda }{2\sqrt{-\lambda }}+\frac{-\lambda -1}{2\sqrt{-\lambda }}J
\end{equation*}%
acts as
\begin{eqnarray*}
c(D)(\gamma J) &=&\gamma J \\
c(D)(\lambda (IJ+I)+\varepsilon \gamma J) &=&(IJ+I)+\varepsilon \gamma J
\end{eqnarray*}
Then, up to conjugation, we can assume
\begin{eqnarray*}
A^{-} &=&\gamma J,\quad \gamma >0 \\
B^{-} &=&(IJ+I)+\gamma ^{\prime }J
\end{eqnarray*}
To relate the parameters to $x,y,$ we compute as usual
\begin{eqnarray*}
A^{-}B^{-} &=&(\gamma J)((IJ+I)+\varepsilon \gamma J)=-\gamma I-\gamma
IJ+\gamma ^{2}\varepsilon \\
&\Longrightarrow &(A^{-}B^{-})^{+}=\gamma ^{2}\varepsilon ,\qquad
(A^{-}B^{-})^{-}=-\gamma I-\gamma IJ
\end{eqnarray*}
\begin{equation*}
\Longrightarrow y=-\varepsilon \gamma ^{2}
\end{equation*}
\begin{equation*}
\varepsilon =\frac{-y}{\gamma ^{2}},\qquad \gamma ^{\prime }=-\frac{y}{%
\gamma }\quad \text{and}\qquad y^{2}=\gamma ^{4}
\end{equation*}
Therefore
\begin{equation}
\fbox{$
\begin{array}{c}
A=x+\sqrt{x^{2}-1}J,\quad \sqrt{x^{2}-1}>0 \\
B=x+(IJ+I)-\frac{y}{\sqrt{x^{2}-1}}J \\
(A^{-}B^{-})^{-}=-\sqrt{(x^{2}-1)}(I+IJ)
\end{array}
$}  \label{eABcaso2C}
\end{equation}

\item[Case 3] $N(A^{-})=N(B^{-})=0$. $(x^{2}=1).$ The vectors $A^{-},$ $B^{-}
$ are in the nullcone. Since $A^{-},$ $B^{-}$ are conjugate in $U_{1}$, they
belong to the same component of $\mathcal{N}-(0).$ Conjugating by $J\in
U_{\pm 1}$ if necessary, assume $A^{-},$ $B^{-}\subset \mathcal{N}_{+}$. Up
to conjugation by a rotation around $I$
 sending $A^{-}$ to $\gamma (I+J)$,
and a parabolic around $A^{-}$, we can assume that $A^{-},$ $B^{-}$ are as follows:
\begin{equation*}
A^{-}=\gamma (I+J),\quad B^{-}=\alpha (I-J)
\end{equation*}%
By using a hyperbolic rotation around $IJ$, we can suppose also that $\gamma
=1$. See Figure \ref{ambm3}.

\begin{figure}[ht]
\epsfig{file=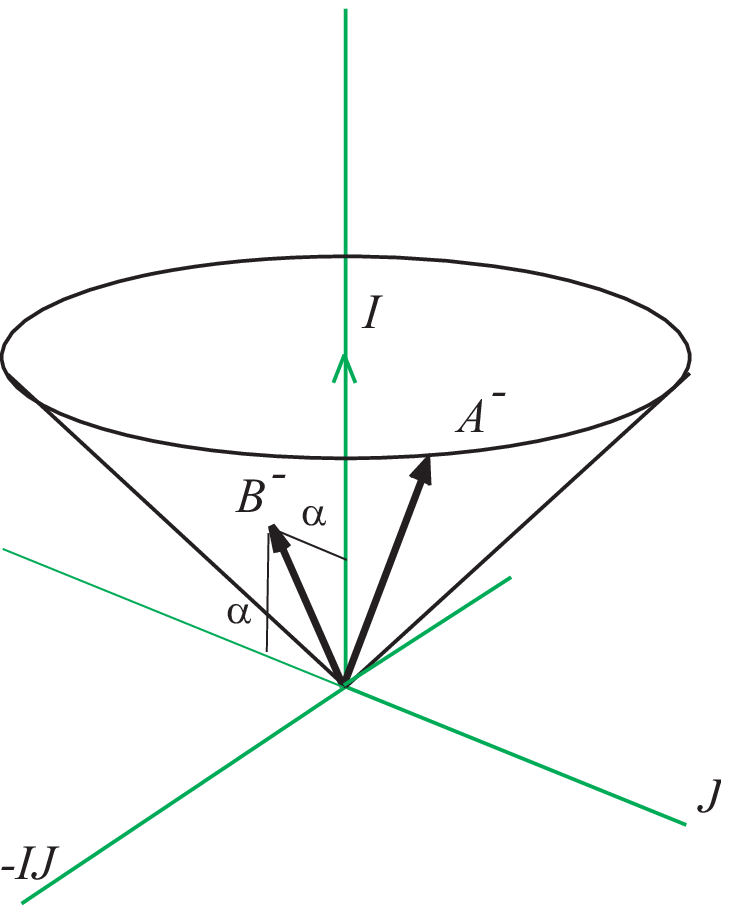,height=4.0cm}\caption{Case 3.}\label{ambm3}
\end{figure}

Then
\begin{eqnarray*}
A^{-}B^{-} &=&\alpha (I+J)(I-J)=-2\alpha -2\alpha IJ \\
&\Longrightarrow &(A^{-}B^{-})^{+}=-2\alpha ,\qquad
(A^{-}B^{-})^{-}=-2\alpha IJ
\end{eqnarray*}%
\begin{equation*}
\Longrightarrow y=2\alpha
\end{equation*}%
Thus%
\begin{equation}
\fbox{$%
\begin{array}{c}
A=x+I+J \\
B=x+\frac{y}{2}(I-J) \\
(A^{-}B^{-})^{-}=-yIJ%
\end{array}%
$}  \label{eABcaso3}
\end{equation}
\end{enumerate}
\end{proof}

\begin{remark}
The multiplicative group generated by the elements $A$ and $B$ of
the above theorem \ref{teorema2A} belong to a quaternion algebra
$\left( \frac{-1,1}{k} \right) $, according to the following cases:

\begin{enumerate}
\item[Case 1] $A$ and $B$ in (\ref{eABcaso1}) $k=Q(x,y,\sqrt{1-x^{2}},\sqrt{
y^{2}-(1-x^{2})^{2}})$

\item[Subcase 2A] $A$ and $B$ in (\ref{eABcaso2A1}) and (\ref{eABcaso2A2}) $
k=Q(x,y,\sqrt{x^{2}-1},\sqrt{y^{2}-(x^{2}-1)^{2}})$

\item[Subcase 2B] $A$ and $B$ in (\ref{eABcaso2B}) $k=Q(x,y,\sqrt{x^{2}-1},
\sqrt{(x^{2}-1)^{2}-y^{2}})$

\item[Subcase 2C] $A$ and $B$ in (\ref{eABcaso2C}) $k=Q(x,y,\sqrt{x^{2}-1})$

\item[Case 3] $A$ and $B$ in (\ref{eABcaso3}) $k=Q(x,y)$
\end{enumerate}
\end{remark}

\begin{proposition}
\label{psubalgebraAB} Let $A,B\in U_{1}$ be conjugate elements in
$U_{1}$, the unit quaternions of a quaternion algebra $H=\left(
\frac{\mu ,\nu }{k} \right) $. The subalgebra of $H$ generated (as
an $k$-algebra) by $A$ and $B$ coincides with the $k$-vector
subspace of $H$ generated by
\begin{equation*}
\left\{ 1,A^{-},B^{-},(A^{-}B^{-})^{-}\right\}.
\end{equation*}
\end{proposition}

\begin{proof}
Set $A=A^{+}+A^{-}$, $B=B^{+}+B^{-}$. Then, since $A$ and $B$ are
conjugate, $A^{+}=B^{+}$ (Lemma \ref{lema1}). Set $x=A^{+}=B^{+}\in
k$. Take a word $r(A,B)$ in the alphabet $ \{A,B,A^{-1},B^{-1}\}$.
Since $A,B\in U_{1}$, $A^{-1}=\overline{A}$ and
$B^{-1}=\overline{B}$. Then because
\begin{eqnarray*}
A &=&x+A^{-} \\
\overline{A} &=&x-A^{-} \\
B &=&x+B^{-} \\
\overline{B} &=&x-B^{-}
\end{eqnarray*}
$r(A,B)$ is a finite product of the expressions $(x+A^{-})$,
$(x-A^{-})$, $(x+B^{-})$, $(x-B^{-})$. This is a $k$-linear
combination of words $ s(A^{-},B^{-})$ in the alphabet
$\{A^{-},B^{-}\}$. Using the formulas
\begin{equation*}
2(A^{-}B^{-})^{+}=A^{-}B^{-}+B^{-}A^{-}\in k
\end{equation*}
\begin{equation*}
A^{-}A^{-}=B^{-}B^{-}=x^{2}-1\in k
\end{equation*}
it is possible to write $s(A^{-},B^{-})$ as a $k$-linear combination
of $\left\{ 1,A^{-},B^{-},A^{-}B^{-}\right\} $. Finally, note that $
A^{-}B^{-}=(A^{-}B^{-})^{+}+(A^{-}B^{-})^{-}.$
\end{proof}

\begin{definition}
Let $A,B\in U_{1}$ be conjugate elements in $U_{1}$, the unit
quaternions of the quaternion algebra $H=\frac{(\mu ,\nu )}{k}$. We
say that $(A,B)$ is an \emph{irreducible} pair iff $A$, $B$
generates $H$ as a $k$-algebra. Otherwise we say that $(A,B)$ is
\emph{reducible}. Thus $(A,B)$ is irreducible iff
\begin{equation*}
\mathcal{B}=\left\{ 1,A^{-},B^{-},(A^{-}B^{-})^{-}\right\}
\end{equation*}%
is a basis of the $k$-vector space $H$. We will say that $(A,B)$ is \emph{%
almost-irreducible} iff $(A,B)$ is reducible but $A^{-}$ and $B^{-}$ are
linearly independent.
\end{definition}

\section{c-Representations in $U_{1}$\label{srepresentationinU}}

Let $G$ be a group given by the presentation
\begin{equation*}
G=\left\vert a,b:w(a,b)\right\vert
\end{equation*}
where $w$ is a word in $a$ and $b$.

A homomorphism
\begin{equation*}
\rho :G\longrightarrow U_{1}
\end{equation*}
such that $\rho (a)$ and $\rho (b)$ are conjugate elements in $U_{1}$ is
called here a \textit{c-representation}. Set
\begin{equation*}
\rho (a)=A=A^{+}+A^{-},\quad \rho (b)=B=B^{+}+B^{-}.
\end{equation*}
Then, by Lemma \ref{lema1}, $A^{+}=B^{+}$ since $A$ and $B$ are conjugate elements. Set
\begin{equation*}
x=A^{+}=B^{+}
\end{equation*}
and
\begin{equation*}
y=-(A^{-}B^{-})^{+}.
\end{equation*}

\begin{definition}
We will say that a $c$-representation $\rho :G\longrightarrow U_{1}$
of a two generator group $G=\left\vert a,b:w(a,b)\right\vert$
\emph{realizes} the point $(x,y)$ if $x=\rho (a)^{+}=\rho (b)^{+}$
and $y=-(\rho (a)^{-}\rho (b)^{-})^{+}$.
\end{definition}

By Proposition \ref{psubalgebraAB} the subalgebra of $H$ generated
(as a $k$-algebra) by $im(\rho )$ coincides with the $k$-vector
subspace of $H$ generated (as a $k$-vector space) by
\begin{equation*}
\left\{ 1,A^{-},B^{-},(A^{-}B^{-})^{-}\right\}
\end{equation*}

We say that $\rho $ is \textit{irreducible} iff $im(\rho )$
generates $H$ as a $k$-algebra. Otherwise we say that $\rho $ is
\textit{reducible}. Thus $\rho $ is irreducible iff $(A,B)$ is an
irreducible pair iff
\begin{equation*}
\mathcal{B}=\left\{ 1,A^{-},B^{-},(A^{-}B^{-})^{-}\right\}
\end{equation*}%
is a basis of the $k$-vector space $H$.

We will say that $\rho $ is \textit{almost-irreducible} iff $(A,B)$
is an almost-irreducible pair; that is: $\rho $ is reducible but
$A^{-}$ and $ B^{-} $ are linearly independent.

\begin{proposition}
\label{pmetabelianas} A c-representation
\begin{equation*}
\rho :G=\left\vert a,b:w(a,b)\right\vert \longrightarrow U_{1}
\end{equation*}
is \textit{almost-irreducible} iff $\rho $ is a reducible metabelian
c-representation.
\end{proposition}

\begin{proof}
Recall that a representation is called reducible metabelian if the image $%
\rho (G)$ is reducible but it is not an abelian subgroup in $U_{1}$.
See \cite{GM1993}.  On the other hand it is clear that the subgroup
generate by  $A$ and  $B$ is abelian if and only if  the two pure
quaternions $A^{-}$ and $B^{-}$ commute, and this occurs if and only
if $A^{-}=\pm B^{-}$.
\end{proof}

\begin{theorem}
\label{Teorema uno} Let $G$ be a group given by the presentation
\begin{equation*}
G=\left\vert a,b:w(a,b)\right\vert
\end{equation*}%
where $w$ is a word in $a$ and $b$ and let $U_{1}$ denote the group
of unit quaternions of the algebra $\left(
\frac{-1,1}{\mathbb{C}}\right)$. Then there is an algorithm to
construct an ideal $\mathcal{I}_{G}^{c}$ generated by four
polynomials
\begin{equation*}
\left\{ p_{1}(x,y),p_{2}(x,y),p_{3}(x,y),p_{4}(x,y)\right\}
\end{equation*}%
with integer coefficients that can be characterized as follows:

\begin{enumerate}
\item The set of points $\{(x,y)\in V(\mathcal{I}_{G}^{c}):y^{2}\neq
(1-x^{2})^{2}\}$ coincides with the pairs $(x,y)$ for which there exists an
irreducible c-representation $\rho :G\longrightarrow U_{1}$, unique up to
conjugation in $U_{1}$, realizing $(x,y)$.

\item The set of points $\{(x,y)\in V(\mathcal{I}%
_{G}^{c}):y^{2}=(1-x^{2})^{2},\quad x^{2}\neq 1\}$ coincides with the pairs $%
(x,y)$ for which there exists an almost- irreducible c-representation $\rho
:G\longrightarrow U_{1}$, unique up to conjugation in $U_{1}$, realizing $%
(x,y)$.

\item The set of points $\{(x,y)\in V(\mathcal{I}_{G}^{c}):y=0,\quad
x^{2}=1\}$ coincides with the pairs $(x,y)$ for which neither irreducible nor
almost-irreducible c-representation $\rho :G\longrightarrow U_{1}$ realizing
$(x,y)$ exist.
\end{enumerate}
\end{theorem}

\begin{proof}
Let $\rho :G\longrightarrow U_{1}$ be an irreducible
c-representation. Set $\rho (a)=A$ and $\rho (b)=B$ and
\begin{equation*}
\fbox{$
\begin{array}{c}
x=A^{+}=B^{+} \\
y=-(A^{-}B^{-})^{+}
\end{array}
$}
\end{equation*}
and consider the vector basis
\begin{equation*}
\mathcal{B}=\left\{ 1,A^{-},B^{-},(A^{-}B^{-})^{-}\right\} .
\end{equation*}
Then
\begin{eqnarray*}
1 &=&AA^{-1}=A\overline{A}=(A^{+}+A^{-})(A^{+}-A^{-})=x^{2}-A^{-}A^{-} \\
&\Longrightarrow &\quad 1-x^{2}=-A^{-}A^{-}=-B^{-}B^{-}
\end{eqnarray*}
Set:
\begin{equation*}
\fbox{$u=1-x^{2}$}
\end{equation*}
The left product by $X\in H$ can be considered as a linear
automorphism of the 4-dimensional vector space $H$.
\begin{equation*}
\begin{array}{cccc}
\mathbb{X}: & H & \longrightarrow & H \\
& C & \rightarrow & XC
\end{array}
\end{equation*}
We want to write the matrix of $\mathbb{X}$ with respect to the
basis
\begin{equation*}
\mathcal{B}=\left\{ 1,A^{-},B^{-},(A^{-}B^{-})^{-}\right\}
\end{equation*}
belonging to the $k$-vector space structure of the quaternion
algebra $H$. Denote by $I_{4\times 4}$ the identity matrix $4\time
4$, and by $m(X)$ the matrix of the automorphism $\mathbb{X}$ with
respect to $\mathcal{B}$. Observe that
\begin{eqnarray*}
m(A) &=&xI_{4\times 4}+m(A^{-}) \\
m(\overline{A}) &=&xI_{4\times 4}-m(A^{-}) \\
m(B) &=&xI_{4\times 4}+m(B^{-}) \\
m(\overline{B}) &=&xI_{4\times 4}-m(B^{-})
\end{eqnarray*}
Computation of $m(A^{-})$:
\begin{eqnarray*}
A^{-}A^{-} &=&-u \\
A^{-}B^{-} &=&(A^{-}B^{-})^{+}+(A^{-}B^{-})^{-}=-y+(A^{-}B^{-})^{-}
\end{eqnarray*}
From the last equation
\begin{equation}
(A^{-}B^{-})^{-}=y+A^{-}B^{-}  \label{ambmm}
\end{equation}
Then%
\begin{equation*}
A^{-}(A^{-}B^{-})^{-}=A^{-}(y+A^{-}B^{-})=yA^{-}+A^{-}A^{-}B^{-}=yA^{-}-uB^{-}
\end{equation*}
\begin{equation}
\fbox{$m(A^{-})=
\begin{pmatrix}
0 & -u & -y & 0 \\
1 & 0 & 0 & y \\
0 & 0 & 0 & -u \\
0 & 0 & 1 & 0
\end{pmatrix}
$}  \label{mam}
\end{equation}
Computation of $m(B^{-})$: Recall that
$2(A^{-}B^{-})^{+}=A^{-}B^{-}+B^{-}A^{-}$, then
\begin{eqnarray*}
B^{-}A^{-} &=&2(A^{-}B^{-})^{+}-A^{-}B^{-}=-2y-(-y+(A^{-}B^{-})^{-}) \\
&=&-y-(A^{-}B^{-})^{-} \\
B^{-}B^{-} &=&-u \\
B^{-}(A^{-}B^{-})^{-} &=&B^{-}(y+A^{-}B^{-})=yB^{-}+B^{-}A^{-}B^{-} \\
&=&yB^{-}+(2(A^{-}B^{-})^{+}-A^{-}B^{-})B^{-} \\
&=&yB^{-}-2yB^{-}+uA^{-}=-yB^{-}+uA^{-}
\end{eqnarray*}
\begin{equation}
\fbox{$m(B^{-})=
\begin{pmatrix}
0 & -y & -u & 0 \\
0 & 0 & 0 & u \\
1 & 0 & 0 & -y \\
0 & -1 & 0 & 0
\end{pmatrix}
$}  \label{mbm}
\end{equation}
Computation of $m((A^{-}B^{-})^{-})$:
\begin{eqnarray*}
(A^{-}B^{-})^{-}A^{-}
&=&(y+A^{-}B^{-})A^{-}=(y+2(A^{-}B^{-})^{+}-B^{-}A^{-})A^{-} \\
&=&-yA^{-}+uB^{-} \\
(A^{-}B^{-})^{-}B^{-} &=&(y+A^{-}B^{-})B^{-}=yB^{-}-uA^{-} \\
(A^{-}B^{-})^{-}(A^{-}B^{-})^{-} &=&(-y-B^{-}A^{-})(y+A^{-}B^{-}) \\
&=&-y^{2}-y(B^{-}A^{-}+A^{-}B^{-})-u^{2}=y^{2}-u^{2}
\end{eqnarray*}
\begin{equation}
\fbox{$m(A^{-}B^{-})^{-})=
\begin{pmatrix}
0 & 0 & 0 & y^{2}-u^{2} \\
0 & -y & -u & 0 \\
0 & u & y & 0 \\
1 & 0 & 0 & 0
\end{pmatrix}
$}  \label{mambmm}
\end{equation}
Observe that $\mathbb{X}(1)=X$. Therefore, if $m(X)$ is the matrix of the
linear map $\mathbb{X}$ then the coordinates of the quaternion $X$ with
respect to $\mathcal{B}$ is the 4-tuple
\begin{equation*}
m(X)
\begin{pmatrix}
1 \\
0 \\
0 \\
0
\end{pmatrix}
\end{equation*}
Next we use the relator of the presentation $\left\vert
a,b:w(a,b)\right\vert $ to obtain necessary conditions on $x,y$ for
$\rho :G\longrightarrow U_{1}$ be a homomorphism. Computing the
matrix $ w(m(A),m(B))$ and imposing the condition that it is the
matrix with respect to $\mathcal{B}$ of the "product by $1"$ linear
map $H\longrightarrow H$, a $ 4\times 1$ matrix
\begin{equation*}
w(m(A),m(B))
\begin{pmatrix}
1 \\
0 \\
0 \\
0
\end{pmatrix}
-
\begin{pmatrix}
1 \\
0 \\
0 \\
0
\end{pmatrix}
\end{equation*}
is obtained whose entries form a set $\left\{
p_{1}(x,y),p_{2}(x,y),p_{3}(x,y),p_{4}(x,y)\right\} $ of four
polynomials with integer coefficients such that $p_{i}(x,y)=0,\quad
i=\{1,2,3,4\}$. This proves the first part of the theorem.

It is clear that if given $x_{0}\in \mathbb{C}$, $y_{0}\in \mathbb{C}$ such
that $p_{i}(x_{0},y_{0})=0$ for $i=\{1,2,3,4\},$ we can find quaternions $A$
and $B$ in $U_{1}\subset \left( \frac{-1,1}{\mathbb{C}}\right) $ such that

\begin{equation*}
\fbox{$
\begin{array}{c}
x_{0}=A^{+}=B^{+} \\
y_{0}=-(A^{-}B^{-})^{+}%
\end{array}
$}
\end{equation*}
and such that
\begin{equation*}
\mathcal{B}=\left\{ 1,A^{-},B^{-},(A^{-}B^{-})^{-}\right\} .
\end{equation*}
is a basis, then $\rho :G\longrightarrow U_{1}$ such that
\begin{equation*}
\begin{array}{c}
\rho (a)=A=x_{0}+A^{-} \\
\rho (b)=B=x_{0}+B^{-}
\end{array}
\end{equation*}
is an irreducible c-representation because the entries of the left side of
\begin{equation*}
w(m(A),m(B))
\begin{pmatrix}
1 \\
0 \\
0 \\
0
\end{pmatrix}
-
\begin{pmatrix}
1 \\
0 \\
0 \\
0
\end{pmatrix}
=
\begin{pmatrix}
0 \\
0 \\
0 \\
0
\end{pmatrix}
\end{equation*}
are $p_{i}(x_{0},y_{0})=0$ for $i=\{1,2,3,4\}$.

Therefore, to continue with the proof, we will consider the algebra
isomorphism
\begin{equation*}
\left( \frac{-1,1}{\mathbb{C}}\right) \longrightarrow
M(2,\mathbb{C})
\end{equation*}
defined by
\begin{equation*}
i\mapsto I=\left(
\begin{array}{cc}
0 & 1 \\
-1 & 0
\end{array}
\right) ,\qquad j\mapsto J=\left(
\begin{array}{cc}
0 & 1 \\
1 & 0
\end{array}
\right) ,\qquad ij\mapsto IJ=\left(
\begin{array}{cc}
1 & 0 \\
0 & -1
\end{array}
\right)
\end{equation*}
sending $U_{1}$ isomorphically onto $SL(2,\mathbb{C})$.

Let $(x_{0},y_{0})\in \mathbb{C\times C}$ be such that
$p_{i}(x_{0},y_{0})=0$ for $i=\{1,2,3,4\}$. A homomorphism $\rho
:G\longrightarrow SL(2,\mathbb{C)}$ is given by the image of the
generators
\begin{equation*}
\rho (a)=A=
\begin{pmatrix}
\alpha & \beta \\
\gamma & \delta
\end{pmatrix}%
,\qquad \rho (b)=B=
\begin{pmatrix}
\lambda & \eta \\
\zeta & \xi
\end{pmatrix}
\end{equation*}
Then the condition $A^{+}=B^{+}=x_{0}$ implies that
\begin{equation*}
A^{-}=
\begin{pmatrix}
\alpha -x_{0} & \beta \\
\gamma & \delta -x_{0}
\end{pmatrix}
,\qquad B^{-}=
\begin{pmatrix}
\lambda -x_{0} & \eta \\
\zeta & \xi -x_{0}
\end{pmatrix}
\end{equation*}
The condition $y_{0}=-(A^{-}B^{-})^{+}$ and the following
computation
\begin{eqnarray*}
(A^{-}B^{-})^{+} &=&\frac{1}{2}tr(A^{-}B^{-})=\frac{1}{2}
tr((A-x_{0}I)(B-x_{0}I)) \\
&=&\frac{1}{2}tr(AB-x_{0}(A+B)+x_{0}^{2}I)=\frac{1}{2}
(tr(AB)-4x_{0}^{2}+2x_{0}^{2}) \\
&=&\frac{1}{2}tr(AB)-x_{0}^{2}
\end{eqnarray*}
implies
\begin{equation*}
-y_{0}=\frac{1}{2}tr(AB)-x_{0}^{2}\quad \Longrightarrow \quad
tr(AB)=2(x_{0}^{2}-y_{0})
\end{equation*}
If $x_{0}^{2}\neq 1$, the elements $A$ and $B$ are not parabolic.
Then up to conjugation in $SL(2,\mathbb{C)}$ we can assume:
\begin{eqnarray*}
A &=&
\begin{pmatrix}
x_{0}+\sqrt{x_{0}^{2}-1} & 0 \\
0 & x_{0}-\sqrt{x_{0}^{2}-1}
\end{pmatrix}
\\
B &=&
\begin{pmatrix}
\lambda & \eta \\
1 & \xi
\end{pmatrix}
,\quad \lambda +\xi =2x_{0},\quad \lambda \xi -\eta =1
\end{eqnarray*}
Therefore
\begin{equation*}
tr(AB)=2(x_{0}^{2}-y_{0})=(x_{0}+\sqrt{x_{0}^{2}-1})\lambda
+(x_{0}-\sqrt{ x_{0}^{2}-1})\xi
\end{equation*}
The three last equations determine unique values for $\lambda $,
$\eta $ and $\xi $, namely:
\begin{equation*}
\begin{array}{ccc}
\lambda =x_{0}-\frac{y_{0}}{\sqrt{x_{0}^{2}-1}} & \xi =x_{0}+\frac{y_{0}}{%
\sqrt{x_{0}^{2}-1}} & \eta =\frac{y_{0}^{2}-(1-x_{0}^{2})^{2}}{1-x_{0}^{2}}%
\end{array}
\end{equation*}
Under the isomorphism $\left( \frac{-1,1}{\mathbb{C}}\right)
\longrightarrow M(2, \mathbb{C})$ we then have:

\begin{equation}\label{eaybalmost1}
\fbox{$
\begin{array}{l}
A =x_{0}+\sqrt{x_{0}^{2}-1}IJ \\
B
=x_{0}+\frac{2-3x_{0}^{2}+x_{0}^{4}-y_{0}^{2}}{2x_{0}^{2}-2}I+\frac{
-x_{0}^{2}+x_{0}^{4}-y_{0}^{2}}{2x_{0}^{2}-2}J-\frac{y}{\sqrt{x_{0}^{2}-1}}IJ
\end{array}
$}
\end{equation}

A calculation gives:
\begin{equation*}
(A^{-}B^{-})^{-}=\frac{-x_{0}^{2}+x_{0}^{4}-y_{0}^{2}}{2\sqrt{x_{0}^{2}-1}}I+%
\frac{2-3x_{0}^{2}+x_{0}^{4}-y_{0}^{2}}{2\sqrt{x_{0}^{2}-1}}J
\end{equation*}
The determinant of the $3\times 3$ matrix with rows the coordinates
of $ A^{-} $, $B^{-}$,$(A^{-}B^{-})^{-}$ is, after a tedious
calculation:
\begin{equation*}
y_{0}^{2}-(1-x_{0}^{2})^{2}
\end{equation*}
Therefore, by Proposition \ref{psubalgebraAB}, $\left\{ 1,A^{-},B^{-},(A^{-}B^{-})^{-}\right\} $ is a
basis iff $y_{0}^{2}\neq (1-x_{0}^{2})^{2}$. Hence if
$(x_{0},y_{0})$ belongs to $ \{(x,y)\in
V(\mathcal{I}_{G}^{c}):y^{2}\neq (1-x^{2})^{2},x^{2}\neq 1\}$ there
exists an irreducible c-representation $\rho
_{(x_{0},y_{0})}:G\longrightarrow U_{1}$, unique up to conjugation
in $U_{1}$, realizing $(x_{0},y_{0})$. Namely the one defined by
$\rho _{(x_{0},y_{0})}(a)=A$, $\rho _{(x_{0},y_{0})}(b)=B$.

Now, if $x_{0}^{2}\neq 1$ and $y_{0}^{2}=(1-x_{0}^{2})^{2}$ then $\left\{
1,A^{-},B^{-},(A^{-}B^{-})^{-}\right\} $ is not a basis but the vectors
\begin{eqnarray*}
A^{-} &=&\sqrt{x_{0}^{2}-1}IJ \\
B^{-} &=&-\frac{1}{2}I+\frac{1}{2}J\pm \sqrt{x_{0}^{2}-1}IJ
\end{eqnarray*}
are linearly independent because the rank of the $2\times 3$ matrix
with rows the coordinates of $A^{-}$, $B^{-}$ has rank $2$. Then the
map $\rho _{(x_{0},y_{0})}:G\longrightarrow U_{1}$, defined by $\rho
(a)=A$, $\rho (b)=B$, is the limit of a sequence of irreducible
c-representations $\rho _{(x_{i},y_{i})}$ where $\{(x_{i},y_{i})\}$
is a sequence of points in the set $V(\mathcal{I}_{G}^{c})\backslash
V(y^{2}-(1-x^{2})^{2})$ converging to $ (x_{0},y_{0})$. By
continuity, $\rho _{(x_{0},y_{0})}:G\longrightarrow U_{1}$ is also a
c-representation. It is almost-irreducible, unique up to conjugation
in $U_{1}$ and realizes $(x_{0},y_{0})$.

Next we study the case $x_{0}^{2}=1$. There are two cases:

(i) $A\in SL(2,\mathbb{C})$ and $B\in SL(2,\mathbb{C})$ are
conjugate parabolic elements with different fixed-points. Then, up
to conjugation in $U_{1}$, we can assume
\begin{equation*}
A=\left(
\begin{array}{cc}
x_{0} & 1 \\
0 & x_{0}
\end{array}
\right) ,\quad B=\left(
\begin{array}{cc}
x_{0} & 0 \\
\lambda & x_{0}
\end{array}
\right) ,\quad \lambda \neq 0
\end{equation*}
Under the isomorphism $\left( \frac{-1,1}{\mathbb{C}}\right)
\longrightarrow M(2, \mathbb{C})$ we have
\begin{equation*}
\begin{array}{ccc}
A^{-}=\frac{1}{2}I+\frac{1}{2}J & B^{-}=\frac{\lambda
}{2}I+\frac{-\lambda }{ 2}J & A^{-}B^{-}=\frac{\lambda
}{2}+\frac{\lambda }{2}IJ
\end{array}
\end{equation*}
and since $y_{0}=-(A^{-}B^{-})^{+}$ we must have $\lambda =-2y_{0}$.
Therefore $y_{0}\neq 0$ and this implies that $\left\{
1,A^{-},B^{-},(A^{-}B^{-})^{-}\right\} $ is a basis.

(ii) $A\in SL(2,\mathbb{C})$ and $B\in SL(2,\mathbb{C})$ are conjugate
parabolic elements with the same fixed-point or both are
\begin{equation*}
\pm \left(
\begin{array}{cc}
1 & 0 \\
0 & 1
\end{array}
\right) .
\end{equation*}
Then, up to conjugation in $U_{1}$, we can assume
\begin{equation*}
A=\left(
\begin{array}{cc}
x_{0} & \lambda \\
0 & x_{0}
\end{array}
\right) ,B=\left(
\begin{array}{cc}
x_{0} & \eta \\
0 & x_{0}
\end{array}
\right)
\end{equation*}
or equivalently
\begin{equation*}
\begin{array}{ccc}
A^{-}=\frac{\lambda }{2}I+\frac{\lambda }{2}J & B^{-}=\frac{\eta
}{2}I+\frac{ \eta }{2}J & A^{-}B^{-}=0
\end{array}
\end{equation*}
Then $y_{0}=-(A^{-}B^{-})^{+}=0$ and moreover $A^{-}$, $B^{-}$ are linearly
dependent.

Thus if $(x_{0},y_{0})\in V(\mathcal{I}_{G}^{c})$, and
$x_{0}^{2}=1$, then either $y_{0}\neq 0$ and there is an irreducible
c-representation $\rho _{(x_{0},y_{0})}:G\longrightarrow U_{1}$,
unique up to conjugation in $U_{1}$, realizing $(x_{0},y_{0})$ (the
one defined by $\rho _{(x_{0},y_{0})}(a)=A$, $\rho
_{(x_{0},y_{0})}(b)=B$) or $y_{0}=0$. In this case none of the
c-representation $\rho :G\longrightarrow U_{1}$ realizing such
$(x_{0},y_{0}) $ are irreducible or almost-irreducible. This
concludes the proof of the theorem.
\end{proof}

\begin{remark}
The common solutions to $\left\{ p_{i}(x,y),i\in \left\{
1,2,3,4\right\} \right\} $ is a algebraic variety
$V(\mathcal{I}_{G}^{c})$, the variety of c-representations of $G$ in
$U_{1}\cong SL(2,\mathbb{C})\subset M(2,\mathbb{C })$. If $a$ and
$b$ are conjugate elements in $G$, every representation is a
c-representation, and the algebraic variety
$V(\mathcal{I}_{G}^{c})=V( \mathcal{I}_{G})$ defines the components
of the $SL(2,\mathbb{C})$-character variety of $G$ containing irreducible
representations. The $SL(2,\mathbb{C})$-character variety of a group has been already
studied in several paper, see for instance, \cite{Riley1984}, \cite{CS1983},
\cite{GM1993},and \cite{HLM1992d}, \cite{HLM1995a}, \cite{HLM2003}
for a knot group.
\end{remark}

\begin{remark}
In practise, to simplify calculations, it is useful to follow a
different approach to compute the polynomials defining the algebraic
variety $V( \mathcal{I}_{G}^{c})$. The idea is to write the relator
of the presentation of $G$ in a balanced way $w_{1}(a,b)=w_{2}(a,b)$
where the words $w_{1},w_{2}$ have more or less half the length of
$w$. Then the condition is
\begin{gather*}
\quad w_{1}(\rho (A),\rho (B)) =w_{2}(\rho (A),\rho (B)) \\
\Leftrightarrow \quad w_{1}(\rho (A),\rho (B))-w_{2}(\rho (A),\rho
(B))=0
\end{gather*}
In this case the 4 polynomials are obtained by setting
\begin{equation}
\fbox{$(w_{1}(m(A),m(B))-w_{2}(m(A),m(B)))
\begin{pmatrix}
1 \\
0 \\
0 \\
0
\end{pmatrix}
$}  \label{relacion}
\end{equation}
This is particularly useful when $G$ is the group of a $2$-bridge knot or
link.
\end{remark}

\begin{example}
\label{ejemplo uno}The trefoil knot $3_{1}$.
\begin{equation*}
G=\pi _{1}(S^{3}-3_{1})=\left\vert a,b:aba=bab\right\vert
\end{equation*}
Using this presentation of $G(3_{1})$ every representation is a
c-representation, because $a$ and $b$ are conjugate elements in $G$. The
computation with Mathematica using (\ref{relacion}) give the polynomials
\begin{equation*}
\begin{pmatrix}
0 \\
-1+2x^{2}-2y \\
1-2x^{2}+2y \\
0%
\end{pmatrix}
\end{equation*}
Therefore
\begin{equation*}
p(x,y)=2x^{2}-2y-1
\end{equation*}
defines the ideal $\mathcal{I}_{G}$. If $(x,y)$ corresponds to a
almost-irreducible c-representation then $y=\pm (1-x^{2})$ and
$x^{2}\neq 1$. The solutions are $x=\frac{\pm \sqrt{3}}{2}$,
$y=\frac{\pm 1}{4}$. They provide the almost-irreducible
c-representations
\begin{equation*}
\rho (a)=\left(
\begin{array}{cc}
\frac{\pm \sqrt{3}+\sqrt{-1}}{2} & 0 \\
0 & \frac{\pm \sqrt{3}-\sqrt{-1}}{2}
\end{array}
\right) ,\qquad \rho (b)=\left(
\begin{array}{cc}
\frac{\pm \sqrt{3}+\sqrt{-1}}{2} & 0 \\
1 & \frac{\pm \sqrt{3}-\sqrt{-1}}{2}
\end{array}
\right) .
\end{equation*}
But note that besides these almost-irreducible c-representation there might
be other reducible c-representations realizing the same values. In fact, for
the same point $(x,y)=($ $\frac{\sqrt{3}}{2},\frac{1}{4})$ the
c-representation
\begin{equation*}
\rho (a)=\left(
\begin{array}{cc}
\frac{\sqrt{3}+\sqrt{-1}}{2} & 0 \\
0 & \frac{\sqrt{3}-\sqrt{-1}}{2}
\end{array}
\right) ,\qquad \rho (b)=\left(
\begin{array}{cc}
\frac{\sqrt{3}+\sqrt{-1}}{2} & 0 \\
0 & \frac{\sqrt{3}-\sqrt{-1}}{2}
\end{array}
\right)
\end{equation*}
is reducible (not almost-reducible).
\end{example}

\begin{remark}
If
\begin{eqnarray*}
\rho (a) &=&A=x_{0}+A^{-} \\
\rho (b) &=&B=x_{0}+B^{-} \\
y_{0} &=&-(A^{-}B^{-})^{+}
\end{eqnarray*}
then $tr(AB)=-2y_{0}+2x_{0}^{2}$. This gives us the formulas to
write the algebraic variety $V\left( \mathcal{I}_{G}^{c}\right) $ in
terms of the variables $tr(A^{2})$ and $tr(AB)$ as in the papers
\cite{HLM1992d} and \cite{HLM1995a}. Lets call $x^{\prime }$ the
variable $x$ used in those papers to avoid mistakes. The variables
in \cite{HLM1992d} and \cite{HLM1995a} are
\begin{equation*}
x^{\prime }=tr(A^{2})=\left( tr(A)\right) ^{2}-2\quad \text{and\quad
} z=tr(AB)
\end{equation*}
Therefore the change from $(x^{\prime },z)$ to our variables $(x,y)$ in the
present work is
\begin{eqnarray*}
x^{\prime } &=&4x^{2}-2 \\
z &=&2x^{2}-2y
\end{eqnarray*}
\end{remark}

In the above example, the polynomial
\begin{equation*}
p(x,y)=2x^{2}-2y-1
\end{equation*}%
in the variables $x^{\prime }$ and $z$ is $p(x^{\prime },z)=z-1$.
\bigskip

The above Theorem \ref{Teorema uno} can be sharpened in case
$x_{0}$, $ y_{0}\in \mathbb{R}$:
\begin{theorem}
\label{Teorema dos}Let $G$ be a group given by the presentation%
\begin{equation*}
G=\left\vert a,b:w(a,b)\right\vert
\end{equation*}
where $w$ is a word in $a$ and $b$. Let $\mathcal{I}_{G}^{c}$ be the ideal
generated by the four polynomials
\begin{equation*}
\left\{ p_{1}(x,y),p_{2}(x,y),p_{3}(x,y),p_{4}(x,y)\right\}
\end{equation*}
with integer coefficients. If $(x_{0}$, $y_{0})$ is a real point of the
algebraic variety $V(\mathcal{I}_{G}^{c})$ we distinguish two cases:

\begin{enumerate}
\item If%
\begin{equation*}
\begin{array}{c}
1-x_{0}^{2}>0, \\
(1-x_{0}^{2})^{2}>y_{0}^{2}
\end{array}
\end{equation*}
there exists an irreducible c-representation $\rho :G\longrightarrow U_{1}$,
$U_{1}=S^{3}\subset \mathbb{H}$,~unique up to conjugation in $U_{1}=S^{3}$,
realizing $(x_{0},y_{0})$. Namely:
\begin{equation*}
\fbox{$
\begin{array}{l}
A=x_{0}+\frac{1}{+\sqrt{1-x_{0}^{2}}}\left( y_{0}i+\sqrt{
(1-x_{0}^{2})^{2}-y_{0}^{2}}j\right) \\
B=x_{0}+\sqrt{1-x_{0}^{2}}i,\quad \sqrt{1-x_{0}^{2}}>0
\end{array}
$}
\end{equation*}
\begin{equation*}
(A^{-}B^{-})^{-}=-\sqrt{(1-x^{2})^{2}-y^{2}}ij
\end{equation*}
and $\rho (G)\subset \left( \frac{-1,-1}{k}\right) $ where
\begin{equation*}
k=Q(x_{0},y_{0},\sqrt{1-x_{0}^{2}},\sqrt{(1-x_{0}^{2})^{2}-y_{0}^{2}}).
\end{equation*}

\item The remaining cases. Then excepting the case
\begin{equation*}
\begin{array}{c}
1-x_{0}^{2}>0, \\
(1-x_{0}^{2})^{2}=y_{0}^{2}%
\end{array}%
\end{equation*}
and the case
\begin{equation*}
x_{0}^{2}=1,y_{0}=0
\end{equation*}
there exists an irreducible (or almost-irreducible) c-representation $\rho
:G\longrightarrow U_{1}$, $U_{1}=SL(2,\mathbb{R})\subset \left( \frac{-1,1}{%
\mathbb{R}}\right) $ realizing $(x_{0},y_{0})$. Moreover two such
homomorphisms are equal up to conjugation in $U_{\pm 1}$. Specifically:

\item[(2.1)] If
\begin{equation*}
\begin{array}{c}
1-x_{0}^{2}>0, \\
(1-x_{0}^{2})^{2}<y_{0}^{2}%
\end{array}%
\end{equation*}%
set
\begin{equation*}
\fbox{$
\begin{array}{c}
A=x_{0}+\sqrt{1-x_{0}^{2}}I,\quad \sqrt{1-x_{0}^{2}}>0 \\
B=x_{0}+\frac{1}{+\sqrt{1-x_{0}^{2}}}\left( y_{0}I+\sqrt{
y_{0}^{2}-(1-x_{0}^{2})^{2}}J\right)
\end{array}
$}
\end{equation*}
Then $\rho :G\longrightarrow U_{1}$ is irreducible, $\rho (G)\subset
\left( \frac{-1,1}{k}\right) $ where
\begin{equation*}
k=Q(x_{0},y_{0},\sqrt{1-x_{0}^{2}},\sqrt{y_{0}^{2}-(1-x_{0}^{2})^{2}})
\end{equation*}
and
\begin{equation*}
(A^{-}B^{-})^{-}=+\sqrt{y_{0}^{2}-(1-x_{0}^{2})^{2}}IJ
\end{equation*}

\item[(2.2)] If
\begin{equation*}
\begin{array}{c}
1-x_{0}^{2}<0, \\
(1-x_{0}^{2})^{2}<y_{0}^{2}%
\end{array}
\end{equation*}
there are two subcases:

\begin{enumerate}
\item[(2.2.1)] $y_{0}<0$. Set
\begin{equation*}
\fbox{$
\begin{array}{c}
A=x_{0}+\sqrt{x_{0}^{2}-1}J,\quad \sqrt{x_{0}^{2}-1}>0 \\
B=x_{0}+\frac{1}{+\sqrt{x_{0}^{2}-1}}\left(
\sqrt{y_{0}^{2}-(x_{0}^{2}-1)^{2} }I-y_{0}J\right)
\end{array}
$}
\end{equation*}
Then $\rho :G\longrightarrow U_{1}$ is irreducible and
\begin{equation*}
(A^{-}B^{-})^{-}=-\sqrt{y_{0}^{2}-(x_{0}^{2}-1)^{2}}IJ
\end{equation*}

\item[(2.2.2)] $y_{0}>0$. Set
\begin{equation*}
\fbox{$
\begin{array}{c}
A=x_{0}+\sqrt{x_{0}^{2}-1}J,\quad \sqrt{x_{0}^{2}-1}>0 \\
B=x_{0}-\frac{1}{+\sqrt{x_{0}^{2}-1}}\left( \sqrt{y_{0}^{2}-(x_{0}^{2}-1)^{2}%
}I+y_{0}J\right)
\end{array}
$}
\end{equation*}
Then $\rho :G\longrightarrow U_{1}$ is irreducible and
\begin{equation*}
(A^{-}B^{-})^{-}=+\sqrt{y_{0}^{2}-(x_{0}^{2}-1)^{2}}IJ
\end{equation*}%
In both cases $\rho (G)\subset \left( \frac{-1,1}{k}\right) $ where
\begin{equation*}
k=Q(x_{0},y_{0},\sqrt{x_{0}^{2}-1},\sqrt{y_{0}^{2}-(x_{0}^{2}-1)^{2}}).
\end{equation*}
\end{enumerate}

\item[(2.3)] If
\begin{equation*}
\begin{array}{c}
1-x_{0}^{2}<0, \\
(1-x_{0}^{2})^{2}>y_{0}^{2}
\end{array}
\end{equation*}
set
\begin{equation*}
\fbox{$
\begin{array}{c}
A=x_{0}+\sqrt{x_{0}^{2}-1}J,\quad \sqrt{x_{0}^{2}-1}>0 \\
B=x_{0}+\frac{1}{+\sqrt{x_{0}^{2}-1}}\left( -y_{0}J+\sqrt{
(x_{0}^{2}-1)^{2}-y_{0}^{2}}IJ\right)%
\end{array}
$}
\end{equation*}
Then $\rho :G\longrightarrow U_{1}$ is irreducible, $\rho (G)\subset
\left( \frac{-1,1}{k}\right) $ where
\begin{equation*}
k=Q(x_{0},y_{0},\sqrt{x_{0}^{2}-1},\sqrt{(x_{0}^{2}-1)^{2}-y_{0}^{2}})
\end{equation*}
and
\begin{equation*}
(A^{-}B^{-})^{-}=-\sqrt{(x_{0}^{2}-1)^{2}-y_{0}^{2}}IJ
\end{equation*}

\item[(2.4)] If
\begin{equation*}
\begin{array}{c}
1-x_{0}^{2}<0, \\
(1-x_{0}^{2})^{2}=y_{0}^{2}%
\end{array}
\end{equation*}
set
\begin{equation*}
\fbox{$
\begin{array}{c}
A=x_{0}+\sqrt{x_{0}^{2}-1}J,\quad \sqrt{x_{0}^{2}-1}>0 \\
B=x_{0}+(IJ+I)-\frac{_{y_{0}}}{\sqrt{x_{0}^{2}-1}}J
\end{array}
$}
\end{equation*}
Then $\rho :G\longrightarrow U_{1}$ is almost-irreducible, $\rho
(G)\subset \left( \frac{-1,1}{k}\right) $ where
\begin{equation*}
k=Q(x_{0},y_{0},\sqrt{x_{0}^{2}-1})
\end{equation*}
and
\begin{equation*}
(A^{-}B^{-})^{-}=-\sqrt{(x_{0}{}^{2}-1)}(I+IJ).
\end{equation*}

\item[(2.5)] If
\begin{equation*}
\begin{array}{c}
1-x_{0}^{2}=0, \\
y_{0}\neq 0
\end{array}
\end{equation*}
set%
\begin{equation*}
\fbox{$
\begin{array}{c}
A=x_{0}+I+J \\
B=x_{0}+\frac{y_{0}}{2}(I-J)
\end{array}
$}
\end{equation*}
\begin{equation*}
(A^{-}B^{-})^{-}=-y_{0}IJ
\end{equation*}
Then $\rho :G\longrightarrow U_{1}$ is irreducible and $\rho
(G)\subset \left( \frac{-1,1}{k}\right) $ where $k=Q(y_{0})$.
\end{enumerate}
\end{theorem}

\begin{proof}
As in the proof of the Theorem \ref{Teorema uno} we only need to find
quaternions $A$ and $B$ in the appropriate $U_{1}$ such that

\begin{equation*}
\fbox{$
\begin{array}{c}
x_{0}=A^{+}=B^{+} \\
y_{0}=-(A^{-}B^{-})^{+}
\end{array}
$}
\end{equation*}
and such that either $\left\{ 1,A^{-},B^{-},(A^{-}B^{-})^{-}\right\}
$ is a basis (case irreducible) or that $A^{-}$ and $B^{-}$ are
linearly independent (case almost-irreducible), because then $\rho
:G\longrightarrow U_{1}$ such that
\begin{equation*}
\begin{array}{c}
\rho (a)=A=x_{0}+A^{-} \\
\rho (b)=B=x_{0}+B^{-}
\end{array}
\end{equation*}
will be the desired homomorphism.

In case $1$, the desired quaternions (unique, up to conjugation in $U_{1}$)
are provided by Theorem \ref{torema1A}.

In case $2$ apply Theorem \ref{teorema2A}.
\end{proof}

\begin{figure}[ht]
\epsfig{file=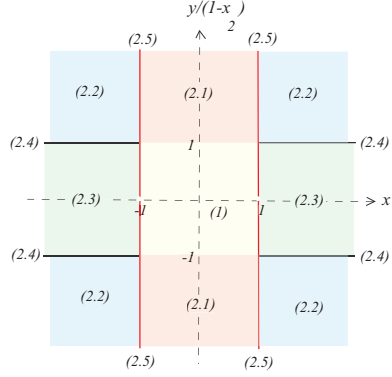,height=5cm}\caption{The pattern for real
points.}\label{fregiones}
\end{figure}

\begin{remark}\label{patron} To apply Theorem \ref{Teorema dos} it is useful to consider the pattern of Figure \ref{fregiones}.
 It shows the real plane with coordinates $x$
and $\frac{y}{1-x^{2}}$ . The plane is divided in several labeled regions by
labeled segments. The label corresponds to the case described in Theorem 4.
Therefore, to apply Theorem 4 to the algebraic variety $V(\mathcal{I}_{G})$
it is enough to study the graphic of $\frac{y}{1-x^{2}}$ as a funtion of $x$
over the pattern.
\end{remark}

The following remark concerns the points in the unlabeled
segments in the pattern between region $(1)$ and region $(2.1)$.
\begin{remark}
In the case that $(x_{0}$, $y_{0})$ is a real point of the algebraic variety
$V(\mathcal{I}_{G}^{c})$ and
\begin{equation*}
\begin{array}{c}
1-x_{0}^{2}>0, \\
(1-x_{0}^{2})^{2}=y_{0}^{2}%
\end{array}%
\end{equation*}%
then, according to Theorem \ref{Teorema uno} there is an almost-irreducible
c-representation $\rho :G\longrightarrow SL(2,\mathbb{C})$ realizing $(x_{0}$
, $y_{0})$ but no almost-irreducible c-representations $\rho
:G\longrightarrow SL(2,\mathbb{R})$ or $\rho :G\longrightarrow S^{3}$
realizing $(x_{0}$, $y_{0})$ exists, because the cases $(1)$ and $(2.1)$ in
Theorem \ref{Teorema dos}  both contain the limiting case
\begin{equation*}
\begin{array}{c}
1-x_{0}^{2}>0, \\
(1-x_{0}^{2})^{2}=y_{0}^{2}
\end{array}%
\end{equation*}%
and then $A^{-}$ and $B^{-}$ coincide, so that $\rho $ is not
almost-irreducible. An example is given above, where the almost-irreducible
c-representation $\rho :\pi _{1}(S^{3}-3_{1})\rightarrow SL(2,\mathbb{C})$
granted by Theorem \ref{Teorema uno} and corresponding to $x_{0}=\frac{\sqrt{
3}}{2},$ $y_{0}=\frac{1}{4}$ is
\begin{equation*}
\rho (a)=\left(
\begin{array}{cc}
\frac{\sqrt{3}+\sqrt{-1}}{2} & 0 \\
0 & \frac{\sqrt{3}-\sqrt{-1}}{2}
\end{array}%
\right) ,\qquad \rho (b)=\left(
\begin{array}{cc}
\frac{\sqrt{3}+\sqrt{-1}}{2} & 0 \\
1 & \frac{\sqrt{3}-\sqrt{-1}}{2}
\end{array}%
\right)
\end{equation*}%
However, since $1-x_{0}^{2}=\frac{1}{4}>0$ and $(1-x_{0}^{2})^{2}=\frac{1}{16
}=y_{0}^{2}$ almost-irreducible c-representations $\rho :\pi
_{1}(S^{3}-3_{1})\longrightarrow SL(2,\mathbb{R})$ or $\rho :\pi
_{1}(S^{3}-3_{1})\longrightarrow S^{3}$ realizing $(x_{0},y_{0})=($ $\frac{
\sqrt{3}}{2},\frac{1}{4})$ do not exist. This is in striking contrast with
the case
\begin{equation*}
\begin{array}{c}
1-x_{0}^{2}<0, \\
(1-x_{0}^{2})^{2}=y_{0}^{2}
\end{array}
\end{equation*}
($(x_{0}$, $y_{0})$ a real point) in which an almost-irreducible
c-representation $\rho :G\longrightarrow SL(2,\mathbb{R})$ realizing $
(x_{0},y_{0})$ does in fact exist (case $(2.4)$ of Theorem \ref{Teorema dos}).
An example follows.
\end{remark}

\begin{example} \label{ejemplo dos}
The figure-eight knot $4_{1}$.
\begin{equation*}
G(4_{1})=\pi _{1}(S^{3}-4_{1})=\left\vert a,b:aw=wb\right\vert ,\,\qquad
w=ba^{-1}b^{-1}a
\end{equation*}%
Using this presentation of $G(4_{1})$ every representation is a
c-representation, because $a$ and $b$ are conjugate elements in
$G(4_{1})$. The computation with Mathematica gives the following
polynomials defining $ \mathcal{I}_{G(4_{1})}$:
\begin{equation*}
\begin{pmatrix}
0 \\
1-6x^{2}+4x^{4}-2y-4y^{2} \\
1-6x^{2}+4x^{4}-2y-4y^{2} \\
0%
\end{pmatrix}
\end{equation*}
Therefore
\begin{equation*}
p(x,y)=1-6x^{2}+4x^{4}-2y-4y^{2}
\end{equation*}
defines the ideal $\mathcal{I}_{G(4_{1})}$. The points $(x_{0},y_{0})\in V(\mathcal{
I}_{G(4_{1})})$ with $y_{0}=\pm (1-x_{0}^{2})$ are $(\pm \frac{\sqrt{5}}{2},-\frac{1
}{4})$. Note that $1-x_{0}^{2}=-\frac{1}{4}<0$. Then, case $(2.4)$
of Theorem \ref{Teorema dos} applies and there exists the following
almost-irreducible c-representation (see \cite{GM1993})
\begin{equation*}
\rho :G(4_{1})\longrightarrow U_{1}\subset \left(
\frac{-1,1}{\mathbb{R}}\right)
\end{equation*}
realizing $(x_{0},y_{0})=(\frac{\sqrt{5}}{2},-\frac{1}{4})$:
\begin{equation*}
\rho (a)=A=\frac{\sqrt{5}}{2}+\frac{1}{2}J,
\end{equation*}
$\qquad $%
\begin{equation*}
\rho (b)=B=\frac{\sqrt{5}}{2}+I+\frac{1}{2}J+IJ
\end{equation*}
and after identifying $\left( \frac{-1,1}{\mathbb{R}}\right) $ with $M(2,
\mathbb{R})$ via
\begin{equation*}
a+bI+cJ+dIJ\longleftrightarrow \left(
\begin{array}{cc}
a+d & b+c \\
-b+c & a-d
\end{array}
\right)
\end{equation*}
we have $\rho :\pi _{1}(S^{3}-4_{1})\longrightarrow SL(2,\mathbb{R})$:
\begin{equation*}
\rho (a)=A=\left(
\begin{array}{cc}
1+\frac{\sqrt{5}}{2} & \frac{3}{2} \\
-\frac{1}{2} & -1+\frac{\sqrt{5}}{2}
\end{array}
\right)
\end{equation*}

\begin{equation*}
\rho (b)=B=\left(
\begin{array}{cc}
\frac{\sqrt{5}}{2} & \frac{1}{2} \\
\frac{1}{2} & \frac{\sqrt{5}}{2}
\end{array}
\right)
\end{equation*}
\end{example}

\subsection{Reducible c-representations}

As we saw before if the c-representation $\rho :G\longrightarrow U_{1}$ is
almost-irreducible then $y^{2}=(1-x^{2})^{2}$. If $\rho $ is reducible and $
B^{-}=\pm A^{-}$ this implies also $y^{2}=(1-x^{2})^{2}$. The algebraic
variety $V(y^{2}-(1-x^{2})^{2})$ is the union of two parabolas of equations $
y=1-x^{2}$ and $y=x^{2}-1$. Given a point $(x_{0},y_{0})\in
V(y^{2}-(1-x^{2})^{2})$ we want to know if there exists a reducible
c-representation $\rho :G\longrightarrow U_{1}$ realizing it.

\textbf{Case 1}. $(x_{0},y_{0})\in V(y+x^{2}-1)$. The putative reducible
c-representation realizing $(x_{0},y_{0})$ must be defined by
\begin{equation*}
\rho (a)=x+A^{-},\quad \rho (b)=x+B^{-},\quad A^{-}=B^{-}\quad \quad
\end{equation*}
because then
\begin{equation*}
y=-(A^{-}B^{-})^{+}=-(A^{-}A^{-})^{+}=-A^{-}A^{-}=1-x^{2}
\end{equation*}
since
\begin{equation*}
1=N(A)=A\overline{A}=(x+A^{-})(x-A^{-})=x^{2}-A^{-}A^{-}
\end{equation*}
A necessary and sufficient condition for the existence of such a
c-representation
\begin{equation*}
\rho :G=\left\vert a,b:w(a,b)\right\vert \longrightarrow U_{1}
\end{equation*}
is
\begin{equation*}
w(x+A^{-},x+A^{-})=1
\end{equation*}
A sufficient condition for this is that the group presented by:
\begin{equation*}
\left\vert a,b:w(a,b)=1;a=b\right\vert
\end{equation*}
be trivial. This happens for the standard presentations of the groups of
2-bridge knots and links ($\left\vert a,b:av(a,b)=v(a,b)b\right\vert $ for
2-bridge knots and $\left\vert a,b:av(a,b)=v(a,b)a\right\vert $ for 2-bridge
links).

We conclude that:
\begin{proposition}\label{prop9}
The points of the algebraic variety $V(y+x^{2}-1)$ can be realized by
reducible c-representations $\rho :G\longrightarrow U_{1}$ if $G$ is the
group of a 2-bridge knot or link.\qed
\end{proposition}

\textbf{Case 2}. $(x_{0},y_{0})\in V(y-x^{2}+1)$. The putative reducible
c-representation realizing $(x_{0},y_{0})$ must be defined by
\begin{equation*}
\rho (a)=x+A^{-},\quad \rho (b)=x+B^{-},\quad B^{-}=-A^{-}\quad \quad
\end{equation*}
because then%
\begin{equation*}
y=-(A^{-}B^{-})^{+}=(A^{-}A^{-})^{+}=A^{-}A^{-}=-1+x^{2}
\end{equation*}
since
\begin{equation*}
1=N(A)=A\overline{A}=(x+A^{-})(x-A^{-})=x^{2}-A^{-}A^{-}
\end{equation*}
A necessary and sufficient condition for the existence of such a
c-representation
\begin{equation*}
\rho :G=\left\vert a,b:w(a,b)\right\vert \longrightarrow U_{1}
\end{equation*}
is
\begin{equation*}
w(x+A^{-},x-A^{-})=1
\end{equation*}
A sufficient condition for this is that the group presented by:
\begin{equation*}
\left\vert a,b:w(a,b)=1;b=a^{-1}\right\vert
\end{equation*}
be trivial. This happens for the standard presentations $\left\vert
a,b:av(a,b)=v(a,b)a\right\vert $ of the group of a 2-bridge link.

For the standard presentations $\left\vert a,b:av(a,b)=v(a,b)b\right\vert $
of the group of a 2-bridge knot, the relation
\begin{equation*}
(x+A^{-})v(x+A^{-},x-A^{-})=v(x+A^{-},x-A^{-})(x-A^{-})
\end{equation*}
implies either
\begin{equation*}
(x+A^{-})(x+A^{-})^{n}=(x+A^{-})^{n}(x-A^{-}),\quad n\geq 0
\end{equation*}
or
\begin{equation*}
(x+A^{-})(x-A^{-})^{n}=(x-A^{-})^{n}(x-A^{-}),\quad n\geq 1
\end{equation*}
because $1=N(A)=A\overline{A}=(x+A^{-})(x-A^{-})$. Both cases imply $
(x+A^{-})=(x-A^{-})$. That is $A^{-}=0,$ which corresponds to the point $
(x_{0},y_{0})=(\pm 1,0)$.

We conclude that:
\begin{proposition}\label{prop10}
The points of the algebraic variety $V(y-x^{2}+1)$ can be realized by
reducible c-representations $\rho :G\longrightarrow U_{1}$ if $G$ is the
group of a 2-bridge link but not if $G$ is the group of a 2-bridge knot.\qed
\end{proposition}

\subsection{Two examples} Next we analyze the real points of the algebraic variety $V(\mathcal{I}_{G})$ for two knot groups.
\subsubsection{The Trefoil knot}
\begin{figure}[ht]
\epsfig{file=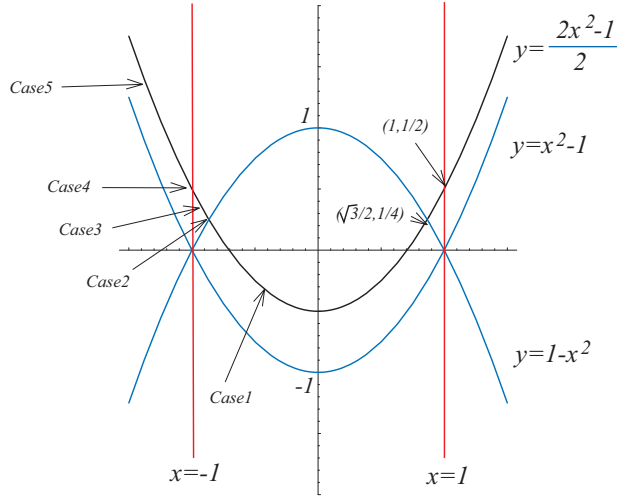,height=7cm}\caption{The real curve
$V(\mathcal{I} _{G(3_{1})})$ and the five cases.}\label{fctrebol}
\end{figure}

As was said in Example \ref{ejemplo uno}, the algebraic variety
$V(\mathcal{I}_{G(3_{1})})$ for the Trefoil knot group
\begin{equation*}
G(3/1)=|a,b;aba=bab|
\end{equation*}
is define by the ideal
 $\mathcal{I}_{G(3_{1})}=(2x^{2}-2y-1)$. The real part of  $V(\mathcal{I}_{G(3_{1})})$ is the parabola
 $y=\frac{2x^{2}-1}{2}$ depicted
in Figure \ref{fctrebol}.

Figure \ref{fctrebol} also exhibits the variety of reducible
representations $ V(y^{2}-(x^{2}-1)^{2})$. This is the union of the
two parabolas $y=\pm (x^{2}-1)$. The trefoil knot is a 2-bridge
knot. Then, by Propositions \ref{prop9} and \ref{prop10}, only the
points of the parabola $ y=1-x^{2}$ can be realized by reducible
c-representations. The intersection of $V(\mathcal{I}_{G(3_{1})})$
with $V(y-1+x^{2})$ define the two almost-irreducible
representations realizing the pairs $(x,y)=(\pm \frac{\sqrt{3}}{2},\frac{1}{4%
})$ described in Example \ref{ejemplo uno}. The remaining points of
$V(\mathcal{I}_{G(3_{1})})$ correspond to irreducible
c-representations.

Theorem \ref{Teorema dos} establishes the different cases of representations associated to real
points of the algebraic variety $V(\mathcal{I}_{G(3_{1})})$. To
apply Theorem \ref{Teorema dos} to the algebraic variety $V(\mathcal{I}_{G})$ it suffices to
study the graphic of $\frac{y}{1-x^{2}}$ as a function of $x$ over the
pattern of Remark \ref{patron}. See Figure \ref{fregiones}.

\begin{figure}[ht]
\epsfig{file=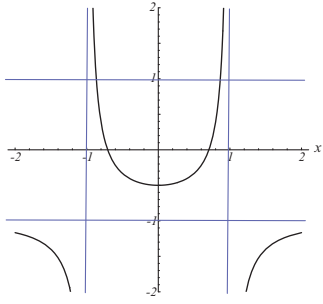,height=4cm}\caption{The function $\frac{y}{1-x^{2}}$
for the Trefoil knot.}\label{fysobreuutrebol}
\end{figure}

Figure \ref{fysobreuutrebol} shows $\frac{y}{1-x^{2}}$ as a function of $x$ for the algebraic variety $V(\mathcal{I}_{G(3_{1})})$.

Then, according with Theorem \ref{Teorema dos}, there are five
cases:

\begin{description}
\item[Case 1] Region (1). $x\in (\frac{-\sqrt{3}}{2},\frac{\sqrt{3}}{2}
)\Longleftrightarrow \left\{ 1-x^{2}>0,(1-x^{2})^{2}>y^{2}\right\}
$, where
\begin{eqnarray*}
y &=&\frac{2x^{2}-1}{2} \\
u &=&1-x^{2}
\end{eqnarray*}%
There exists an irreducible c-representation $\rho
_{x}:G(3/1)\longrightarrow S^{3}$ realizing $(x,y)$, unique up to
conjugation in $S^{3}$, such that%
\begin{equation}
\begin{array}{l}
\rho _{x}(a)=A=x+\frac{2x^{2}-1}{2\sqrt{1-x^{2}}}i+\frac{1}{2}\sqrt{\frac{
3-4x^{2}}{1-x^{2}}}j \\
\rho _{x}(b)=B=x+\sqrt{1-x^{2}}i,\quad \sqrt{1-x^{2}}>0
\end{array}
\label{ecase1}
\end{equation}
The composition of $\rho _{x}$ with $c:S^{3}\rightarrow SO(3)$, where $c(X)$
, $X\in S^{3}$, acts on $P\in H_{0}\cong E^{3}$ by conjugation,
defines the representation $\rho _{x}^{\prime }=c\circ \rho
_{x}:G(3/1)\longrightarrow SO(3)$. In linear notation, where
$\left\{ X,Y,Z\right\} $ is the coordinate system in $E^{3}$
associated to the basis $\left\{ -ij,j,i\right\} $ we have
\begin{equation*}
\begin{array}{l}
\rho _{x}^{\prime }(a)=m_{x}(a)\left(
\begin{array}{c}
X \\
Y \\
Z
\end{array}
\right) =\left(
\begin{array}{c}
X^{\prime } \\
Y^{\prime } \\
Z^{\prime }
\end{array}
\right)  \\
\rho _{x}^{\prime }(b)=m_{x}(b)\left(
\begin{array}{c}
X \\
Y \\
Z%
\end{array}
\right) =\left(
\begin{array}{c}
X^{\prime } \\
Y^{\prime } \\
Z^{\prime }
\end{array}
\right)
\end{array}
\end{equation*}
where, by (\ref{emtxyz}),
\begin{equation*}
m_{x}(a)=\left(
\begin{array}{ccc}
2x^{2}-1 & \frac{x-2x^{3}}{\sqrt{1-x^{2}}} & x\sqrt{\frac{4x^{2}-3}{x^{2}-1}}
\\
\frac{x(2x^{2}-1)}{\sqrt{1-x^{2}}} & \frac{-4x^{4}+2x^{2}+1}{2-2x^{2}} &
\frac{(1-2x^{2})\sqrt{3-4x^{2}}}{2(x^{2}-1)} \\
-x\sqrt{\frac{4x^{2}-3}{x^{2}-1}} & \frac{(1-2x^{2})\sqrt{3-4x^{2}}}{%
2(x^{2}-1)} & \frac{1-2x^{2}}{2x^{2}-2}
\end{array}%
\right)
\end{equation*}
and
\begin{equation*}
m_{x}(b)=%
\begin{pmatrix}
2x^{2}-1 & -2x\sqrt{1-x^{2}} & 0 \\
2x\sqrt{1-x^{2}} & 2x^{2}-1 & 0 \\
0 & 0 & 1%
\end{pmatrix}
.
\end{equation*}
The maps $\rho _{x}^{\prime }(a)$ and $\rho _{x}^{\prime }(b)$ are
right rotations of angle $\alpha $ around the axes $A^{-}$ and
$B^{-}$ where $ x=A^{+}=B^{+}=\cos \frac{\alpha }{2}$. See Figure
\ref{AmBm}. The angle $\omega $ between the axes of $\rho
_{x}^{\prime }(a)$ and $\rho _{x}^{\prime }(b)$ is given by
\begin{equation*}
\cos \omega =\frac{y}{u}=\frac{x^{2}-\frac{1}{2}}{1-x^{2}}
\end{equation*}

\item[Case 2]  $(x,y)=(\pm \frac{\sqrt{3}}{2},\frac{1}{4}
)\Longleftrightarrow \left\{ 1-x^{2}>0,(1-x^{2})^{2}=y^{2}\right\}
$.

There exists an almost-irreducible c-representation $\rho
_{x}:G(3/1)\longrightarrow U_{1}\subset \left(
\frac{-1,1}{\mathbb{C}} \right) $ realizing $(x,y)$, unique up to
conjugation in $U_{1}$ (Theorem 3, (\ref{eaybalmost1})), such that :
\begin{eqnarray}
\rho _{\pm \sqrt{3}/2}(a) &=&A=\frac{\pm \sqrt{3}}{2}+\frac{\sqrt{-1}}{2}IJ
\label{ecaso2} \\
\rho _{\pm \sqrt{3}/2}(b) &=&B=\frac{\pm \sqrt{3}}{2}-\frac{1}{2}I+\frac{1}{2%
}J+\frac{\sqrt{-1}}{2}IJ  \notag
\end{eqnarray}%
This representation cannot be conjugated to any real representation. Under
the isomorphism $U_{1}\approx SL(2,\mathbb{C})$ we have
\begin{equation*}
\begin{array}{cccc}
\rho _{\pm \sqrt{3}/2}: & G(3/1) & \longrightarrow & SL(2,\mathbb{C}) \\
& a & \rightarrow & A=
\begin{pmatrix}
\frac{\pm \sqrt{3}}{2}+\frac{\sqrt{-1}}{2} & 0 \\
0 & \frac{\pm \sqrt{3}}{2}-\frac{\sqrt{-1}}{2}
\end{pmatrix}
\\
& b & \rightarrow & B=
\begin{pmatrix}
\frac{\pm \sqrt{3}}{2}+\frac{\sqrt{-1}}{2} & 0 \\
1 & \frac{\pm \sqrt{3}}{2}-\frac{\sqrt{-1}}{2}
\end{pmatrix}
\end{array}
\end{equation*}

\item[Case 3] Region (2.1). $x\in (-1,\frac{-\sqrt{3}}{2})\cup (\frac{\sqrt{3}%
}{2},1)\Longleftrightarrow \left\{ 1-x^{2}>0,(1-x^{2})^{2}<y^{2}\right\} .$

There exists an irreducible c-representation $\rho
_{x}:G(3/1)\longrightarrow SL(2,\mathbb{R})=U_{1}\subset \left( \frac{-1,1}{%
\mathbb{R}}\right) $ realizing $(x,y)$, unique up to conjugation in $SL(2,%
\mathbb{R})$, such that%
\begin{equation*}
\begin{array}{l}
\rho _{x}(a)=A=x+\sqrt{1-x^{2}}I,\quad \sqrt{1-x^{2}}>0 \\
\rho _{x}(b)=B=x+\frac{2x^{2}-1}{2\sqrt{1-x^{2}}}I+\frac{1}{2}\sqrt{\frac{%
-3+4x^{2}}{1-x^{2}}}J%
\end{array}%
\end{equation*}%
The composition of $\rho _{x}$ with $c:SL(2,\mathbb{R})\rightarrow
SO^{+}(1,2)\cong Iso^{+}(\mathbb{H}^{2})$, where $c(X)$, $X\in SL(2,\mathbb{R%
})$, acts on $P\in H_{0}\cong E^{1,2}$ by conjugation, defines the
representation $\rho _{x}^{\prime }=c\circ \rho _{x}:G(3/1)\longrightarrow
SO^{+}(1,2)$ in affine linear notation, where $\left\{ X,Y,Z\right\} $ is
the coordinate system associated to the basis $\left\{ -IJ,J,I\right\} $
\begin{equation*}
\begin{array}{l}
\rho _{x}^{\prime }(a)=m_{x}(a)\left(
\begin{array}{c}
X \\
Y \\
Z
\end{array}
\right) =\left(
\begin{array}{c}
X^{\prime } \\
Y^{\prime } \\
Z^{\prime }
\end{array}
\right) \\
\rho _{x}^{\prime }(b)=m_{x}(b)\left(
\begin{array}{c}
X \\
Y \\
Z
\end{array}
\right) =\left(
\begin{array}{c}
X^{\prime } \\
Y^{\prime } \\
Z^{\prime }
\end{array}
\right)
\end{array}
\end{equation*}
where the matrices of $\rho _{x}^{\prime }(a)$ and $\rho
_{x}^{\prime }(b)$ are respectively:
\begin{equation*}
m_{x}(a)=\left(
\begin{array}{ccc}
2x^{2}-1 & -2x\sqrt{1-x^{2}} & 0 \\
2x\sqrt{1-x^{2}} & 2x^{2}-1 & 0 \\
0 & 0 & 1
\end{array}
\right)
\end{equation*}
and
\begin{equation*}
m_{x}(b)=
\begin{pmatrix}
2x^{2}-1 & \frac{x-2x^{3}}{\sqrt{1-x^{2}}} & x\sqrt{\frac{3-4x^{2}}{x^{2}-1}}
\\
\frac{x(2x^{2}-1)}{\sqrt{1-x^{2}}} & \frac{1+2x^{2}-4x^{4}}{2-2x^{2}} &
\frac{\sqrt{4x^{2}-3}(2x^{2}-1)}{2(1-x^{2})} \\
x\sqrt{\frac{3-4x^{2}}{x^{2}-1}} & \frac{\sqrt{4x^{2}-3}(2x^{2}-1)}{
2(1-x^{2})} & \frac{1-2x^{2}}{2x^{2}-2}
\end{pmatrix}
\end{equation*}
The maps $\rho _{x}^{\prime }(a)$ and $\rho _{x}^{\prime }(b)$ are
right (spherical) rotations of $H_{0}\cong E^{1,2}$ of angle $\alpha
$ around the time-like axes $A^{-}$ and $B^{-}$ where
$x=A^{+}=B^{+}=\cos \frac{\alpha }{2 }$. See Figure \ref{AmBm2}.

The distance $d$ (measured in the hyperbolic plane) between the axes
of $ \rho _{x}^{\prime }(a)$ and $\rho _{x}^{\prime }(b)$ is given
by
\begin{equation*}
\cosh d=\frac{y}{u}=\frac{x^{2}-\frac{1}{2}}{-x^{2}+1}
\end{equation*}

\item[Case 4] Segment (2.5). $(x,y)=(\pm 1,\frac{1}{2})\Longleftrightarrow
\left\{ 1-x^{2}=0,(1-x^{2})^{2}<y^{2}\right\} $

There exists an irreducible c-representation $\rho
_{x}:G(3/1)\longrightarrow SL(2,\mathbb{R})=U_{1}\subset \left( \frac{-1,1}{
\mathbb{R}}\right) $ realizing $(x,y)$, unique up to conjugation in $SL(2,
\mathbb{R})$, such that
\begin{equation*}
\begin{array}{l}
\rho _{x}(a)=A=\pm 1+I+J \\
\rho _{x}(b)=B=\pm 1+\frac{1}{4}(I-J)
\end{array}%
\end{equation*}
The composition of $\rho _{x}$ with $c:SL(2,\mathbb{R})\rightarrow
SO^{+}(1,2)\cong Iso^{+}(\mathbb{H}^{2})$, defines the representation $\rho
_{x}^{\prime }=c\circ \rho _{x}:G(3/1)\longrightarrow SO^{+}(1,2)$ such that
the matrices of $\rho _{x}^{\prime }(a)$ and $\rho _{x}^{\prime }(b)$ are
respectively:
\begin{equation*}
m_{x}(a)=m(-1,1;\pm 1,1,1,0)=\left(
\begin{array}{ccc}
1 & -2 & 2 \\
2 & -1 & 2 \\
2 & -2 & 3
\end{array}
\right)
\end{equation*}
and%
\begin{equation*}
m_{x}(b)=m(-1,1;\pm 1,\frac{1}{4},-\frac{1}{4},0)=%
\begin{pmatrix}
1 & -\frac{1}{2} & -\frac{1}{2} \\
\frac{1}{2} & \frac{7}{8} & -\frac{1}{8} \\
-\frac{1}{2} & \frac{1}{8} & \frac{9}{8}
\end{pmatrix}
\end{equation*}
where $\left\{ X,Y,Z\right\} $ is the coordinate system associated to the
basis $\left\{ -IJ,J,I\right\} $. The maps $\rho _{x}^{\prime }(a)$ and $%
\rho _{x}^{\prime }(b)$ are parabolic rotations on $H_{0}\cong
E^{1,2}$ around the nullvector axes $A^{-}$ and $B^{-}$. See Figure
\ref{ambm3}.

\item[Case 5] Region (2.2). $x\in (-\infty ,-1)\cup (1,\infty )\Longleftrightarrow
\left\{ 1-x^{2}<0,(1-x^{2})^{2}<y^{2}\right\} $ ( $y>0$).

There exists an irreducible c-representation $\rho
_{x}:G(3/1)\longrightarrow SL(2,\mathbb{R})=U_{1}\subset \left( \frac{-1,1}{\mathbb{R}}\right) $ realizing $(x,y)$, unique up to conjugation in $SL(2,
\mathbb{R})$, such that
\begin{equation*}
\begin{array}{l}
\rho _{x}(a)=A=x+\sqrt{x^{2}-1}J,\quad \sqrt{x^{2}-1}>0 \\
\rho _{x}(b)=B=x-\frac{1}{2}\sqrt{\frac{4x^{2}-3}{x^{2}-1}}I-\frac{2x^{2}-1}{%
2\sqrt{x^{2}-1}}J%
\end{array}
\end{equation*}
The composition of $\rho _{x}$ with $c:SL(2,\mathbb{R})\rightarrow
SO^{+}(1,2)\cong Iso^{+}(\mathbb{H}^{2})$, define the representations $\rho
_{x}^{\prime }=c\circ \rho _{x}:G(3/1)\longrightarrow SO^{+}(1,2)$ such that
the matrices of $\rho _{x}^{\prime }(a)$ and $\rho _{x}^{\prime }(b)$ are
respectively:

\begin{equation*}
m_{x}(a)=\left(
\begin{array}{cccc}
2x^{2}-1 & 0 & 2x\sqrt{x^{2}-1} & 0 \\
0 & 1 & 0 & \frac{\left( 3-4x^{2}\right) \sqrt{x^{2}-1}}{4x} \\
2x\sqrt{x^{2}-1} & 0 & 2x^{2}-1 & 0
\end{array}
\right)
\end{equation*}

\begin{equation*}
m_{x}(b)=
\begin{pmatrix}
2x^{2}-1 & x\sqrt{\frac{-3+4x^{2}}{x^{2}-1}} & -\frac{x(2x^{2}-1)}{\sqrt{%
x^{2}-1}} \\
-x\sqrt{\frac{-3+4x^{2}}{x^{2}-1}} & \frac{1-2x^{2}}{2x^{2}-2} & \frac{\sqrt{%
-3+4x^{2}}(2x^{2}-1)}{2(x^{2}-1)} \\
-\frac{x(2x^{2}-1)}{\sqrt{x^{2}-1}} & \frac{\sqrt{-3+4x^{2}}(2x^{2}-1)}{%
2(x^{2}-1)} & \frac{1+2x^{2}-4x^{4}}{2-2x^{2}}%
\end{pmatrix}
\end{equation*}

The maps $\rho _{x}^{\prime }(a)$ and $\rho _{x}^{\prime }(b)$ are
hyperbolic rotations on $H_{0}\cong E^{1,2}$ moving $\delta $ along the
polars of the space-like vectors $A^{-}$ and $B^{-}$ where $
x=A^{+}=B^{+}=\cosh \frac{\delta }{2}$. See Figure
\ref{fsubcase2a2}. The distance $d$ (measured in the hyperbolic
plane) between the polars of the axes of $\rho _{x}^{\prime }(a)$
and $\rho _{x}^{\prime }(b)$ is given by
\begin{equation*}
\cosh d=\frac{y}{x^{2}-1}=\frac{x^{2}-\frac{1}{2}}{x^{2}-1}
\end{equation*}
\end{description}

\subsubsection{The Figure Eight knot}

\begin{figure}[ht]
\epsfig{file=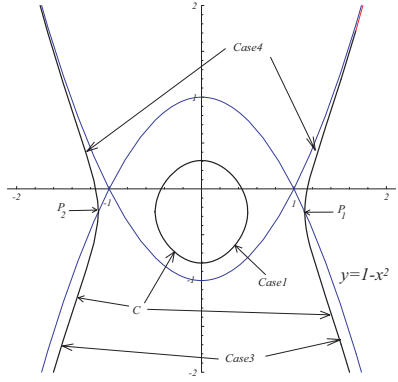,height=5cm}\caption{The real part of
$V(\mathcal{I}_{G(4_{1})})$.}\label{fcurvaocho}
\end{figure}

The algebraic variety $V(\mathcal{I}_{G(4_{1})})$ for the
Figure Eight knot group
\begin{equation*}
G(4_{1})=|a,b;aw=wb|,\qquad w=ba^{-1}b^{-1}a,
\end{equation*}
is defined by the ideal
 $\mathcal{I}_{G(4_{1})}=(1-6x^{2}+4x^{4}-2y-2y^{2})$. See Example \ref{ejemplo dos}. The real part of
the algebraic variety $V(\mathcal{I}_{G(4_{1})})$ is the curve
$\mathcal{C}$ depicted in Figure \ref{fcurvaocho}, together with the two
parabolas $y=\pm (x^{2}-1)$ of the variety of reducible
representations $ V(y^{2}-(x^{2}-1)^{2})$.

\begin{figure}[ht]
\epsfig{file=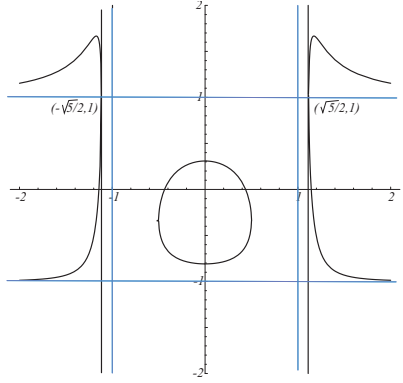,height=5cm}\caption{The function
$\frac{y}{1-x^{2}}$ for the Figure Eight
knot.}\label{fcurvaysobreuuocho}
\end{figure}

Figure \ref{fcurvaysobreuuocho} shows $\frac{y}{1-x^{2}}$ as a
function of $x$ for the algebraic variety
$V(\mathcal{I}_{G(4_{1})})$. As before, it is enough to consider
this Figure over the pattern of Figure \ref{fregiones} to classify
the different classes of representation of the group $G(4_{1})$ in
$S^{3}$ or $SL(2,\mathbb{R})$ according to Theorem \ref{Teorema
dos}.

There are four cases:
\begin{description}
\item[Case 1]  Region (1). $x\in \lbrack -\frac{1}{2},\frac{1}{2}
]\Longleftrightarrow \left\{ 1-x^{2}>0,(1-x^{2})^{2}>y^{2}\right\}
$.There are two different $y$ values for each value of $x$, except
for $x=\pm 1/2$ where they coincide.
\begin{eqnarray*}
y_{1} &=&\frac{1}{4}\left( -1+\sqrt{5-24x^{2}+16x^{4}}\right) \\
y_{2} &=&\frac{1}{4}\left( -1-\sqrt{5-24x^{2}+16x^{4}}\right)
\end{eqnarray*}
There exists an irreducible c-representation $\rho
_{xi}:G(4_{1})\longrightarrow S^{3}$ realizing $(x,y_{i})$, $i=1,2,$
unique up to conjugation in $S^{3}$, such that
\begin{equation}
\begin{array}{l}
\rho
_{x1}(a)=A=x+\frac{-1+\sqrt{5-24x^{2}+16x^{4}}}{4\sqrt{1-x^{2}}}i+\frac{
1}{2}\sqrt{\frac{5-4x^{2}+\sqrt{5-24x^{2}+16x^{4}}}{2-x^{2}}}j \\
\rho _{x1}(b)=B=x+\sqrt{1-x^{2}}i,\quad \sqrt{1-x^{2}}>0
\end{array}
\label{ecase11}
\end{equation}
\begin{equation*}
\begin{array}{l}
\rho
_{x2}(a)=A=x+\frac{-1-\sqrt{5-24x^{2}+16x^{4}}}{4\sqrt{1-x^{2}}}i+\frac{
1}{2}\sqrt{\frac{5-4x^{2}-\sqrt{5-24x^{2}+16x^{4}}}{2-x^{2}}}j \\
\rho _{x2}(b)=B=x+\sqrt{1-x^{2}}i,\quad \sqrt{1-x^{2}}>0
\end{array}
\end{equation*}
The composition of $\rho _{xi}$ with $c:S^{3}\rightarrow SO(3)$,
where $c(X)$ , $X\in S^{3}$, acts on $P\in H_{0}\cong E^{3}$ by
conjugation, defines the representation $\rho _{xi}^{\prime }=c\circ
\rho _{xi}:G(3/1)\longrightarrow SO(3)$. The maps $\rho
_{xi}^{\prime }(a)$ and $\rho _{xi}^{\prime }(b)$ are right
rotations of angle $\alpha $ around the axes $A^{-}$ and $B^{-}$
where $ x=A^{+}=B^{+}=\cos \frac{\alpha }{2}$.

The angle $\omega _{i}$ between the axes of $\rho _{xi}^{\prime
}(a)$ and $ \rho _{xi}^{\prime }(b)$ is given by $\cos \omega
_{i}=\frac{y_{i}}{u}$
\begin{eqnarray*}
\cos \omega _{1} &=&\frac{y_{1}}{u}\,=\frac{\frac{1}{4}\left(
-1+\sqrt{
5-24x^{2}+16x^{4}}\right) }{1-x^{2}} \\
\cos \omega _{2} &=&\frac{y_{2}}{u}\,=\frac{\frac{1}{4}\left(
-1-\sqrt{ 5-24x^{2}+16x^{4}}\right) }{1-x^{2}}
\end{eqnarray*}

\item[Case 2] Segment (2.4) $(x,y)=(\pm \frac{\sqrt{5}}{2},-\frac{1}{4}%
)\Longleftrightarrow \left\{ 1-x^{2}<0,(1-x^{2})^{2}=y^{2}\right\}
$. There exists an almost-irreducible c-representation $\rho
_{x}:G(4_{1})\longrightarrow U_{1}\subset \left( \frac{-1,1}{\mathbb{R}}%
\right) $ realizing $(x,y)$, unique up to conjugation in $U_{1}$,
such that:
\begin{eqnarray}
\rho _{\pm \sqrt{5}/2}(a) &=&A=\frac{\pm \sqrt{5}}{2}+\frac{1}{2}J
\label{ecaso22} \\
\rho _{\pm \sqrt{5}/2}(b) &=&B=\frac{\pm
\sqrt{5}}{2}+I+\frac{1}{2}J+IJ \notag
\end{eqnarray}%
Under the isomorphism $U_{1}\approx SL(2,\mathbb{R})$ we have
\begin{equation*}
\begin{array}{cccc}
\rho _{\pm \sqrt{5}/2}: & G(4_{1}) & \longrightarrow & SL(2,\mathbb{R}) \\
& a & \rightarrow & A=
\begin{pmatrix}
1+\frac{\pm \sqrt{5}}{2} & \frac{3}{2} \\
-\frac{1}{2} & -1+\frac{\pm \sqrt{5}}{2}
\end{pmatrix}
\\
& b & \rightarrow & B=
\begin{pmatrix}
\frac{\pm \sqrt{5}}{2} & \frac{1}{2} \\
\frac{1}{2} & \frac{\pm \sqrt{5}}{2}
\end{pmatrix}
\end{array}
\end{equation*}
The composition of $\rho _{\pm \sqrt{5}/2}$ with $c:SL(2,\mathbb{R}
)\rightarrow SO^{+}(1,2)\cong Iso^{+}(\mathbb{H}^{2})$, where
$c(X)$, $X\in SL(2,\mathbb{R})$, acts on $P\in H_{0}\cong E^{1,2}$
by conjugation, defines the representation $\rho _{\pm
\sqrt{5}/2}^{\prime }=c\circ \rho _{\pm \sqrt{
5}/2}:G(3/1)\longrightarrow SO^{+}(1,2)$ in affine linear notation,
where $ \left\{ X,Y,Z\right\} $ is the coordinate system associated
to the basis $ \left\{ -IJ,J,I\right\} $
\begin{equation*}
\begin{array}{l}
\rho _{\pm \sqrt{5}/2}^{\prime }(a)=m_{\pm \sqrt{5}/2}(a)\left(
\begin{array}{c}
X \\
Y \\
Z%
\end{array}
\right) =\left(
\begin{array}{c}
X^{\prime } \\
Y^{\prime } \\
Z^{\prime }
\end{array}
\right) \\
\rho _{\pm \sqrt{5}/2}^{\prime }(b)=m_{\pm \sqrt{5}/2}(b)\left(
\begin{array}{c}
X \\
Y \\
Z
\end{array}
\right) =\left(
\begin{array}{c}
X^{\prime } \\
Y^{\prime } \\
Z^{\prime }
\end{array}
\right)
\end{array}
\end{equation*}
where the matrices of $\rho _{\pm \sqrt{5}/2}^{\prime
}(a)$ and $ \rho _{\pm \sqrt{5}/2}^{\prime }(b)$ are respectively:
\begin{equation*}
m_{\pm \sqrt{5}/2}(a)=\left(
\begin{array}{ccc}
\frac{3}{2} & 0 & \frac{\sqrt{5}}{2} \\
0 & 1 & 0 \\
\frac{\sqrt{5}}{2} & 0 & \frac{3}{2}
\end{array}
\right)
\end{equation*}
and
\begin{equation*}
m_{\pm \sqrt{5}/2}(b)=
\begin{pmatrix}
-\frac{1}{2} & 1-\sqrt{5} & \frac{1}{2}(\sqrt{5}-4) \\
1+\sqrt{5} & 1 & 1+\sqrt{5} \\
\frac{1}{2}(\sqrt{5}+4) & -1+\sqrt{5} & \frac{7}{2}
\end{pmatrix}
\end{equation*}
The maps $\rho _{x}^{\prime }(a)$ and $\rho _{x}^{\prime }(b)$ are
hyperbolic rotations on $H_{0}\cong E^{1,2}$ around the space-like
axes $ A^{-}$ and $B^{-}$, where $x=A^{+}=B^{+}=\cosh \frac{d}{2}$.
See Figure \ref {subcase2c1}.

\item[Case 3] Region (2.2) $\left\vert x\right\vert >1,(1-x^{2})^{2}<y^{2}.$
This case consists of the real points $(x,y)\in V(\mathcal{I}_{G})$
where $ \left\vert x\right\vert >1$ and $y=\frac{1}{4}\left(
-1-\sqrt{ 5-24x^{2}+16x^{4}}\right) \leq -\frac{1}{4}$.

There exists an irreducible c-representation $\rho
_{x}:G(3/1)\longrightarrow SL(2,\mathbb{R})=U_{1}\subset \left(
\frac{-1,1}{ \mathbb{R}}\right) $ realizing $(x,y)$, unique up to
conjugation in $SL(2, \mathbb{R})$, such that
\begin{equation*}
\begin{array}{l}
\rho _{x}(a)=A=x+\sqrt{x^{2}-1}J,\quad \sqrt{x^{2}-1}>0 \\
\rho
_{x}(b)=B=x+\frac{1}{2}\sqrt{\frac{5-4x^{2}-\sqrt{5-24x^{2}+16x^{4}}}{
2x^{2}-2}}I+\frac{1+\sqrt{5-24x^{2}+16x^{4}}}{4\sqrt{x^{2}-1}}J
\end{array}
\end{equation*}
The composition of $\rho _{x}$ with $c:SL(2,\mathbb{R})\rightarrow
SO^{+}(1,2)\cong Iso^{+}(\mathbb{H}^{2})$, where $c(X)$, $X\in
SL(2,\mathbb{R })$, acts on $P\in H_{0}\cong E^{1,2}$ by
conjugation, defines the representation $\rho _{x}^{\prime }=c\circ
\rho _{x}:G(3/1)\longrightarrow SO^{+}(1,2)$.

The maps $\rho _{x}^{\prime }(a)$ and $\rho _{x}^{\prime }(b)$ are
hyperbolic rotations on $H_{0}\cong E^{1,2}$ around the space-like
axes $ A^{-}$ and $B^{-}$, where $x=A^{+}=B^{+}=\cosh \frac{d}{2}$.
See Figure \ref {fsubcase2a2}.

The distance $\delta $ between the polars of the axes of $\rho
_{x}^{\prime }(a)$ and $\rho _{x}^{\prime }(b)$ (measured in the
hyperbolic plane) is given by
\begin{equation*}
\cosh \delta =\frac{-y}{u}=\frac{-\frac{1}{4}\left( -1-\sqrt{
5-24x^{2}+16x^{4}}\right) }{x^{2}-1}>1
\end{equation*}

\item[Case 4] Region (2.3) $\left\vert x\right\vert >1,(1-x^{2})^{2}>y^{2}.$
This case consists of the real points $(x,y)\in V(\mathcal{I}_{G})$
where $ \left\vert x\right\vert >1$ and $y=\frac{1}{4}\left(
-1+\sqrt{ 5-24x^{2}+16x^{4}}\right) \geq -\frac{1}{4}$.

There exists an irreducible c-representation $\rho
_{x}:G(3/1)\longrightarrow SL(2,\mathbb{R})=U_{1}\subset \left(
\frac{-1,1}{ \mathbb{R}}\right) $ realizing $(x,y)$, unique up to
conjugation in $SL(2,
\mathbb{R})$, such that
\begin{equation*}
\begin{array}{l}
\rho _{x}(a)=A=x+\sqrt{x^{2}-1}J \\
\rho _{x}(b)=B=x+\frac{1-\sqrt{5-24x^{2}+16x^{4}}}{4\sqrt{x^{2}-1}}J+\frac{1
}{2}\sqrt{\frac{5-4x^{2}+\sqrt{5-24x^{2}+16x^{4}}}{2x^{2}-2}}IJ
\end{array}
\end{equation*}
The composition of $\rho _{x}$ with $c:SL(2,\mathbb{R})\rightarrow
SO^{0}(1,2)\cong Iso^{+}(\mathbb{H}^{2})$, where $c(X)$, $X\in
SL(2,\mathbb{R })$, acts on $P\in H_{0}\cong E^{1,2}$ by
conjugation, defines the representation $\rho _{x}^{\prime }=c\circ
\rho _{x}:G(3/1)\longrightarrow SO^{0}(1,2)$.

The maps $\rho _{x}^{\prime }(a)$ and $\rho _{x}^{\prime }(b)$ are
hyperbolic rotations on $H_{0}\cong E^{1,2}$ around the space-like
vectors $ A^{-}$ and $B^{-}$, where $x=A^{+}=B^{+}=\cosh
\frac{d}{2}$. See Figure \ref {subcase2b2}.

The angle $\theta $ (measured in the hyperbolic plane) between the
polars of the axes of $\rho _{x}^{\prime }(a)$ and $\rho
_{x}^{\prime }(b)$ is given by
\begin{equation*}
\cos \theta
=\frac{-y}{x^{2}-1}=\frac{1-\sqrt{5-24x^{2}+16x^{4}}}{4(x^{2}-1)}
\end{equation*}

\end{description}

We have omitted the expression of $\rho _{x}^{\prime }(a)$ and
$\rho _{x}^{\prime }(b)$ in affine linear notation,
where $ \left\{ X,Y,Z\right\} $ is the coordinate system associated
to the basis $ \left\{ -ij,j,i\right\} $ for the above cases 1, 3 and 4,
but they can be easily obtained using the equation (\ref{emtxyz}) of \S 2.3.1.

\section{Groups of isometries in a quaternion algebra H
\label{Sec4}}

Let $H=\left( \frac{\mu ,\nu }{k}\right) $ be a quaternion algebra. The pure
quaternions form a 3-dimensional vector space $H_{0}$.

The following map is a left action of the group $U$ on the 3-dimensional
vector space $H_{0}$.
\begin{equation*}
\begin{array}{llll}
\Phi : & U\times H_{0} & \longrightarrow & H_{0} \\
& (A,B^{-}) & \rightarrow & A\circ B^{-}:=AB^{-}\overline{A}
\end{array}%
\end{equation*}

The restriction of $\Phi $ to the subgroup $U_{1}$ (the unit quaternions,
norm 1) is also a left action on the 3-dimensional vector space $H_{0}$.%
\begin{equation*}
\begin{array}{llll}
\Phi _{0}: & U_{1}\times H_{0} & \longrightarrow  & H_{0} \\
& (A,B^{-}) & \rightarrow  & A\circ B^{-}:=AB^{-}\overline{A}
\end{array}
\end{equation*}

The \emph{equiform group} or \emph{group of similarities} $\mathcal{E}q(H)$
of a quaternion algebra  $H$, is the semidirect product $H_{0}\rtimes U$.
This is the group whose underlying space is $H_{0}\times U$ and the product
is
\begin{equation*}
\begin{array}{lll}
\mathcal{E}q(H)\times \mathcal{E}q(H) & \longrightarrow & \mathcal{E}q(H) \\
((v,A),(w,B)) & \rightarrow & (v+A\circ w,AB)%
\end{array}%
\end{equation*}

The group of \emph{affine isometries} $A(H)$ of a quaternion algebra is the
subgroup of $\mathcal{E}q(H)\ $which is the semidirect product $H_{0}\rtimes
U_{1}$. This is the group whose underlying space is $H_{0}\times U_{1}$ and
the product is
\begin{equation*}
\begin{array}{lll}
A(H)\times A(H) & \longrightarrow & A(H) \\
((v,A),(w,B)) & \rightarrow & (v+A\circ w,AB)
\end{array}
\end{equation*}

The group $\mathcal{E}q(H)$ defines a left action on the 3-dimensional vector
space $H_{0}$,
\begin{equation*}
\begin{array}{llll}
\Psi : & \mathcal{E}q(H_{0})\times H_{0} & \longrightarrow & H_{0} \\
& ((v,A),u) & \rightarrow & (v,A)u:=v+A\circ u
\end{array}
\end{equation*}
because
\begin{equation*}
\Psi ((0,1),u)=(0,1)u=0+1u\overline{1}=u
\end{equation*}
\begin{eqnarray*}
\Psi ((v_{2},A_{2}),\Psi ((v_{1},A_{1}),u))
&=&(v_{2},A_{2})(v_{1}+A_{1}\circ u)=v_{2}+A_{2}\circ (v_{1}+A_{1}\circ u) \\
&=&v_{2}+A_{2}\circ v_{1}+A_{2}A_{1}\circ u=\Psi ((v_{2}+A_{2}\circ
v_{1},A_{2}A_{1}),u) \\
&=&\Psi ((v_{2},A_{2})(v_{1},A_{1}),u)
\end{eqnarray*}

The restriction of $\Psi $ to the subgroup $A(H)$ is also a left action on
the 3-dimensional vector space $H_{0}$.

For an element $(v,A)\in \mathcal{E}q(H)$, $A$ is the \emph{linear
part} of $(v,A)$, $N(A)$ is the \emph{homothetic factor}, and $v$
is the \emph{translational part}. Note that if $(v,A)\in A(H)$, then
the homothetic factor is 1. The action of each element $(v,A)\in A(H)$ in $H_{0}$ is called an \emph{affine isometry}.

\begin{example}
1.- $H=\left( \frac{-1,-1}{\mathbb{R}}\right) =\mathbb{H}$. Then
$\mathcal{E} q(H)$ is the group of affine isometries and
similarities of the Euclidean 3-dimensional space $E^{3}$. The
linear part of $\mathcal{E}q(H)$ is the multiplicative group
$(\mathbb{H}\backslash \{0\})$ of the algebra $\mathbb{H }$ . The
subgroup $A(H)$ is the group of orientation preserving affine
isometries of $E^{3}$. It is called the \emph{Euclidean group} and
we denote it by $\mathcal{E}(\mathbb{R}^{3})$. The linear part of
$A(H)$ is $SO(3,\mathbb{R})$.

2.- $H=\left( \frac{-1,1}{\mathbb{R}}\right) =M(2,\mathbb{R})$. Then
$ \mathcal{E}q(H)$ is the group of affine isometries and
similarities of the Minkowski 3-dimensional space $E^{1,2}$. The
linear part of $\mathcal{E}q(H)$ is isomorphic to
$Gl(2,\mathbb{R})$. The subgroup $A(H)$ is the group of orientation
preserving affine isometries of $E^{1,2}$. It is called the
\emph{Lorentz group} and we denote it by
$\mathcal{L}(\mathbb{R}^{3})$. The linear part is $SO^{0}(1,2)$.

3.- $H=\left( \frac{-1,1}{\mathbb{C}}\right) =M(2,\mathbb{C})$. Then
$ \mathcal{E}q(H)$ is the group of affine isometries and
similarities of the complex 3-dimensional space $\mathbb{C}^{3}$.
The linear part of $\mathcal{E} q(H)$ is isomorphic to
$Gl(2,\mathbb{C})$. Here $A(H)$ is the subgroup of orientation
preserving affine isometries of $\mathbb{C}^{3}$. The linear part of
$A(H)$ is $SO(3,\mathbb{C}).$
\end{example}

\subsubsection{Axis and shift of an element $(v,A)\in A(H)$}

Consider $(v,A)\in A(H).$ The vector $v\in H_{0}$ can be decomposed
in a unique way as the orthogonal sum of two vectors, one of them in
the $A^{-}$ direction:
\begin{equation*}
v=sA^{-}+v^{\perp },\quad \qquad \left\langle v^{\perp },A^{-}\right\rangle
=0
\end{equation*}
Then
\begin{equation*}
(v,A)=(v^{\perp },1)(sA^{-},A)
\end{equation*}
The element $(v^{\perp },1)$ is a translation in $H_{0}$. The
restriction of the action of $(sA^{-},A)$ on the line generated by
$A^{-}$ is a translation with vector $sA^{-}$:
\begin{equation*}
(sA^{-},A)(\lambda A^{-})=sA^{-}+A(\lambda A^{-})\overline{A}=sA^{-}+\lambda
A^{-}=(s+\lambda )A^{-}
\end{equation*}
We define $sA^{-}$ as the \emph{vector shift of the element
}$(v,A)$. The length $\sigma $ of the vector shift will be called
the \emph{shift of the element }$(A,v)$.

The action of $(v,A)\ $leaves (globally) invariant an affine line parallel
to $A^{-}$ and its action on this line is a translation with vector the shift $%
sA^{-}$ . This invariant affine line will be call \emph{the axis of} $(v,A)$
. Then the action of $(v,A)$ on the axis of $(v,A)$ is a translation by $%
\sigma $. It is easy to see that the axis of $(v,A)$ is $\left\{ u+\mu
A^{-}:\mu \in k\right\} $ where $u$ is defined by $u-A\circ u=v^{\perp }$.
In fact,
\begin{equation}
\begin{array}{c}
(v,A)(u+\mu A^{-}) =u+(\mu +s)A^{-}\quad   \\
\Longrightarrow \quad v+A\circ (u+\mu A^{-})=u+(\mu +s)A^{-}\quad \\
\Longrightarrow \quad sA^{-}+v^{\perp }\quad +A\circ u+\mu
A^{-}=u+(\mu
+s)A^{-} \\
\Longrightarrow \quad u-A\circ u=v^{\perp }+(s+\mu -(\mu
+s))A^{-}=v^{\perp }\label{eaxis}
\end{array}
\end{equation}

\begin{figure}[ht]
\epsfig{file=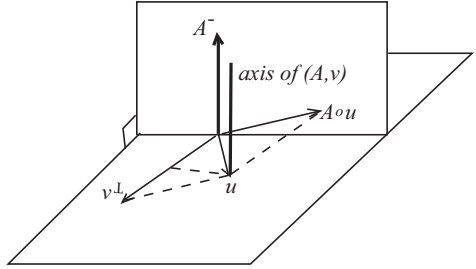,height=3.5cm}\caption{The axis of $(A,v)$.}\label{faxixav}
\end{figure}

\begin{remark}
\label{sameaxis} Observed that $(v_{1},A)$ and $(v_{2},A)$ have the
same axis if and only if $v_{1}^{\perp }=v_{2}^{\perp }$. Therefore,
given an element $(v,A)$, the element $(v^{\perp },A)$, where
$v=sA^{-}+v^{\perp }$, $\left\langle v^{\perp },A^{-}\right\rangle
=0$, has the same axis that $(v,A)$ but shift zero.
\end{remark}

\subsection{Two conjugate elements of $A(H)$.}

If $(v,A),(w,B)\in A(H)$ are conjugate elements in $A(H)$ then $A$
and $B$ are conjugate in $U_{1}$ and the vector shifts of $(v,A)$
and $(w,B)$ are respectively $sA^{-}$ and $sB^{-}$. In fact, a pair
$(v,A),(w,B)\in A(H)$ of conjugate elements in $A(H),$ is conjugate
to the pair $(sA^{-},A),$ $ (sB^{-}+w'^{\perp },B)$, where
$\left\langle w'^{\perp },B^{-}\right\rangle =0$. This conjugation
is made by a translation (change of the origin point in $ H_{0}$).

\begin{lemma}
\label{lema} If the elements $(sA^{-},A)$,  $(sB^{-}+w^{\perp },B)$,
where $ \left\langle w^{\perp },B^{-}\right\rangle =0$, are
conjugate elements in $ A(H)$ or $\mathcal{E}q(H)$ by an element
$(u,C)$ then, for any $\lambda \in k $, $(\lambda sA^{-},A)$,
$(\lambda sB^{-}+w^{\perp },B)$ are also conjugate by the same
element.
\end{lemma}

\begin{proof}
Let $(u,C)\in A(H)$ or $\mathcal{E}q(H)$ be the element such that
\begin{equation*}
(u,C)(sA^{-},A)(u,C)^{-1}=(sB^{-}+w^{\perp },B)
\end{equation*}
Then
\begin{eqnarray*}
(u,C)(sA^{-},A)(u,C)^{-1} &=&(u,C)(sA^{-},A)(-C^{-1}\circ u,C^{-1}) \\
&=&(u,C)(sA^{-}-(AC^{-1})\circ u,AC^{-1}) \\
&=&(u+C\circ sA^{-}-(CAC^{-1})\circ u,CAC^{-1}) \\
&=&(u+sCA^{-}C^{-1}-(CAC^{-1})\circ u,CAC^{-1}) \\
&=&(sB^{-}+w^{\perp },B)
\end{eqnarray*}
implies that
\begin{eqnarray*}
CAC^{-1} &=&B \\
u-CAC^{-1}\circ u &=&w^{\perp }
\end{eqnarray*}
Then
\begin{equation*}
(u,C)(\lambda sA^{-},A)(u,C)^{-1}=(\lambda sB^{-}+w^{\perp },B)
\end{equation*}
\end{proof}

\begin{proposition}
\label{paxesayb} If $(A,B)$ is a irreducible pair, that is $\{A^{-},B^{-},(A^{-}B^{-})^{-}\}$is a basis of $H_{0}$,
 then the axes of $(sA^{-},A)$ and $
(sB^{-}+(A^{-}B^{-})^{-},B)$ do not intersect.
\end{proposition}

\begin{proof}
By Remark \ref{sameaxis}, the axes of $(sA^{-},A)$ and $(0,A)$ are
both the vector line generated by $A^{-}$. By the same reason the
axes of $ (sB^{-}+(A^{-}B^{-})^{-},B)$ and $((A^{-}B^{-})^{-},B)$ also coincide. By (\ref{eaxis}), the axis of $((A^{-}B^{-})^{-},B)$
is $\left\{ u+\mu B^{-},\mu \in k\right\} $, such that $u-B\circ
u=(A^{-}B^{-})^{-}$. Let us prove by contradiction that $u$ is not
contained in the plane generated by $\left\{ A^{-},B^{-}\right\}$,
where the axis of $(0,A)$ lives. Assume that $u$ belongs to the
plane generated by $\left\{ A^{-},B^{-}\right\} $, then by Corollary
\ref{cambmmortoambm}  $u\perp (A^{-}B^{-})^{-}$, $\left\langle
u,(A^{-}B^{-})^{-}\right\rangle =0$. Thus, since $B$ is an isometry,
\begin{eqnarray*}
\left\langle u,u\right\rangle  &=&\left\langle B\circ u,B\circ
u\right\rangle =\left\langle
u-(A^{-}B^{-})^{-},u-(A^{-}B^{-})^{-}\right\rangle  \\
&=&\left\langle u,u\right\rangle +\left\langle
(A^{-}B^{-})^{-},(A^{-}B^{-})^{-}\right\rangle  \\
&\Longrightarrow &\left\langle
(A^{-}B^{-})^{-},(A^{-}B^{-})^{-}\right\rangle =N((A^{-}B^{-})^{-})=0
\end{eqnarray*}
Then $(A^{-}B^{-})^{-}\in ((A^{-}B^{-})^{-})^{\perp }$. But the
plane $ ((A^{-}B^{-})^{-})^{\perp }$ is the plane generated by
$\left\{ A^{-},B^{-}\right\}$, which is impossible because $\left\{
A^{-},B^{-},(A^{-}B^{-})^{-}\right\} $ is a basis. Therefore the
axis $ \left\{ \lambda A^{-}:\lambda \in k\right\} $ is contained in
the plane $\Pi $ generated by $\left\{ A^{-},B^{-}\right\}$, and
the axis $\left\{ u+\mu B^{-}:\mu \in k\right\} $ is contained in
the plane $u+\Pi $ parallel to $ \Pi $ but different from $\Pi $.
Then both axes do not intersect.
\end{proof}

Next, we will prove that if $\left\{ A^{-},B^{-},(A^{-}B^{-})^{-}\right\} $
is a basis in $H_{0}$ and $(v,A),(w,B)\in A(H)$ is a pair of conjugate
elements whose axes do not intersect, it is possible to conjugate them in $%
\mathcal{E}q(H)$ to the standard elements of the above Proposition \ref%
{paxesayb}.

\begin{theorem}
\label{teo2}Let $(v,A),(w,B)\in A(H)$ be a pair of conjugate elements in $%
A(H)$ whose axes do not intersect. Assume $(A,B)$ is a irreducible pair. Then the pair $%
((v,A),(w,B))\in A(H)$ is conjugate in $\mathcal{E}q(H)$ to a pair of the
form $((sA^{-},A)$ $,$ $(sB^{-}+(A^{-}B^{-})^{-},B))$ or $((sA^{-},A)$ $,$ $%
(sB^{-}-(A^{-}B^{-})^{-},B))$.
\end{theorem}

\begin{proof}
Up to conjugation by a translation (change of the origin point in $H_{0}$)
we can assume that $(v,A)=(tA^{-},A)$, and $(w,B)=(tB^{-}+w^{\perp },B)$,
where $w^{\perp }\perp B^{-}$.

The idea is to apply a conjugation by an element of
$\mathcal{E}q(H)$ to a pair of the form $((sA^{-},A)$,
$(sB^{-}\pm (A^{-}B^{-})^{-},B))$ to obtain the pair
$((tA^{-},A)$,$(tB^{-}+w^{\perp },B))$. Note that conjugating by the
translations $\left\{ (\lambda A^{-},1),\lambda \in k\right\} $  does not change the axis of $(sA^{-},A)$. However the axis of $
(sB^{-}\pm (A^{-}B^{-})^{-},B))$ generates a family $\mathcal{F}$ of
non intersecting lines whose union is two planes $\mathcal{P}$
parallel to the plane $\Pi =\{A^{-},B^{-}\}$ and placed on both
sides of it. Conjugating by the set of homothetic transformations
$\left\{ (0,rI),r\in k^{\ast }\right\} $ the axis of $(sA^{-},A)$
does not change, but the family $\mathcal{F}$ contained in
$\mathcal{P}$ generates a family of non intersecting lines whose
union is the complement of the plane $\Pi =\{A^{-},B^{-}\}.$ See
Figure \ref{fteorema}.

\begin{figure}[ht]
\epsfig{file=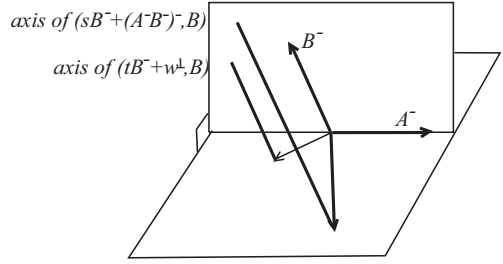,height=3.5cm}\caption{Axes of $ (w,B)$ and
$(sB^{-}+(A^{-}B^{-})^{-},B))$.}\label{fteorema}
\end{figure}

Actually, we can find an element $(\lambda A^{-},rI)\in \mathcal{E}q(H)$
such that
\begin{eqnarray}
(\lambda A^{-},rI)(sA^{-},A)(\lambda A^{-},rI)^{-1} &=&(r^{2}sA^{-},A)
\label{conj1} \\
(\lambda A^{-},rI)(sB^{-}\pm (A^{-}B^{-})^{-},B)(\lambda A^{-},rI)^{-1}
&=&(r^{2}sB^{-}+w^{\perp },B)  \label{conj2}
\end{eqnarray}
The first equation (\ref{conj1}) is always true. From the second one (\ref%
{conj2}) we can obtain the values of $r$ and $\lambda $:
\begin{eqnarray*}
&&(\lambda A^{-},rI)(sB^{-}\pm (A^{-}B^{-})^{-},B)(\lambda A^{-},rI)^{-1} \\
&=&(\lambda A^{-},rI)(sB^{-}\pm (A^{-}B^{-})^{-},B)(-(\frac{1}{r}I)\circ
(\lambda A^{-}),\frac{1}{r}I) \\
&=&(\lambda A^{-},rI)\left( sB^{-}\pm (A^{-}B^{-})^{-}-B\circ (\frac{1}{r^{2}%
}\lambda A^{-}),\frac{1}{r}B\right) \\
&=&\left( \lambda A^{-}+r^{2}(sB^{-}\pm (A^{-}B^{-})^{-}-\frac{\lambda }{%
r^{2}}BA^{-}\overline{B},B\right) =(v,B)
\end{eqnarray*}
where%
\begin{equation*}
v=\lambda A^{-}+r^{2}(sB^{-}\pm (A^{-}B^{-})^{-}-\frac{\lambda
}{r^{2}}BA^{-} \overline{B}
\end{equation*}

Let us compute $BA^{-}\overline{B}$ in the basis $\left\{
A^{-},B^{-},(A^{-}B^{-})^{-}\right\} $
\begin{equation*}
BA^{-}\overline{B}
=(B^{+}+B^{-})A^{-}(B^{+}-B^{-})=(B^{+})^{2}A^{-}+B^{+}(B^{-}A^{-}-A^{-}B^{-})-B^{-}A^{-}B^{-}
\end{equation*}
Using the equation
\begin{equation*}
A^{-}B^{-}+B^{-}A^{-}=2(A^{-}B^{-})^{+}
\end{equation*}
and the notation
\begin{eqnarray*}
x &=&B^{+} \\
u &=&-B^{-}B^{-}=1-x^{2} \\
y &=&-(A^{-}B^{-})^{+}
\end{eqnarray*}
We have
\begin{equation*}
BA^{-}\overline{B}
=x^{2}A^{-}-2x(A^{-}B^{-})^{-}+2yB^{-}-uA^{-}=(2x^{2}-1)A^{-}+2yB^{-}-2x(A^{-}B^{-})^{-}
\end{equation*}
Therefore
\begin{equation*}
v=\lambda (2-2x^{2})A^{-}+(r^{2}s-2\lambda y)B^{-}+(\pm r^{2}-2\lambda
x)(A^{-}B^{-})^{-}
\end{equation*}
and if
\begin{equation*}
(v,B)=(r^{2}sB^{-}+w^{\perp },B)
\end{equation*}
then
\begin{eqnarray*}
v &=&r^{2}sB^{-}+w^{\perp }\Longrightarrow \\
w^{\perp } &=&\lambda (2-2x^{2})A^{-}-2\lambda yB^{-}+(\pm r^{2}-2\lambda
x)(A^{-}B^{-})^{-}
\end{eqnarray*}
where $w^{\perp }=aA^{-}+bB^{-}+c(A^{-}B^{-})^{-}$ is such that $
\left\langle w^{\perp },B^{-}\right\rangle =0.$ Thus%
\begin{eqnarray*}
0 &=&\left\langle aA^{-}+bB^{-}+c(A^{-}B^{-})^{-},B^{-}\right\rangle \\
&=&a\left\langle A^{-},B^{-}\right\rangle +b\left\langle
B^{-},B^{-}\right\rangle +c\left\langle (A^{-}B^{-})^{-},B^{-}\right\rangle
=ay+b(1-x^{2})
\end{eqnarray*}

We obtain the system
\begin{eqnarray*}
a &=&\lambda (2-2x^{2}) \\
b &=&-2\lambda y \\
c &=&\pm r^{2}-2\lambda x
\end{eqnarray*}
Then
\begin{eqnarray*}
&&\fbox{$\lambda =\frac{a}{2-2x^{2}}$} \\
&&\fbox{$\pm r^{2}=c+\frac{ax}{1-x^{2}}$}
\end{eqnarray*}
Here, if $k=\mathbb{R}$, two cases are possible:%
\begin{eqnarray*}
c+\frac{ax}{1-x^{2}} &>&0,\text{take the }+\text{ sign,
}r=\sqrt{c+\frac{ax}{
1-x^{2}}} \\
c+\frac{ax}{1-x^{2}} &<&0,\text{take the }-\text{ sign,
}r=\sqrt{-\left( c+ \frac{ax}{1-x^{2}}\right) }
\end{eqnarray*}

We have proven that, in the first case, the pair $((sA^{-},A)$ $,$ $
(sB^{-}+(A^{-}B^{-})^{-},B))$ is conjugate in $\mathcal{E}q(H)$ to
the pair $ \left( (r^{2}sA^{-},A),(r^{2}sB^{-}+w^{\perp },B)\right)
$, and that in the second case, the pair $((sA^{-},A)$,
$(sB^{-}-(A^{-}B^{-})^{-},B))$ is conjugate in $\mathcal{E}q(H)$ to
the pair $\left( (r^{2}sA^{-},A),(r^{2}sB^{-}+w^{\perp },B)\right)
$. By Lemma \ref{lema} we can deduce that the pair
($(v,A)=(tA^{-},A),(w,B)=(tB^{-}+w^{\perp },B))$ is conjugate to the
pair $((tA^{-},A),(tB^{-}\pm (A^{-}B^{-})^{-},B))$.
\end{proof}

\begin{corollary}
Let $(v,A),(w,B)\in A(H)$ be a pair of conjugate elements in $%
A(H)$ whose axes do not intersect. Assume $(A,B)$ is a irreducible pair. Then $%
(v,A),(w,B)$ are determined up to similarity in $H_{0}$ by the parameters $%
x,y,s$.\qed
\end{corollary}

The geometrical meaning of the parameters depends on the geometry of $H_{0}$%
, that is, on the quaternion algebra $H=\left( \frac{\mu ,\nu }{k}\right) $.

\subsection{c-Representations in $A(H)$}

Let $G=\left\vert a,b;w(a,b)\right\vert $ be a presentation of a
group. For example, $G=\left\vert a,b;w(a,b)\right\vert $ can be a
presentation of the group of a 2-bridge knot, where $a$ and $b$ are
represented by coherently oriented meridians of the knot. We want to
study the \emph{c-representations of }$G$\emph{\ in the affine group
}$A(H)$ of a quaternion algebra $H$, that is representations of $G$
in the affine group $A(H)$ of a quaternion algebra $H$ such that the
generators $a$ and $b$ go to conjugate elements, up to conjugation
in $\mathcal{E}q(H)$. We have already studied the case of the
c-representation whose image lies in the subgroup $U_{1}$ of unit
quaternions:%
\begin{equation*}
\begin{array}{llll}
\rho : & G & \longrightarrow & A(H) \\
& a & \rightarrow & (0,A) \\
& b & \rightarrow & (0,B)%
\end{array}%
\end{equation*}

From now on we will assume that at least one of the elements $\rho (a),$ $%
\rho (b)$ has translational part different from $0$, and that
$(A,B)$ is a irreducible pair of conjugate unit quaternions. By
Theorem \ref{teo2} we may assume that
\begin{eqnarray*}
\rho (a) &=&(sA^{-},A) \\
\rho (b) &=&(sB^{-}+(A^{-}B^{-})^{-},B)
\end{eqnarray*}

Because $\rho $ is an homomorphism of the semidirect product $
H_{0}\rtimes U_{1}=A(H),$ we have
\begin{equation}
\rho (w(a,b))=(\frac{\partial w}{\partial a}|_{\phi }\circ v+\frac{\partial w%
}{\partial b}|_{\phi }\circ u,w(A,B))=(0,I)  \label{erelfox}
\end{equation}
where $\frac{\partial w}{\partial a}|_{\phi }$ is the Fox derivative
of the word $w(a,b)$ with respect to $a$, and evaluated by $\phi $
such that $\phi (a)=A$, $\phi (b)=B$. (See \cite{CF1963}.) The
equation (\ref{erelfox}) yields two relations between the parameters
\begin{eqnarray*}
x &=&A^{+}=B^{+} \\
y &=&-(A^{-}B^{-})^{+} \\
s &=&vector\ shift\ parameter
\end{eqnarray*}
the relations are%
\begin{equation}
w(A,B)=I  \label{erel1}
\end{equation}%
\begin{equation}
\frac{\partial w}{\partial a}|_{\phi }\circ v+\frac{\partial w}{\partial b}%
|_{\phi }\circ u=0  \label{erel2}
\end{equation}

 The relation (\ref{erel1}) yields the ideal $\mathcal{I}
_{G}^{c}=\{p_{i}(x,y)\mid i\in \left\{ 1,2,3,4\right\} \}$, as we proved in
Section \ref{srepresentationinU}. It defines $V(\mathcal{I}_{G}^{c})$ the
algebraic variety of c-representations of $G$ in $SL(2,\mathbb{C})$.

The relation (\ref{erel2}) produces four polynomials in $x,y,s$: $
\{q_{j}(x,y,s)\mid j\in \left\{ 1,2,3,4\right\} \}.$ The ideal
\begin{equation*}
\mathcal{I}_{aG}^{c}=\{p_{i}(x,y),q_{j}(x,y,s)\mid i,j\in \left\{
1,2,3,4\right\} \}
\end{equation*}
defines an algebraic variety, that we call
$V_{a}(\mathcal{I}_{aG}^{c})$ the \emph{variety of affine
c-representations of }$G$\emph{\ in }$A(H)$ up to conjugation in
$\mathcal{E}q(H_{0})$.

Let
\begin{equation*}
\begin{array}{llll}
\rho : & G & \longrightarrow & A(H) \\
& a & \rightarrow & \rho (a)=(sA^{-},A) \\
& b & \rightarrow & \rho (b)=((sB^{-}+(A^{-}B^{-})^{-},B)
\end{array}
\end{equation*}
be a representation of $G$ in the affine group of a quaternion
algebra $H$. The composition of $\rho $ with the projection $\pi
_{2}$ on the second factor of $A(H)=H_{0}\rtimes U_{1}$ gives the
linear part of $\rho $ and it is a representation $\widehat{\rho }$
on the group of unit quaternions.
\begin{equation*}
\begin{array}{llll}
\widehat{\rho }=\pi _{2}\circ \rho : & G & \longrightarrow & U_{0} \\
& a & \rightarrow & A \\
& b & \rightarrow & B
\end{array}
\end{equation*}
The composition of $\rho $ with the projection $\pi _{1}$on the first factor
of $A(H)=H_{0}\rtimes U_{1}$ gives the translational part of $\rho :$
\begin{equation*}
\begin{array}{llll}
v_{\rho }=\pi _{1}\circ \rho : & G & \longrightarrow & H_{0}%
\end{array}
\end{equation*}
which is called a \emph{cocycle} because it satisfy the \emph{cocycle
condition}
\begin{equation}
v_{\rho }(g_{1}g_{2})=v_{\rho }(g_{1})+\widehat{\rho }(g_{1})v_{\rho }(g_{2})
\label{ecocycle}
\end{equation}
Therefore $\widehat{\rho }$ corresponds to a point in the character variety
of c-representations of $G$ and it is determined by the characters $x$ and $%
y $. Reciprocally, given a representation
\begin{equation*}
\begin{array}{llll}
\widehat{\rho }: & G & \longrightarrow & U_{1}%
\end{array}
\end{equation*}
an \emph{affine deformation} $\rho $ of $\widehat{\rho }$ is a
homomorphism $
\begin{array}{llll}
\rho : & G & \longrightarrow & A(H)%
\end{array}
$ such that
\begin{eqnarray*}
(i)\quad \rho (g) &=&(v_{g},\widehat{\rho }(g)),\forall g\in G \\
(ii)\quad \quad v_{g} &\neq &0\quad \text{for some }g\in G
\end{eqnarray*}
Observe that every representation $\rho $ of $G$ in the affine group of a
quaternion algebra $H$ is an affine deformation of $\widehat{\rho }=\pi
_{2}\circ \rho $ if and only if the cocycle $v_{\rho }=\pi _{1}\circ \rho
:G\longrightarrow H_{0}$ is not constant.

We are interested in the classes of affine deformations up to
conjugation in $\mathcal{E}q(H)$. Each of these classes is
determined by the parameter $s$.

\begin{example}[The trefoil knot]
Consider the group of the trefoil knot $3_{1}$
\begin{equation*}
G(3_{1})=|a,b;aba=bab|
\end{equation*}%
Let
\begin{equation*}
\begin{array}{llll}
\rho : & G(3_{1}) & \longrightarrow  & A(H) \\
& a & \rightarrow  & \rho (a)=(sA^{-},A) \\
& b & \rightarrow  & \rho (b)=(sB^{-}+(A^{-}B^{-})^{-},B)%
\end{array}%
\end{equation*}%
be a representation of $G$ in the affine group of a quaternion
algebra $H$. Assume that $\left\{
A^{-},B^{-},(A^{-}B^{-})^{-}\right\} $ is a basis of $H_{0}$.
\end{example}

Then the parameters $x,y,s$ satisfy the equations:
\begin{eqnarray*}
2y-(2x^{2}-1) &=&0  \label{epolytrebol} \\
4x^{2}+4sx-3 &=&0  \label{epolytrebol2}
\end{eqnarray*}

Therefore the variety of affine representation of $G(3_{1})$ is
\begin{equation*}
V_{a}(\mathcal{I}_{aG(3_{1})})=<2y-(2x^{2}-1),4x^{2}+4sx-3>.
\end{equation*}

We have studied the representations corresponding to the points in
this variety $V_{a}(\mathcal{I}_{aG(3_{1})})$, in \cite{HLM2009}.

\begin{example}[The Figure Eight knot]
Consider the group of the Figure Eight knot $4_{1}$.
\begin{equation*}
G(4_{1})=\pi _{1}(S^{3}-4_{1})=\left\vert
a,b:aba^{-1}b^{-1}a=ba^{-1}b^{-1}ab\right\vert
\end{equation*}
Let
\begin{equation*}
\begin{array}{llll}
\rho : & G(4_{1}) & \longrightarrow  & A(H) \\
& a & \rightarrow  & \rho (a)=(sA^{-},A) \\
& b & \rightarrow  & \rho (b)=(sB^{-}+(A^{-}B^{-})^{-},B)%
\end{array}
\end{equation*}
be a representation of $G$ in the affine group of a quaternion
algebra $H$. Assume that $\left\{
A^{-},B^{-},(A^{-}B^{-})^{-}\right\} $ is a basis of $ H_{0}$.
\end{example}

Then, using the computer program Mathematica, we found that the
parameters $ x,y,s$ satisfy the equations:

\begin{eqnarray*}
p_{1}(x,y)&=&1-6x^{2}+4x^{4}-2y-4y^{2} =0 \\
q_{1}(x,y,s)&=&5+22sx-9x^{2}-16sx^{3}+4x^{4}+15y-12x^{2}y =0 \\
q_{2}(x,y,s)&=&-5-10sx+19x^{2}-12x^{4}-5y-16sxy+4x^{2}y =0
\end{eqnarray*}

Therefore the variety of affine representation of $G(4_{1})$ is
\begin{equation*}
V_{a}(\mathcal{I}_{aG(4_{1})})=<p_{1}(x,y),q_{1}(x,y,s),q_{2}(x,y,s)>.
\end{equation*}


\begin{thebibliography}{10}

\bibitem{BH1995}
G.~W. Brumfiel and H.~M. Hilden.
\newblock {\em {${\rm SL}(2)$} representations of finitely presented groups},
  volume 187 of {\em Contemporary Mathematics}.
\newblock American Mathematical Society, Providence, RI, 1995.

\bibitem{CDGM2003}
V.Charette, T. Drumm, W. Goldman, and M. Morrill.
\newblock Complete flat affine and {L}orentzian manifolds.
\newblock {\em Geom. Dedicata}, 97:187--198, 2003.
\newblock Special volume dedicated to the memory of Hanna Miriam Sandler
  (1960--1999).

\bibitem{CL1996}
D.~Cooper and D.~D. Long.
\newblock Remarks on the {$A$}-polynomial of a knot.
\newblock {\em J. Knot Theory Ramifications}, 5(5):609--628, 1996.

\bibitem{CF1963}
R.~H. Crowell and R.~H. Fox.
\newblock {\em Introduction to knot theory}.
\newblock Based upon lectures given at Haverford College under the Philips
  Lecture Program. Ginn and Co., Boston, Mass., 1963.

\bibitem{CS1983}
M. Culler and P.~B. Shalen.
\newblock Varieties of group representations and splittings of {$3$}-manifolds.
\newblock {\em Ann. of Math. (2)}, 117(1):109--146, 1983.

\bibitem{Rham1967}
G. de~Rham.
\newblock Introduction aux polyn\^omes d'un n\oe ud.
\newblock {\em Enseignement Math. (2)}, 13:187--194 (1968), 1967.

\bibitem{FG1983}
D. Fried and W.~M. Goldman.
\newblock Three-dimensional affine crystallographic groups.
\newblock {\em Adv. in Math.}, 47(1):1--49, 1983.

\bibitem{G1985}
W.~M. Goldman.
\newblock Nonstandard {L}orentz space forms.
\newblock {\em J. Differential Geom.}, 21(2):301--308, 1985.

\bibitem{GM1993}
F.~Gonz\'{a}lez-Acu\~{n}a and J.~M. Montesinos-Amilibia.
\newblock On the character variety of group representations in {${\rm
  SL}(2,{\bf C})$} and {${\rm PSL}(2,{\bf C})$}.
\newblock {\em Math. Z.}, 214(4):627--652, 1993.

\bibitem{HLM2005}
H.M. Hilden, M.~T. Lozano, and J.~M. Montesinos-Amilibia.
\newblock Peripheral polynomials of hyperbolic knots.
\newblock {\em Topology Appl.}, 150(1-3):267--288, 2005.

\bibitem{HLM2009}
H.M. Hilden, M.~T. Lozano, and J.~M. Montesinos-Amilibia.
\newblock On the affine representations of the trefoil knot group.
\newblock {\em preprint}, 2009.

\bibitem{HLM2003}
H.M. Hilden, M.T. Lozano, and J.~M. Montesinos-Amilibia.
\newblock Character varieties and peripheral polynomials of a class of knots.
\newblock {\em J. Knot Theory Ramifications}, 12(8):1093--1130, 2003.

\bibitem{HLM1992d}
H.~M. Hilden, M.~T. Lozano, and J.~M.
  Montesinos-Amilibia.
\newblock On the character variety of group representations of a {$2$}-bridge
  link {$p/3$} into {${\rm PSL}(2,{\bf C})$}.
\newblock {\em Bol. Soc. Mat. Mexicana (2)}, 37(1-2):241--262, 1992.
\newblock Papers in honor of Jos\'e Adem (Spanish).

\bibitem{HLM1995a}
H.~M. Hilden, M.~T. Lozano, and J.~M.
  Montesinos-Amilibia.
\newblock On the arithmetic {$2$}-bridge knots and link orbifolds and a new
  knot invariant.
\newblock {\em J. Knot Theory Ramifications}, 4(1):81--114, 1995.

\bibitem{Lam}
T.Y.Lam.
\newblock {\em The algebraic theory of quadratic forms}.
\newblock Mathematics Lecture Note Series. W. A. Benjamin, Inc., Reading,
  Mass., 1973.

\bibitem{Riley1984}
R.Riley.
\newblock {\em Nonabelian representations of $2$-bridge knot groups}.
\newblock Quart. J. Math. Oxford Ser. (2)  35,  no. 138, 191--208, 1984.

\end{thebibliography}
\end{document}